\documentclass[reqno]{amsart} 
\usepackage{amssymb}

\newtheorem{theorem}{Theorem}
\newtheorem*{conjecture}{Conjecture}
\newtheorem{corollary}[theorem]{Corollary}

 \numberwithin{equation}{section}

\newcommand{\be}{\begin{equation}}
\newcommand{\ee}{\end{equation}}
\newcommand{\ba}{\begin{eqnarray}}
\newcommand{\ea}{\end{eqnarray}}

\newcommand{\lab}[1]{\label{#1}}

\newcommand{\C}{\mathbb C}
\newcommand{\R}{\mathbb R}
\newcommand{\Z}{\mathbb Z}

\newcommand{\T}{\mathbb T}
\newcommand{\PP}{\mathcal P}

\newcommand{\m}{\omega}
\newcommand{\ve}{\varepsilon}
\newcommand{\eg}{\Gamma}
\newcommand{\e}{\text e}

\begin{document}

\title[Elliptic hypergeometric functions]
{Essays on the theory of elliptic\\  hypergeometric functions}

\author{V. P. Spiridonov}

 \address{Laboratory of Theoretical Physics, JINR, Dubna, 141980 Russia;
e-mail address: spiridon@theor.jinr.ru}

\begin{abstract}
We give a brief review of the main results of the theory of
elliptic hypergeometric functions --- a new class of special functions of
mathematical physics. We prove the most general univariate exact
integration formula generalizing Euler's beta integral, which is
called the elliptic beta integral. An elliptic analogue of the
Gauss hypergeometric function is constructed together with the
elliptic hypergeometric equation for it. Biorthogonality relations for
this function and its particular subcases are described.
We list known elliptic beta integrals on root systems and
consider symmetry transformations for the corresponding elliptic hypergeometric
functions of the higher order.
\end{abstract}

\thanks{To be published in {\em Russian Math. Surveys}}

\maketitle

\tableofcontents

\section{Introduction}

Theory of special functions is widely used in theoretical and
mathematical physics as a handbook collection of exact mathematical
formulae together with the methods of their derivation. This concerns
the series summation formulae, exactly computable integrals,
symmetry transformations for functions, differential or other
equations solvable in terms of ``simple" functions, and so on.
An impetuous buildup of such a database, which was taking place
in the XIX century and which was on a top of priorities
of the mathematics of that time, has changed in the XX century
by an essential  deterioration of interest to special functions,
investigations of which started to be considered as a
pursuit of a secondary importance. Such an attitude to investigations
in this field was justified
by an opinion that all principle types of interesting functions
``with classical properties" (elliptic, hypergeometric,
automorphic, and some other functions) have been found already,
and it remains only to investigate them in more detail.

The theory of special functions of hypergeometric type
was developing during several centuries, starting from
the fundamental results obtained by Euler \cite{aar}.
Gauss, Jacobi, Riemann, Kummer and other prominent mathematicians
were contributing to its foundations.  The Gauss hypergeomeric
function $_2F_1(a,b;c;x)$ is a canonical example of functions
of such a type. According to the approach by Pochhammer and Horn \cite{ggr},
the generalized plain hypergeometric series can be defined as the
sums $\sum_n c_n$ for which the ratio $c_{n+1}/c_n$ is a rational function of $n$.
In 1847 Heine has introduced a $q$-analogue of the $_2F_1$-series
$_2\varphi_1(a,b;c;q;x)$ \cite{gas-rah:basic}. The general series
of such type $\sum_n c_n$ are characterized by the property that
for them the ratio $c_{n+1}/c_n$ is a rational function of
$q^n,$ where $q$ is some complex parameter. Until recent time
only these two classes of functions of hypergeometric type were known
(including integral representations for them), and
many papers were devoted to their investigations.

Unexpectedly, around ten years ago it became clear that there exist
hypergeometric functions of the third type, which are related to elliptic curves.
Such objects appeared for the first time within the quantum inverse
scattering method developed for exactly solvable models of statistical
mechanics \cite{stf,tf} as elliptic solutions of the Yang--Baxter
equation \cite{bax:partition,djkmo}. As shown by Frenkel and Turaev \cite{ft},
these solutions (called elliptic $6j$-symbols) are expressed
in terms of an elliptic generalization of the terminating very-well poised
balanced $q$-hypergeometric series $_{10}\varphi_9$ with discrete
values of the parameters. A generalized $(1+1)$-dimensional
integrable chain similar to the discrete time Toda chain
was constructed in the paper \cite{spi-zhe:spectral}, and it was
shown that the same terminating series with arbitrary parameters
appears as a particular solution of the corresponding Lax pair equations.

The general formal definition of elliptic hypergeometric series,
which was suggested and investigated in detail in \cite{spi:theta1},
describes them as the series $\sum_n c_n$, for which the
ratio $c_{n+1}/c_n$ is equal to an elliptic function of $n$. Within this scheme,
the case  \cite{ft} is characterized by the presence of several interesting
structural restrictions upon the coefficients  $c_n$. There are certain
difficulties in the description of infinite series connected with their
convergency. Therefore the general elliptic hypergeometric functions
are defined by the integral representations \cite{spi:theta2}.

The Meijer function can be considered as the most general plain
hypergeometric function \cite{erd:higher}. It is defined by
a contour integral of some ratio of the Euler gamma functions.
For integral representations of more complicated functions of
the hypergeometric type one needs the generalized gamma functions
the theory of which was developed by Barnes \cite{bar:multiple}
and Jackson \cite{jac:basic} more than a century ago. The Jackson
$q$-gamma function is necessary for the description of
$q$-hypergeometric functions at $|q|<1$. A more complicated
function is needed when $q$ lies on the unit circle, $|q|=1$
\cite{shi:kronecker,kur,fad:discrete1,rui:first}.
These functions, related to the Barnes gamma function
of the second order, are actively used in the modern
mathematical physics in the description of
quantum integrable models and representations of
quantum algebras \cite{jm,tv,fkv,pon-tes:clebsch,kls:unitary,volk}.

For the definition of elliptic hypergeometric integrals one needs the
elliptic gamma function related to the Barnes gamma function of the third
order. Importance of the elliptic gamma function was stressed by Ruijsenaars
 \cite{rui:first}, who gave to it this name and considered some of its
properties.
Modular transformations of this function were described in \cite{fel-var:elliptic}.
In \cite{spi:theta2}, a modified elliptic gamma function was built
which remains also well defined in the case when one of the base parameters
lies on the unit circle. Other aspects of this function were
investigated  in \cite{FR,nar,rai:limits}.

The first exact integration formula, which uses the elliptic
gamma function, was constructed by the author in \cite{spi:umn}.
It represents the most general known univariate integral
generalizing Euler's beta integral \cite{aar}. This elliptic beta
integral serves as a basis for building general very-well poised
elliptic hypergeometric functions. The first direction of
generalizations consists in increasing the number of free parameters
entering the integrand, which leads to the elliptic analogues of
the functions $_{s+1}F_s$ \cite{spi:theta2,spi:bailey2}.
In particular, in this way one builds an elliptic analogue
of the Gauss hypergeometric function and derives its properties
 \cite{spi:thesis}. The second direction of generalizations
increases the number of integrations in such a way that the
integrands acquire symmetries in the integration variables
related to the root systems  \cite{spi:theta2,die-spi:elliptic,
die-spi:selberg,rai:trans,rai:abelian,spi-war:inversions}.
One of such generalizations leads to an elliptic analogue of the
Selberg integral.
Over a short period of time, in the papers mentioned above and
cited below in the list of references, a systematic theory
of elliptic hypergeometric functions of one and many variables
has been built. The present review is devoted to a brief description of
this theory.

\section{Generalized gamma functions}

We use the symbols: $\Z=0,\pm 1,\pm2,\ldots;$
$\C$ -- the open complex plane; 
${\C}^* =\C/\{0\}$; $\R$ -- the real axis; $i=\sqrt{-1}$.

\subsection{The Barnes multiple gamma function}\

The Euler gamma function is a cornerstone of the theory of ordinary
hypergeometric functions \cite{aar}.
Different analogues of this function with many parameters
have been considered in the mathematical literature of the
beginning of the twentieth century. The most complete investigations
of the generalized gamma functions belong to Barnes \cite{bar:multiple}.
As a starting point in his work serves a function
generalizing the Hurwitz zeta-function \cite{aar}
$$
\zeta(s,u)=\sum_{n=0}^\infty\frac{1}{(u+n)^s}, \quad \text{Re}(s)>1,
$$
which reduces to the Riemann zeta-function for  $u=1$.
Because of the grown interest the Barnes theory was
sufficiently widely discussed in the recent literature,
for instance, in \cite{jm,FR,nar,rai:limits,rui:barnes}.

Let us take $m$ quasiperiods $\omega_j\in\C$, which we assume to be linearly
independent over $\Z$ for simplicity (the condition of incommensurability).
For $s, u\in\C$ the Barnes zeta function is defined by the $m$-fold series
$$
\zeta_m(s,u;\mathbf{\omega})=\sum_{n_1,\ldots,n_m=0}^\infty
\frac{1}{(u+\Omega)^s}, \qquad\Omega=n_1\omega_1+\ldots+n_m\omega_m,
$$
converging for $\text{Re}(s)>m$ and under the condition that
all $\m_j$ lie in one half-plane defined by some line passing through
the zero point of the coordinates axes. Because of the latter requirement,
the sequences $n_1\omega_1+\ldots+n_m\omega_m$ do not have
accumulation points on the finite plane for any $n_j\to+\infty$.
It is convenient to assume for definiteness that Re$(\m_j)>0$ or Im$(\m_j)> 0$.

The function $\zeta_m(s,u;\mathbf{\omega})$ satisfies the following set
of finite difference equations
\be
\zeta_m(s,u+\omega_j;\mathbf{\omega})-\zeta_m(s,u;\mathbf{\omega})
=-\zeta_{m-1}(s,u;\mathbf{\omega}(j)), \quad j=1,\ldots,m,
\lab{zeta-eq}\ee
where $\mathbf{\omega}(j)=(\omega_1,\ldots,\omega_{j-1},\omega_{j+1},\ldots,
\omega_m)$ and $\zeta_0(s,u;\mathbf{\omega})=u^{-s}$.
It may be continued analytically (meromorphically) to the whole
complex plane $\text{Re}(s)\leq m$ with the simple poles at
the points $s=1,2,\ldots,m$. The Barnes multiple gamma function
is defined by the equality
$$
\Gamma_m(u;\mathbf{\omega})
=\exp(\partial \zeta_m(s,u;\mathbf{\omega})/\partial s)\big|_{s=0}.
$$
It has an infinite product representation of the form
\begin{equation}
\frac{1}{\Gamma_m(u;\mathbf{\omega})} =
e^{\sum_{k=0}^m\gamma_{mk}\frac{u^k}{k!}}\;u\makebox[-1em]{} \prod_{n_1,\ldots,
n_m=0}^\infty\makebox[-1,2em]{}'\makebox[0,5em]{} \left(1+\frac{u}{\Omega}\right)
e^{\sum_{k=1}^m(-1)^k\frac{u^k}{k\Omega^k}},
\label{m-gamma}\end{equation}
where $\gamma_{mk}$ are some constants analogous to the Euler constant.
The prime of the product means that the point $n_1=\ldots=n_m=0$ is
skipped in it.  In particular, the function $\Gamma_1(u;\omega)$
is directly related to the Euler gamma function $\Gamma(u)$,
$$
\Gamma_1(u;\omega)=\frac{\omega^{u/\omega}}{\sqrt{2\pi\omega}}
\Gamma\left(\frac{u}{\omega}\right),
$$
which can be checked by straightforward manipulations with the Hurwitz zeta function.
We note that Barnes used in \cite{bar:multiple} a different normalization
of the $\Gamma_m$-function, in which one has $\gamma_{m0}=0$.

The function $\Gamma_m(u;\mathbf{\omega})$ satisfies $m$ finite difference
equations of the first order, obtained by differentiation of equalities
\eqref{zeta-eq} at the point $s=0$:
\begin{equation}
\Gamma_m(u+\omega_j;\mathbf{\omega})
=\frac{1}{\Gamma_{m-1}(u;\mathbf{\omega}(j))}\, \Gamma_m(u;\mathbf{\omega}),
\qquad j=1,\ldots,m,
\label{bar-eq}\end{equation}
where  $\Gamma_0(u;\omega):=u^{-1}$.

The following integral representations for the Euler gamma function
are well known
$$
\Gamma(s)=\int_0^\infty t^{s-1}e^{-t}dt=\frac{i}{2\sin(\pi s)}
\int_{C_H}(-t)^{s-1}e^{-t}dt,
$$
where in the first case Re$(s)>0$, and in the second
expression $|\arg(-t)|<\pi$ and the Hankel contour
$C_H$ starts and finishes near the $+\infty$ point, turning around
the half-axis $[0,\infty)$ counterclockwise. One can write with their help
\begin{eqnarray*} &&
\zeta_m(s,u;\mathbf{\omega})=\frac{1}{\Gamma(s)}
\int_0^\infty \frac{t^{s-1}e^{-ut}}{\prod_{k=1}^m(1-e^{-\omega_k t)} } dt
\\ && \makebox[2em]{}
=\frac{i\Gamma(1-s)}{2\pi}
\int_{C_H} \frac{(-t)^{s-1}e^{-ut}}{\prod_{k=1}^m(1-e^{-\omega_k t}) } dt
\end{eqnarray*}
and analytically continue this function in $s$ to the whole complex plane.
Using the latter expression, Barnes has derived the following integral
representation for the multiple gamma functions
\begin{eqnarray}
\Gamma_m(u;\mathbf{\omega})=\exp\left(
\frac{1}{2\pi i}
\int_{C_H} \frac{e^{-ut}(\log(-t)+\gamma)}
{t\prod_{k=1}^m(1-e^{-\omega_k t})} dt \right),
\label{bar-gam-int}\end{eqnarray}
where $\gamma$ is the Euler constant.

The values $\zeta_m(0,u;\mathbf{\omega})$ are expressed in terms
of the multiple Bernoulli polynomials
$$
\zeta_m(0,u;\mathbf{\omega})=\frac{(-1)^m}{m!}B_{m,m}(u|\mathbf{\omega}),
$$
which are defined by the generating function
\be
\frac{x^m e^{xu}}{\prod_{k=1}^m(e^{\m_k x}-1)}
=\sum_{n=0}^\infty B_{m,n}(u|\m_1,\ldots,\m_m)\frac{x^n}{n!}.
\lab{ber}\ee
We shall need in the following the first three diagonal polynomials
\begin{eqnarray*} &&
B_{1,1}(u|\m_1)=\frac{u}{\m_1}-\frac{1}{2}, \quad
\\ &&
B_{2,2}(u|\m_1,\m_2)=\frac{1}{\m_1\m_2}\left(u^2-(\m_1+\m_2)u+\frac{\m_1^2+\m_2^2}{6}+
\frac{\m_1\m_2}{2}\right),
\\ &&
B_{3,3}(u|\m_1,\m_2,\m_3)=\frac{1}{\m_1\m_2\m_3}\Biggl(u^3-\frac{3u^2}{2}\sum_{k=1}^3\m_k
+\frac{u}{2}\left(\sum_{k=1}^3\m_k^2+3\sum_{j<k}\m_j\m_k\right)
\\ && \makebox[9em]{}
-\frac{1}{4}\left(\sum_{k=1}^3\m_k\right)\sum_{j<k}\m_j\m_k\Biggr).
\end{eqnarray*}

Theory of the plain hypergeometric functions is built with the help
of the Euler gamma function or $\Gamma_1(u;\omega_1)$; the
$q$-hypergeometric functions are tied to
$\Gamma_2(u;\omega_1,\omega_2)$, and the elliptic hypergeometric
functions ``live" at the level of the Barnes multiple gamma function of
the third order, respectively.

\subsection{The elliptic gamma functions}\

Let $\omega_{1},\omega_{2},\omega_{3}$ denote complex parameters
linearly independent over $\Z$ and lying in the right half-plane.
We define with their help the base variables $p,q,r\in\C$:
\be
q= e^{2\pi i\frac{\omega_1}{\omega_2}}, \quad
p=e^{2\pi i\frac{\omega_3}{\omega_2}}, \quad  r=e^{2\pi i\frac{\omega_3}{\omega_1}}
\lab{bases}\ee
and their modular transformed ($\tau\to -1/\tau$) partners
\be
\tilde q= e^{-2\pi i\frac{\omega_2}{\omega_1}}, \quad
\tilde p=e^{-2\pi i\frac{\omega_2}{\omega_3}},   \quad
\tilde r=e^{-2\pi i\frac{\omega_1}{\omega_3}}.
\lab{mod-bases}\ee
For $|p|,|q|<1$, the infinite products
$$
(z;q)_\infty=\prod_{j=0}^\infty (1-zq^j), \quad
(z;p,q)_\infty=\prod_{j,k=0}^\infty (1-zp^jq^k)
$$
are well defined and satisfy $q$-difference equations
\begin{eqnarray}\nonumber  &&
(qz;q)_\infty=\frac{(z;q)_\infty}{1-z}, \qquad
\\ &&
(qz;q,p)_\infty=\frac{(z;q,p)_\infty}{(z;p)_\infty}, \qquad
(pz;q,p)_\infty= \frac{(z;q,p)_\infty}{(z;q)_\infty}.
\label{eq-1}\end{eqnarray}

The shortened theta function $\theta(z;p)$
\begin{equation}
\theta(z;p):=(z;p)_\infty (pz^{-1};p)_\infty
 =\frac{1}{(p;p)_\infty} \sum_{k\in\Z} (-1)^kp^{k(k-1)/2}z^k
\label{short-theta}\end{equation}
plays a key role in our considerations. It obeys the following simple
symmetry transformations:
\begin{equation}
\theta(pz;p)=\theta(z^{-1};p)=-z^{-1}\theta(z;p),
\label{theta-trafo}\end{equation}
and has zeros, $\theta(z;p)=0$, at $z=p^k,\; k\in\Z.$
Evidently, $\theta(z;0)=1-z$. For $k>0$, we have
$$
\theta(p^kz;p)=\frac{\theta(z;p)}{(-z)^kp^{\binom{k}{2}}}, \qquad
\theta(p^{-k}z;p)=\frac{(-z)^k\theta(z;p)}{p^{\binom{k+1}{2}}}.
$$

We shall need the $\tau\to -1/\tau$ modular transformation rule for
 $\theta(z;p)$, for the description of which it is necessary to use
the exponential parameterization of variables:
\begin{equation}
\theta\left(e^{-2\pi i\frac{u}{\omega_1}};
e^{-2\pi i\frac{\omega_2}{\omega_1}}\right)
=e^{\pi iB_{2,2}(u|\omega_1,\omega_2)} \theta\left(e^{2\pi i\frac{u}{\omega_2}};
e^{2\pi i\frac{\omega_1}{\omega_2}}\right),
\label{mod-theta}\end{equation}
where $B_{2,2}(u|\omega_1,\omega_2)$ is the second Bernoulli polynomial.
In the following it is convenient to use the compact notations
$$
\theta(a_1,\ldots,a_k;p):=\theta(a_1;p)\cdots\theta(a_k;p), \quad
\theta(at^{\pm 1};p):=\theta(at;p)\theta(at^{-1};p).
$$

The simplest gamma function can be defined as a special meromorphic
solution of the functional equation $f(u+\omega_1)=uf(u)$.
Following Jackson's approach  \cite{jac:basic}, we shall be
connecting $q$-gamma functions with the meromorphic
solutions of the equation
\be
f(u+\omega_1)=(1-e^{2\pi i u/\omega_2})f(u),
\lab{q-eqn}\ee
where $q=e^{2\pi i\omega_1/\omega_2}$. Introducing the variable
$z=e^{2\pi iu/\omega_2}$, this equation can be replaced by
$f(qz)=(1-z)f(z)$. For $|q|<1$ its particular solution,
analytical at the point $z=0$, is determined by a simple iteration.
This yields the standard $q$-gamma function
$\gamma_q(z)=1/(z;q)_\infty$ (which can be considered also
as a $q$-exponential function \cite{gas-rah:basic}).
This expression differs from the Jackson $q$-gamma function
\be
\Gamma_q^{(J)}(u)=\frac{(q;q)_\infty}
{(q^{u};q)_\infty} (1-q)^{1-u},
\lab{jac}\ee
satisfying the equation $\Gamma_q^{(J)}(u+1)/\Gamma_q^{(J)}(u)=(1-q^u)/(1-q)$,
by a simple change of the argument and by a simple multiplier.
The limiting transition to the ordinary gamma function has the form
 $\lim_{q\to1}\Gamma_q^{(J)}(u)=\Gamma(u)$
\cite{aar,gas-rah:basic}, but for a simplification of $q$-hypergeometric
formulae it is more convenient to use the function $\gamma_q(z)$.

The modified $q$-gamma function, which remains well defined at
 $|q|=1$ as well, has the form
\be
\gamma(u;\omega_1,\omega_2)=
\exp\left(-\int_{\R+i0}\frac{e^{ux}}
{(1-e^{\omega_1 x})(1-e^{\omega_2 x})}\frac{dx}{x}\right),
\label{mod-q-gamma}\ee
where the contour $\R+i0$ passes along the real axis turning
over the point $x=0$ from above in an infinitesimal way. This
function appeared in the number theory \cite{shi:kronecker}
and in the theory of completely integrable systems
\cite{fad:discrete1,rui:first,jm}. It figures in the literature under
the different names: ``the double sine" \cite{kur},
``the non-compact quantum dilogarithm" \cite{fkv},
``the hyperbolic gamma function" \cite{rui:first}.

Let $\text{Re}(\omega_1), \text{Re}(\omega_2)>0$. Then integral
\eqref{mod-q-gamma} is convergent for
$0<\text{Re}(u)< \text{Re}(\omega_1+\omega_2)$.
Under appropriate restrictions  on $u$ and $\omega_{1,2}$,
integral \eqref{mod-q-gamma} can be computed as a convergent sum
of residues of the poles in the upper half-plane. For
$\text{Im}(\omega_1/\omega_2)>0$ this leads to the expression
\begin{equation}
\gamma(u;\omega_1,\omega_2)=\frac{(e^{2\pi iu/\omega_1}\tilde q;\tilde q)_\infty}
{(e^{2\pi i u/\omega_2}; q)_\infty},
\label{mod-gam-prod}\end{equation}
which is continued analytically to the whole complex plane of $u$
(this expression satisfies equation \eqref{q-eqn} in an evident way).
This  $q$-gamma function serves as a main ``brick" in the construction
of analytical $q$-hypergeometric functions at $|q|=1$,
which were not considered in the literature until the recent time.

The modified $q$-gamma function is proportional to a ratio
of two Barnes gamma functions of the second order. The general
relation of such a type is derived with the help of integral
representation \eqref{bar-gam-int} and has the form \cite{nar}
$$
\exp\left(-\int_{\R+i0}\frac{e^{ux}}
{\prod_{k=1}^m(e^{\omega_k x}-1)}\frac{dx}{x}\right)
=e^{\frac{\pi i}{m!}B_{m,m}(u|\mathbf{\omega})}
\frac{\Gamma_m(u;\mathbf{\omega})^{(-1)^m}}
{\Gamma_m(\sum_{k=1}^m\omega_k-u;\mathbf{\omega})},
$$
where Re$(\omega_k)>0$ and $0<\text{Re}(u)< \text{Re}(\sum_{k=1}^m\omega_k)$.
In \cite{nar} there were obtained also infinite product representations
of these functions  analogous to
\eqref{mod-gam-prod}, which we are not describing here.

Already in the works of Barnes it was noticed that
the Jacobi $\theta_1(u|\tau)$-function
can be decomposed as a product of four multiple gamma functions of
the second order with different arguments. The exact form of such
a relation is (see, for instance, \cite{FR})
\be
\theta(e^{2\pi i u};p)=\frac{e^{-\pi i B_{2,2}(u|1,\tau)}}
{\Gamma_2(u;1,\tau)\Gamma_2(1+\tau-u;1,\tau)
\Gamma_2(u-\tau;1,-\tau)\Gamma_2(1-u;1,-\tau)},
\lab{theta-gam2}\ee
where $p=e^{2\pi i\tau}$.
The general relation between multiple gamma functions and infinite
products of the Jackson type has the form  \cite{FR}
$$
\prod_{n_1,\ldots,n_m=0}^\infty(1-e^{2\pi i(u+\Omega)})
=\frac{e^{-\pi i\zeta_{m+1}(0,u;1,\alpha)}}
{\Gamma_{m+1}(u;1,\alpha)\Gamma_{m+1}(1-u;1,-\alpha)},
$$
where $\alpha=(\alpha_1,\ldots,\alpha_m)$,
$\Omega=n_1\alpha_1+\ldots+n_m\alpha_m$, and Im$(\alpha_j)>0$.

Following the logic of definitions of the $q$-gamma functions,
we connect the elliptic gamma functions with meromorphic solutions
of the finite difference equation
\be
f(u+\omega_1)=\theta(e^{2\pi i u/\omega_2};p)f(u),
\lab{e-gamma-eq}\ee
which passes to \eqref{q-eqn} for $p\to 0$.
Using the factorization \eqref{short-theta} and equalities \eqref{eq-1},
it is not difficult to see that the ratio
\begin{equation}
 \eg(z;p,q) = \frac{(pqz^{-1};p,q)_\infty}{(z;p,q)_\infty}
=\prod_{j,k=0}^\infty\frac{1-z^{-1}p^{j+1}q^{k+1}}{1-zp^{j}q^{k}},
\label{ell-gamma}\end{equation}
where $|p|,|q|<1$ and $z\in\C^*$, satisfies the equations
\begin{equation}
\eg(qz;p,q)=\theta(z;p)\eg(z;p,q),\quad \eg(pz;p,q)=\theta(z;q)\eg(z;p,q).
\label{e-gam-eq}\end{equation}
Thus, the function $f(u)=\eg(e^{2\pi i u/\omega_2};p,q)$ defines
a solution of  equation \eqref{e-gamma-eq} at $|q|, |p|<1$,
and it is called the (standard) elliptic gamma function \cite{rui:first}.
Because non-trivial triply periodic functions do not exist,
it can be defined  uniquely  as the meromorphic solution of the system
of three equations:
$$
f(u+\omega_2)=f(u),\qquad f(u+\omega_3)=\theta(e^{2\pi i u/\omega_2};q)f(u)
$$
and equation \eqref{e-gamma-eq} with the normalization of the solution
$f(\sum_{k=1}^3\omega_k/2)=1$. The reflection equation for this
generalized gamma function has the form  $\eg(z;p,q) \eg(pq/z;p,q)=1$.
For $p= 0$, we have $\eg(z;0,q)=\gamma_q(z)$. 

The modified elliptic gamma function, which is well defined for $|q|=1$
as well, has the form \cite{spi:theta2}
\be
G(u;\mathbf{\omega})=
\eg(e^{2\pi i \frac{u}{\omega_2}};p,q)
\Gamma(re^{-2\pi i \frac{u}{\omega_1}};\tilde q,r).
\lab{unit-e-gamma}\ee
It defines the unique solution of three equations:
$$
f(u+\omega_2) =\theta(e^{2\pi i u/\omega_1};r) f(u),
\qquad   f(u+\omega_3) =e^{-\pi iB_{2,2}(u|\mathbf{\omega})} f(u)
$$
and equation \eqref{e-gamma-eq} with the normalization of
the solution $f(\sum_{k=1}^3\omega_k/2)=1$.

It is easy to check \cite{die-spi:unit}, that the function
\be
G(u;\mathbf{\omega})
= e^{-\frac{\pi i}{3}B_{3,3}(u|\mathbf{\omega})}
\Gamma(e^{-2\pi i \frac{u}{\omega_3}};\tilde r,\tilde p),
\lab{mod-e-gamma}\ee
where  $|\tilde p|,|\tilde r|<1$,
satisfies the same three equations and the normalization
as the function \eqref{unit-e-gamma}. Therefore these functions coincide,
and their equality constitutes one of the laws of the modular
transformations for the elliptic gamma function related to the
$SL(3;\Z)$-group \cite{fel-var:elliptic}. From the
expression \eqref{mod-e-gamma} it follows, that
$G(u;\mathbf{\omega})$ is a meromorphic function of $u$ for
$\omega_1/\omega_2>0$, i.e. $|q|=1$.

Because of the antisymmetry condition
$B_{3,3}(\sum_{k=1}^3\omega_k-u|\mathbf{\omega})
=-B_{3,3}(u|\mathbf{\omega})$,
the reflection formula for the $G$-function has the form
$G(a;{\bf \omega})G(b;{\bf \omega})=1,$ $a+b=\sum_{k=1}^3\omega_k.$
For $|q|<1$ in the limit $p,r\to0$ (i.e., Im$(\omega_3/\omega_1)$,
Im$(\omega_3/\omega_2)\to +\infty$) expression \eqref{unit-e-gamma}
passes in an evident way to the modified $q$-gamma function
$\gamma(u;\omega_1,\omega_2)$. The representation \eqref{mod-e-gamma}
provides an alternative way for the reduction
to this function (such a limiting transition was
rigorously justified in a different way in \cite{rui:first}).
As follows from the results of \cite{rai:limits}, for a fixed
domain of values of parameters the function $G(u;\mathbf{\omega})$
converges in this limit to $\gamma(u;\omega_1,\omega_2)$
exponentially fast and
uniformly on compact subsets of this domain. This result is
important for a rigorous justification of the corresponding
degeneration of the elliptic hypergeometric integrals.

Using the theta function factorization \eqref{theta-gam2}, one can consider
equation \eqref{e-gamma-eq} as a composition of four equations for
$\Gamma_3(u;\mathbf{\omega})$ with different arguments and quasiperiods.
This permits to represent the elliptic gamma function as a ratio
of four Barnes gamma functions of the third order \cite{FR}
$$
\Gamma(e^{2\pi i u};e^{2\pi i\tau},e^{2\pi i\sigma})
=\frac{e^{-\frac{\pi i}{3}B_{3,3}(u|1,\sigma,\tau)}\Gamma_3(u;1,\sigma,\tau)
\Gamma_3(1-u;1,-\sigma,-\tau)}
{\Gamma_3(1+\sigma+\tau-u;1,\sigma,\tau)
\Gamma_3(u-\sigma-\tau;1,-\sigma,-\tau)}.
$$
For $0<\text{Im}(u)<\text{Im}(\tau+\sigma)$, one has the representation
$$
\Gamma(e^{2\pi i u};e^{2\pi i\tau},e^{2\pi i\sigma})
=\exp\left( -\frac{i}{2}\sum_{k=1}^\infty\frac{\sin(\pi k (2u-\tau-\sigma))}
{k\sin(\pi k \tau)\sin(\pi k \sigma)}
\right),
$$
through which the elliptic gamma function appeared implicitly in the
work of Baxter on the eight-vertex model \cite{bax:partition}
(see also \cite{tf,fel-var:elliptic}).

\section{The elliptic beta integral}

As a first example of  elliptic hypergeometric functions
we describe the elliptic beta integral, which was discovered
by the author in \cite{spi:umn}.

\begin{theorem}
We consider six complex parameters $t_j,\; j=1,\ldots,6$,
and two base variables $p$ and $q$ satisfying the constraints
$|p|, |q|, |t_j|<1$ and $\prod_{j=1}^6 t_j=pq$ (the balancing condition).
Then the following equality is true
\be
\kappa\int_\T\frac{\prod_{j=1}^6
\eg(t_jz^{{\pm 1}};p,q)}{\eg(z^{\pm 2};p,q)}\frac{dz}{z}
=\prod_{1\leq j<k\leq6}\eg(t_jt_k;p,q),
\label{ell-int}\ee
where $\T$ denotes the unit circle with positive orientation and
$$
\kappa=\frac{(p;p)_\infty(q;q)_\infty}{4\pi i}.
$$
\end{theorem}

Here and below we use the compact notation
\begin{eqnarray*}
&& \eg(a_1,\ldots,a_k;p,q):=\eg(a_1;p,q)\cdots \eg(a_k;p,q),\quad
\\ &&
\Gamma(tz^{\pm 1};p,q):=\Gamma(tz;p,q)\Gamma(tz^{-1};p,q),\quad
\Gamma(z^{\pm2};p,q):=\Gamma(z^2;p,q)\Gamma(z^{-2};p,q).
\end{eqnarray*}
\begin{proof}
We take variables $z,q,p\in\C$, $|q|,|p|<1$, and five complex
parameters $t_m, m=1, \dots,5,$ and compose the function
\begin{equation}
\rho(z,t_1,\ldots,t_5)=\frac{\prod_{m=1}^5\eg(t_mz^{\pm 1},At_m^{-1};p,q)}
{\eg(z^{\pm 2},Az^{\pm 1};p,q)\prod_{1\leq m<s\leq 5}\eg(t_mt_s;p,q)},
\label{kernel}\end{equation}
where $A=\prod_{m=1}^5t_m$. This function has sequences of
poles converging to zero along the points
$$
\PP=\{t_mq^ap^b,\, A^{-1}q^{a+1}p^{b+1}
\}_{m=1,\ldots, 5,\, a,b=0,1,\ldots}
$$
and diverging to infinity along their $z\to 1/z$ reciprocals $\PP^{-1}$.
Let $C$ denote a contour on the complex plane with positive orientation,
which separates the sets $\PP$ and $\PP^{-1}$ (the existence
of such a contour is the only restriction on the parameters
$t_m$). For instance, for $|t_m|<1,\; |pq|<|A|,$
the contour $C$ can coincide with the unit circle $\mathbb{T}$.
Let us prove now that
\begin{equation}\label{ell-int'}
\int_C\rho(z,t_1,\ldots,t_5)\frac{dz}{z}=\frac{4\pi i}
{(q;q)_\infty(p;p)_\infty},
\end{equation}
where from the needed formula will follow after the substitution
$A=pq/t_6$.

The first step consists in the derivation of the following
$q$-difference equation for the kernel:
\begin{equation}
\rho(z,qt_1,t_2,\ldots,t_5)-\rho(z,t_1,\ldots,t_5)
=g(q^{-1}z,t_1,\ldots,t_5)-g(z,t_1,\ldots,t_5),
\label{eqn}\end{equation}
where
\begin{equation}
g(z,t_1,\ldots,t_5)=\rho(z,t_1,\ldots,t_5)
\frac{\prod_{m=1}^5\theta(t_mz;p)}{\prod_{m=2}^5\theta(t_1t_m;p)}
\frac{\theta(t_1A;p)}{\theta(z^2,Az;p)}\frac{t_1}{z}.
\label{g-ker}\end{equation}
After the division of equation (\ref{eqn}) by $\rho(z,t_1,\ldots,t_5)$,
it takes the form
\begin{eqnarray}\nonumber
\lefteqn{\frac{\theta(t_1z,t_1z^{-1};p)}{\theta(Az,Az^{-1};p)}
\prod_{m=2}^5\frac{\theta(At_m^{-1};p)}{\theta(t_1t_m;p)}-1 } &&
\\ && \makebox[-1em]{}
=\frac{t_1\theta(t_1A;p)}{z\theta(z^2;p)\prod_{m=2}^5\theta(t_1t_m;p)}
\left(\frac{z^4\prod_{m=1}^5\theta(t_mz^{-1};p)}{\theta(Az^{-1};p)} -
\frac{\prod_{m=1}^5\theta(t_mz;p)}{\theta(Az;p)}\right).
\label{eqn-exp}\end{eqnarray}

Both sides of this equality define elliptic functions of $\log z$
(that is they are invariant under the transformation $z\to pz$)
with equal sets of poles and their residues. For example,
$$
\lim_{z\to A}\theta(Az^{-1};p)\, \Big({\text{left-hand}\atop \text{side}}\Big) =
\frac{\theta(t_1A,t_1A^{-1};p)}{\theta(A^2;p)}
\prod_{m=2}^5\frac{\theta(At_m^{-1};p)}{\theta(t_1t_m;p)},
$$
with the same result for the right-hand side. Therefore the difference
of expressions in two sides of equality (\ref{eqn-exp})
defines an elliptic function without poles, that is a constant.
This constant is equal to zero because equation (\ref{eqn-exp})
is checked in a trivial way for the choice $z=t_1$.

We integrate now (\ref{eqn}) over the variable $z\in C$ and obtain
\begin{equation}
I(qt_1,t_2,\ldots,t_5)-I(t_1,\ldots,t_5)=
\left(\int_{q^{-1}C}-\int_C\right)g(z,t_1,\ldots,t_5)
\frac{dz}{z},
\label{int-eqn}\end{equation}
where $I(t_1,\ldots,t_5)=\int_C\rho(z,t_1,\ldots,t_5)dz/z$, and
$q^{-1}C$ denotes the contour $C$ dilated with respect to the point
$z=0$ by the factor $q^{-1}$. Function (\ref{g-ker}) has sequences of poles
converging to zero along the points $z=\{t_mq^ap^b, A^{-1}q^ap^{b+1}\}$
and diverging to infinity at $z=\{ t_m^{-1}q^{-a-1}p^{-b},$
$Aq^{-a-1}p^{-b-1}\}$ for $m=1,\ldots,5$ and  $a,b=0,1,\ldots$.
For the choice $C=\mathbb{T}$, it is seen that at $|t_m|<1$ and $|p|<|A|$
there are no poles in the annulus  $1\leq |z|\leq |q|^{-1}$.
Therefore we can deform  $q^{-1}\mathbb{T}$ to $\mathbb{T}$ in
(\ref{int-eqn}) and obtain zero on the right-hand side.
As a result,  $I(qt_1,t_2,\ldots,t_5)=I(t_1,\ldots,t_5)$.

Requiring $|p|,|q|<|A|$, we have by symmetry in $p$ and $q$ that
$I(pt_1,t_2,\ldots,t_5)=$ $I(t_1,\ldots,t_5)$. Further transformations
$t_1\to q^{\pm1}t_1$ and $t_1\to p^{\pm1}t_1$ can be performed only if
they do not take parameters outside of the annulus of analyticity
of the function $I(t_1,\ldots,t_5)$.

Let us suppose temporarily that $p$ and $q$ are real, $p<q$, and
$p^n\neq q^k$ for any $n,k=0,1,\ldots$. Impose also the constraint that
the arguments of $t_m^{\pm1},\, m=1,\ldots,5,$ and $A^{\pm 1}$
differ from each other. Let $C$ denotes now a contour encircling
$\PP$ and two cuts $c_1=[t_1, t_1p^2]$, $c_2=[(pq/A)p^{-2},pq/A]$
and excluding their $z\to 1/z$ reciprocals. Then we can make
transformations $t_1\to t_1q^k$, $k=1,2,\ldots,$ until the moment
when $t_1q^k$ enters the interval $[t_1p,t_1p^2]$, after which we
replace $t_1\to t_1p^{-1}$; this does not take out
needed parameters outside of the intervals $c_1$ or $c_2$.
In this way we obtain
$I(q^jp^{-k}t_1,t_2,\ldots,t_5)=$ $I(t_1,\ldots,t_5)$ for all $j,k=0,1,\ldots$
such that $q^jp^{-k}\in [1,p]$. Since the set of such points
is dense, we come to the conclusion that $I$ does not
depend on $t_1$ and, by symmetry, on any $t_m$.

Alternatively, we can use the $p$-expansion $I(t_1,\ldots,t_5)=$
$\sum_{n=0}^\infty I_n(t_1,\ldots,t_5)p^n$
and check validity of the equalities $I_n(qt_1,\ldots,t_5)=I_n(t_1,\ldots,t_5)$
termwise. The coefficients $I_n$ are analytical in parameters
near the points $t_m=0$ (the constraints on the absolute
values of parameters from below appear from the requirement of
convergency of the $p$-expansion). Therefore, we can simply iterate
the dilations $t_1\to qt_1$ until reaching the limiting point.
As a result, $I_n$ and the integral $I$ itself do not depend on $t_1$
and, consequently, on all $t_m$.

We conclude thus that the integral $I$ is a constant depending only on
$p$ and $q$. In order to find its value, which is given by the
right-hand side of (\ref{ell-int'}), it is sufficient to consider the limit
in parameters $t_1t_2\to 1$. In this case two poles
approach the contour of integration, and it is necessary to deform this
contour so that it crosses over these poles. Then it appears that
in the limit  $t_1t_2\to 1$ only the residues of these poles
have finite values and the integral itself vanishes. (This procedure
is described in more detail below.)
After proving the integration formula in a restricted domain
of parameter values, it can be continued analytically to the domain
permissible by the contour of integration $C$.
\end{proof}

There are many ways to degenerate the elliptic beta integral.
In the simplest case it is necessary to substitute  $t_6=pq/t_1\ldots t_5$,
use the reflection formula for $\Gamma(z;p,q)$, and take the limit
$p\to 0$. After this, the elliptic beta integral degenerates to the
``trigonometric" $q$-beta-integral of Rahman \cite{rah:integral}
(which is connected with the integral representation for a $_8\varphi_7$-series
\cite{NR85}):
\begin{eqnarray*}
&& \frac{(q;q)_\infty}{4\pi i}\int_{{\mathbb T}}
\frac{(z^2;q)_\infty(z^{-2};q)_\infty(Az;q)_\infty(Az^{-1};q)_\infty}
{\prod_{m=1}^5(t_mz;q)_\infty(t_mz^{-1};q)_\infty}
\frac{dz}{z}
\\ && \makebox[4em]{}
=\frac{\prod_{m=1}^5(At_m^{-1};q)_\infty}{\prod_{1\le m<s\le 5}
(t_mt_s;q)_\infty},
\end{eqnarray*}
where $A=\prod_{m=1}^5t_m$, $|t_m|<1$. Further simplification of this equality
by taking the limit $t_5\to 0$ leads to the famous Askey--Wilson integral
$$
\frac{(q;q)_\infty}{4\pi i}\int_{{\mathbb T}}
\frac{(z^2;q)_\infty(z^{-2};q)_\infty}
{\prod_{m=1}^4(t_mz;q)_\infty(t_mz^{-1};q)_\infty}
\frac{dz}{z}
=\frac{(t_1t_2t_3t_4;q)_\infty}{\prod_{1\le m<s\le 4}(t_mt_s;q)_\infty},
$$
which serves as a measure for the Askey--Wilson polynomials \cite{aw}
-- the most general orthogonal polynomials
obeying the classical properties.
The first proof of formula \eqref{ell-int} was based on an elliptic
extension of the Askey approach \cite{ask:beta}
to computation of the Rahman integral.
Here we presented the proof obtained in the paper \cite{spi:short},
which generalizes the method of computation of the Askey--Wilson integral
from \cite{wz:invent}.

If we express in the given $q$-beta integrals infinite products $(a;q)_\infty$
in term of the Jackson $q$-gamma function and pass to the limit
$q\to 1$, then we obtain ``rational" beta integrals over non-compact
contours containing the ordinary Euler gamma function \cite{aar}.
Their further simplification by special choices of parameters leads
to the classical Euler beta-integral
$$
B(\alpha,\beta)=\int_0^1t^{\alpha-1}(1-t)^{\beta-1}dt=\frac{\Gamma(\alpha)\Gamma(\beta)}
{\Gamma(\alpha+\beta)},\quad
\text{Re}(\alpha),\text{Re}(\beta)>0.
$$
The elliptic beta integral \eqref{ell-int} gives thus
the most general (from known ones) exact univariate integration formula
including into itself the Euler beta integral
as a particular case.

\begin{corollary}
The following Frenkel and Turaev summation formula \cite{ft} is true:
\begin{equation}
\sum_{n=0}^N\frac{\theta(t_5^2q^{2n};p)}{\theta(t_5^2;p)}
\prod_{m=0}^5 \frac{\theta(t_mt_5)_n} {\theta(qt_m^{-1}t_5)_n}\, q^n=
\frac{\theta(qt_5^2,\frac{q}{t_1t_2},
\frac{q}{t_1t_3},\frac{q}{t_2t_3})_N }
{\theta(\frac{q}{t_1t_2t_3t_5},\frac{qt_5}{t_1},
\frac{qt_5}{t_2},\frac{qt_5}{t_3})_N},
\label{ft-sum}\end{equation}
where $t_4t_5=q^{-N},$ $\prod_{m=0}^5t_m=q$, and the compact notation
$$
\theta(t_1,\ldots,t_k)_n:=\prod_{j=1}^k\theta(t_j)_n
$$
is used for products of the elliptic Pochhammer symbols
$$
\theta(t)_n=\prod_{j=0}^{n-1}\theta(tq^j;p)=\frac{\eg(tq^n;p,q)}{\eg(t;p,q)}.
$$
\end{corollary}

\begin{proof}
We replace in integral \eqref{ell-int} $\T$ by a contour $C$
separating sequences of the poles $z=t_jq^ap^b,\, j=1,\ldots,6,
\, a,b=0,1,\ldots$ converging to zero from their $z\to 1/z$ reciprocals
going to infinity. This permits us to remove the constraints $|t_j|<1$
without changing the right-hand side of \eqref{ell-int}.
Substitute now $t_6=pq/A$, $A=\prod_{m=1}^5t_m$, and suppose that
$|t_m|<1,\, m=1,\ldots ,4,$ $|pt_5|<1<|t_5|$, $|pq|<|A|$, and that the
arguments of $t_m,\, m=1,\ldots,5,$ and $p,q$ are linearly independent
over $\Z$. Then the following equality is true \cite{die-spi:elliptic}:
\be
\kappa\int_C \Delta_E(z,\underline{t})\frac{d z}{z} =
\kappa\int_\mathbb{T} \Delta_E(z,\underline{t})\frac{d z}{z}
+c_0(\underline{t}) \sum_{|t_5q^n|>1,\, n \geq 0}
\nu_n(\underline{t}),
\lab{res}\ee
where $\Delta_E(z,\underline{t})=\prod_{m=1}^5\eg(t_mz^{\pm 1};p,q)
/\eg(z^{\pm2},Az^{\pm 1};p,q)$ and
\begin{eqnarray*} \nonumber
c_0(\underline{t}) =
\frac{\prod_{m=1}^4\eg(t_mt_5^{\pm 1};p,q)}
{\eg(t_5^{-2},A t_5^{\pm 1};p,q)}, \qquad
\nu_n(\underline{t}) =
 \frac{\theta(t_5^2q^{2n};p)}{\theta(t_5^2;p)}
\prod_{m=0}^5 \frac{\theta(t_mt_5)_n} {\theta(qt_m^{-1}t_5)_n}\, q^n.
\end{eqnarray*}
Here we introduced the new parameter  $t_0$ with the help of the
relation $\prod_{m=0}^5t_m $ $=q$.
The multiplier $\kappa$ is absent in the coefficient $c_0$ because of
the relation $\lim_{z\to 1}(1-z)\eg(z;p,q)=1/(p;p)_\infty(q;q)_\infty$
and due to doubling of the number of residues (the latter follows from the
symmetry of the kernel $z\to z^{-1}$).

In the limit $t_5t_4\to q^{-N},\, N=0,1,\ldots$,
the integral on the left-hand side of \eqref{res}
(coinciding with \eqref{ell-int}) and the multiplier
$c_0(\underline{t})$ in front of the sum of residues
in the right-hand side diverge. But the integral over the unit circle
$\mathbb{T}$ on the right-hand side remains finite.
After dividing all the terms by $c_0(\underline{t})$ and passing
to the limiting equality, we obtain the summation formula
\eqref{ft-sum}, which was obtained for the first time in \cite{ft}
by a completely different method.
\end{proof}

Other proofs of formula \eqref{ft-sum} are given in
\cite{gas-rah:basic,ros:proof,spi-zhe:rims,e-comb,e-taylor,chu}.
For $p\to 0$ and fixed parameters, equality  \eqref{ft-sum} is reduced
to the Jackson sum for the terminating very-well-poised balanced
$q$-hypergeometric series $_8\varphi_7$ \cite{aar}. The left-hand side
of formula  \eqref{ft-sum} represents thus an elliptic analogue
of this $q$-series.

Using the modified elliptic gamma function, it is not difficult to
construct the modified elliptic beta integral \cite{die-spi:unit},
one of the base variables
for which can lie on the unit circle, say,  $|q|=1$.

\begin{theorem}
Let Im$(\omega_1/\omega_2)\geq 0$, Im$(\omega_3/\omega_1)>0$,
Im$(\omega_3/\omega_2)>0$, and six parameters $g_j\in\C$, $j=1,\ldots,6$,
satisfy the restrictions Im$(g_j/\omega_3)<0$ and
$\sum_{j=1}^6g_j=\sum_{k=1}^3\omega_k.$ Then
\be
\tilde\kappa\int_{-\omega_3/2}^{\omega_3/2} \frac{\prod_{j=1}^6 G(g_j\pm
u;{\bf \omega})} {G(\pm 2u;{\bf \omega})} \frac{du}{\omega_2}
= \prod_{1\leq j<m\leq 6}G(g_j+g_m;{\bf \omega}),
\lab{circle-int}\ee
where
$$
\tilde\kappa= -\frac
{(q;q)_\infty(p;p)_\infty(r;r)_\infty}{2(\tilde q;\tilde q)_\infty}.
$$
Here the integration goes along the cut with the end points
$-\omega_3/2$ and $\omega_3/2$ and the convention
 $G(a\pm b;{\bf \omega})\equiv G(a+b;{\bf \omega})G(a-b;{\bf \omega})$ is used.
\end{theorem}
\begin{proof}
We substitute relation (\ref{mod-e-gamma}) into the left-hand side
of (\ref{circle-int}) and obtain
\begin{equation}
\tilde \kappa e^{\pi ia/3}
\int_{-\omega_3/2}^{\omega_3/2}\frac{\prod_{j=1}^6
\Gamma(e^{-2\pi i\frac{g_j\pm u}{\omega_3}};\tilde r,\tilde p)}
{\Gamma(e^{\mp 4\pi i\frac{u}{\omega_3}};\tilde r,\tilde p)}
\frac{du}{\omega_2}.
\label{cir-int}\end{equation}
where $a=2B_{3,3}(0|\omega)-2\sum_{j=1}^6B_{3,3}(g_j|\omega)$.
Taken restrictions on the parameters permit us to use formula
(\ref{ell-int}) with the substitutions
\begin{equation*}
z\to e^{\frac{2\pi i}{\omega_3}u}, \quad  t_j\to e^{-\frac{2\pi i
}{\omega_3}g_j},\quad p\to e^{-2\pi i\frac{\omega_1}{\omega_3}},\quad
q\to e^{-2\pi i\frac{\omega_2}{\omega_3}} ,
\end{equation*}
which yields for (\ref{cir-int})
\begin{eqnarray*}
&& \frac{2\tilde\kappa \omega_3\omega_2^{-1} e^{\pi ia/3}}{(\tilde r; \tilde r)_\infty
(\tilde p;\tilde p)_\infty}\prod_{1\leq j<m\leq 6}
\Gamma(e^{-2\pi i\frac{g_j+g_m}{\omega_3}};\tilde r,\tilde p)
 \\ && \makebox[2em]{}
= \frac{2\tilde\kappa \omega_3\omega_2^{-1} e^{\pi i(a+b)/3}}{(\tilde r; \tilde r)_\infty
(\tilde p;\tilde p)_\infty}\prod_{1\leq j<m\leq 6}
G(g_j+g_m;\boldsymbol{\omega}), \nonumber
\end{eqnarray*}
where $b=\sum_{1\leq j<m\leq 6}B_{3,3}(g_j+g_m|\omega)$.
A straightforward computation shows that
$$
a+b=\frac{1}{4}\left(\sum_{k=1}^3\omega_k\right)\left(\sum_{k=1}^3\omega_k^{-1}\right).
$$
Therefore for the choice
$$
\tilde\kappa^{-1}= \frac{2\omega_3e^{\frac{\pi i}{12}(\sum_{k=1}^3\omega_k)
(\sum_{k=1}^3\omega_k^{-1}) }}
{\omega_2(\tilde r; \tilde r)_\infty(\tilde p;\tilde p)_\infty}
$$
we obtain the needed result. After application of the modular
transformation law for the Dedekind function
\begin{equation}
e^{-\frac{\pi i}{12\tau}} \left(e^{-2\pi i/\tau};e^{-2\pi
i/\tau}\right)_\infty =(-i\tau)^{1/2}e^{\frac{\pi i\tau}{12}}
\left(e^{2\pi i\tau};e^{2\pi i\tau}\right)_\infty
\label{ded}\end{equation}
to infinite products entering definition of $\tilde\kappa$, we obtain
$$
\tilde\kappa^{-1}= -2\sqrt{\frac{\omega_1}{i\omega_2}}
\frac{e^{\frac{\pi i}{12}(\frac{\omega_1}{\omega_2}
+\frac{\omega_2}{\omega_1}) }}{(r;r)_\infty (p;p)_\infty}.
$$
One more application of the relation \eqref{ded}
permits us to replace the exponential function
by a ratio of infinite products, and this leads to the needed
form of $\kappa$.
\end{proof}

If we take the limit Im$(\omega_3)\to\infty$ in such a way that $p,r\to 0$,
then the modified elliptic beta integral reduces to a $q$-beta integral
of the Mellin--Barnes type. More precisely, for
$\omega_{1,2}$ such that
Im$(\omega_1/\omega_2)\geq 0$ and Re$(\omega_1/\omega_2)> 0$,
we substitute $g_6=\sum_{k=1}^3\omega_k
-\mathcal{A}$, where $\mathcal{A}=\sum_{j=1}^5g_j$ and apply the
inversion formula for $G(u;{\bf \omega})$.
After that we set $\omega_3=it\omega_2$, $t\to +\infty$, and obtain formally
\be
\int_{\mathbb{L}}\frac{\prod_{j=1}^5 \gamma(g_j\pm
u;{\bf \omega})}{\gamma(\pm 2u, \mathcal{A}\pm
u;{\bf\omega})} \frac{du}{\omega_2}= -2\frac
{(\tilde q;\tilde q)_\infty}{(q;q)_\infty}
\frac{ \prod_{1\leq j<m\leq 5}\gamma(g_j+g_m;{\bf \omega})}
{\prod_{j=1}^5 \gamma(\mathcal{A}-g_j;{\bf \omega})},
\label{unit-rah}\ee
where $\gamma(u;\omega_1,\omega_2)$ denotes the modified $q$-gamma
function, and the integration is taken along the line
$\mathbb{L}\equiv i\omega_2\mathbb{R}$.
This result is true provided the parameters satisfy
the constraints Re$(g_j/\omega_2)>0$ and
Re$((\mathcal{A}-\omega_1)/\omega_2)<1$. This integration formula represents
a ``hyperbolic" analogue of the Rahman integral; it was proved
for the first time by Stokman in \cite{sto:hyperbolic}.
Because of the non-compactness of the integration contour,
the described method of derivation of  \eqref{unit-rah} is rigorous
under the condition of uniform convergence of the function $G(u;{\bf \omega})$
to $\gamma(u;\omega_1,\omega_2)$, which follows from the results
obtained by Rains in \cite{rai:limits}. One can establish also formula
\eqref{unit-rah} by the method, which was used above for proving
the elliptic beta integral. The limit $q\to1$ leads to the same rational
beta integral as the ``trigonometric" Rahman integral.

Summarizing the consideration of the present section, we see that the
elliptic beta integral includes into itself the whole hierarchy of
exactly computable integrals: two types of the $q$-beta integrals,
the rational class of beta integrals, whose kernels are expressed in
terms of the Euler gamma function, and the classical Euler beta integral.
This scheme reflects the general picture of degenerations of the
elliptic hypergeometric functions which was rigorously considered
in  \cite{rai:limits}.

\section{General elliptic hypergeometric series and integrals}

In the papers \cite{spi:theta1} and \cite{spi:theta2}, the author
has proposed definitions of general elliptic hypergeometric series
and integrals, which will be considered in this section.

\subsection{An elliptic analogue of the Meijer function}\

Univariate contour integrals $\int_C\Delta(u)du$
are called the elliptic hypergeometric integrals, if the
meromorphic function $\Delta(u)$ satisfies the following system of
three equations
\be
\Delta(u+\omega_k)=h_k(u)\Delta(u),\quad k=1,2,3,
\lab{ell-gip}\ee
where $\omega_{1,2,3}\in \C$ are linearly independent over $\Z$
parameters, and $h_k(u)$ are some elliptic functions with the periods
$\omega_k,$ $\omega_{k+1}$ (we set $\omega_{k+3}=\omega_k$).

The general elliptic function of the order $s$ with the periods
$\omega_{2}$ and $\omega_3$ can be represented in the form (see the Appendix)
$$
h_1(u)=y_1\prod_{j=1}^s\frac{\theta(t_je^{2\pi iu/\omega_2};p)}
{\theta(w_je^{2\pi iu/\omega_2};p)},
$$
where $y_1$ is an arbitrary constant and $t_j, w_j$ are some parameters
satisfying the balancing condition $\prod_{j=1}^st_j=\prod_{j=1}^s w_j$
(we remind that $p=e^{2\pi i \omega_3/\omega_2}$). Using properties
of the function $\Gamma(z;p,q)$, it is not difficult to build the general
solution of the $k=1$ equation in \eqref{ell-gip} for $|q|<1$:
$$
\Delta(u)=y_1^{u/\omega_1}\varphi(u)\prod_{j=1}^s\frac{\Gamma(t_jz;p,q)}
{\Gamma(w_jz;p,q)},\qquad z=e^{2\pi iu/\omega_2},
$$
where $\varphi(u+\omega_1)=\varphi(u)$ is an arbitrary periodic
function. So, if we would restrict ourselves to a single
equation for $\Delta(u)$, then our definition of the integrals
would be highly non-unique.

An arbitrary elliptic function of the order $\ell$ with the periods $\omega_1$
and $\omega_3$ has the form
$$
h_2(u)=y_2\prod_{j=1}^{\ell }\frac{\theta(\tilde t_je^{-2\pi iu/\omega_1};r)}
{\theta(\tilde w_je^{-2\pi iu/\omega_1};r)},
$$
where $|r|<1$, $y_2$ is an arbitrary constant,  and the
parameters satisfy the constraint
$\prod_{j=1}^{\ell }\tilde t_j=\prod_{j=1}^{\ell }\tilde w_j$.
The $k=2$ equation from \eqref{ell-gip} serves now as a constraint
for the function $\varphi(u)$. It is easy to prove that
the common solution of these two equations for $|q|<1$ has the form
\begin{equation}\label{delta-ell}
\Delta(u)=\phi(u)\prod_{j=1}^s\frac{\eg(t_je^{2\pi i u/\omega_2 };p,q)}
{\eg(w_je^{2\pi i u/\omega_2 };p,q)}
\prod_{j=1}^{\ell }\frac{\Gamma(\tilde t_je^{-2\pi i u/\omega_1 };
\tilde q,r)} {\Gamma(\tilde w_je^{-2\pi i u/\omega_1 }; \tilde q,r)},
\end{equation}
where $\phi(u)$ is an arbitrary function satisfying the equations
$ \phi(u+\omega_1)=y_1\phi(u)$ and $\phi(u+\omega_2)=y_2\phi(u)$.
$\phi(u)$ is thus a meromorphic theta function with the special
quasiperiodicity multipliers the general form of which is
easily established (see the Appendix):
\begin{eqnarray*}
&& \phi(u)= e^{cu+d}\prod_{k=1}^m \frac{\theta(a_ke^{2\pi i u/\omega_2 };q)}
{\theta(b_ke^{2\pi i u/\omega_2 };q)}
\\ && \makebox[2em]{}
=e^{cu+d}\prod_{k=1}^{m}\frac{\Gamma(pa_ke^{2\pi i u/\omega_2 },
b_ke^{2\pi i u/\omega_2 };p,q)}{\Gamma(a_ke^{2\pi i u/\omega_2 },
pb_ke^{2\pi i u/\omega_2 };p,q)},
\end{eqnarray*}
where $m$ is an arbitrary integer, the parameter
$d$ is arbitrary, and the parameters $a_k, b_k, c$
are connected with $y_1$ and $y_2$ by the relations
$y_2=e^{c\omega_2}$ and $y_1=e^{c\omega_1} \prod_{k=1}^mb_ka_k^{-1}$.

Due to the representation in terms of the elliptic gamma functions,
the function $\phi(u)$ can be reduced to the pure exponential factor
by the replacement of $s$ in $\Delta(u)$ by $s+2m$ and the choice
of parameters $t_{s+k}=pa_k, t_{s+m+k}=b_k, w_{s+k}=a_k, w_{s+m+k}=pb_k$,
$k=1,\ldots,m$, which does not violate the balancing condition
$\prod_{j=1}^{s+2m}t_j=\prod_{j=1}^{s+2m}w_j$.
Since $s$ was arbitrary from the very beginning,
we can set $\phi(u)=e^{cu+d}$ without loss of generality.

As a result, two equations determine already the kernel $\Delta(u)$,
i.e.  for the elliptic function $h_3(u)$ we automatically obtain
$$
h_3(u)=e^{c\omega_3}\prod_{j=1}^s\frac{\theta(t_je^{2\pi iu/\omega_2};q)}
{\theta(w_je^{2\pi iu/\omega_2};q)}
\prod_{j=1}^{\ell }\frac{\theta(r^{-1}\tilde w_je^{-2\pi iu/\omega_1};\tilde q)}
{\theta(r^{-1}\tilde t_je^{-2\pi iu/\omega_1};\tilde q)}.
$$

To summarize, for $|q|<1$ the most general elliptic hypergeometric
integral has the form  \cite{spi:theta2,spi:thesis}
\begin{equation}\label{e-meijer}
\int_C e^{cu+d} \prod_{j=1}^s\frac{\eg(t_je^{2\pi i u/\omega_2 };p,q)}
{\eg(w_je^{2\pi i u/\omega_2 };p,q)}
\prod_{j=1}^{\ell }\frac{\Gamma(\tilde t_je^{-2\pi i u/\omega_1 };
\tilde q,r)} {\Gamma(\tilde w_je^{-2\pi i u/\omega_1 }; \tilde q,r)}\, du
\end{equation}
with two balancing conditions for the parameters indicated above
and some integration contour $C$.

Consider the definition of integrals for $|q|=1$.
It appears that now even the function $h_2(u)$
cannot be arbitrary. In this case it is necessary to take
$\ell =s$ and choose the parameters  $t_j, \tilde t_j$,
and $w_j,\tilde w_j$ in such a way that all $\Gamma$-functions
are combined to the modified elliptic gamma functions $G(u;\mathbf{\omega})$
(it is in this way that this function was built in
\cite{spi:theta2}). This leads to the integrals of the form
\begin{equation}\label{delta-ell-unit}
\int_'  e^{cu+d} \prod_{j=1}^s\frac{G(u+g_j;\mathbf{\omega})}
{G(u+v_j;\mathbf{\omega})}\, du,
\end{equation}
where the parameters $g_j$ and $v_j$ satisfy the balancing condition
$\sum_{j=1}^s(g_j-v_j)=0$ together with the relations
$t_j=e^{2\pi ig_j/\omega_2},$ $w_j=e^{2\pi i v_j/\omega_2}$,
 $\tilde t_j=re^{-2\pi ig_j/\omega_1},$
$\tilde w_j=re^{-2\pi i v_j/\omega_1}$,
and $y_{1,2}=e^{c\omega_{1,2}}$.

The case $|q|>1$ appears to be equivalent to the case $|q|<1$ after a
change of parameters and leads to the integrals
\be
\int_C e^{cu+d} \prod_{j=1}^{s}\frac{\eg(q^{-1}w_je^{2\pi i u/\omega_2 };p,q^{-1})}
{\eg(q^{-1}t_je^{2\pi i u/\omega_2 };p,q^{-1})}
\prod_{j=1}^{\ell }
\frac{\Gamma({\tilde q}^{-1}\tilde w_je^{-2\pi i u/\omega_1 }; {\tilde q}^{-1},r)}
{\Gamma({\tilde q}^{-1}\tilde t_je^{-2\pi i u/\omega_1 }; {\tilde q}^{-1},r)}
\, du.
\lab{e-meijer'}\ee
Functions \eqref{e-meijer}, \eqref{delta-ell-unit}, and \eqref{e-meijer'}
can be called as elliptic analogues of the Meijer function, because for
some particular choice of parameters and of the integration contour
$C$ they degenerate to that function \cite{erd:higher}. During this
degeneration procedure, at the intermediate steps there appear various
$q$-analogues of the Meijer function, including the cases considered in
\cite{slater}. The more general theta hypergeometric analogues of the
Meijer function, for which the kernels $\Delta(u)$
satisfy the system of equations \eqref{ell-gip} with $h_k(u)$
given by arbitrary meromorphic theta functions,
are built in \cite{spi:theta2}; we do not consider them here.

\subsection{Well poised and very-well poised integrals}\

We consider integrals  \eqref{e-meijer} with $\ell=c=d=0$ and replace
the integration variable by  $z=e^{2\pi iu/\omega_2}$.
Until now we did not fix the integration contour $C$. Let us choose
it as the unit circle $\T$ oriented counterclockwise. As a result, we
obtain integrals of the form
$$
\int_\T\Delta(z)\frac{dz}{z},\qquad \Delta(z)=\prod_{k=1}^s\frac{\Gamma(t_kz;p,q)}
{\Gamma(w_kz;p,q)}
$$
with $\prod_{k=1}^st_k=\prod_{k=1}^sw_k$.
In the case when the conditions $w_kt_k=pq$, $k=1,\ldots,s$ are
satisfied, the integrals take the form
$$
\int_\T\Delta^{(s)}(z)\frac{dz}{z},
\qquad \Delta^{(s)}(z)= \prod_{k=1}^{s} \Gamma(t_kz,t_k/z; p,q),
$$
and are called well poised integrals. The balancing condition
for them takes the form $\prod_{k=1}^{s} t_k^2= (pq)^{s}$ or
$\prod_{k=1}^{s} t_k= \mu (pq)^{s/2}$ with the ambiguity in the
sign choice $\mu =\pm 1$. The reflection formula
$\Gamma(a,b; p,q)=1$,  $ab=pq$, shows that the choice of parameters
$t_jt_k=pq$ plays an essential role, since it reduces the number
of parameters in $\Delta^{(s)}(z)$. In particular, for $t_k^2=pq$
the variable  $t_k$ drops out of formulae completely. The function
$$
h^{(p)}(z):=\frac{\Delta^{(s)}(qz)}{\Delta^{(s)}(z)}
=\prod_{k=1}^{s}\frac{\theta(t_kz;p)}{\theta(pqz/t_k;p)}
$$
is evidently $p$-elliptic, that is $h^{(p)}(pz)=h^{(p)}(z)$.
Denoting $u_k:=t_kz,\, v_k:=pqz/t_k,\, \lambda:=pqz^2$, we can
rewrite $h^{(p)}(z)$ as
\be
h^{(p)}(u_1,\ldots,u_{s};\lambda)
=\prod_{k=1}^{s}\frac{\theta(u_k;p)}{\theta(v_k;p)}
\lab{wp-ef}\ee
with the conditions of well poisedness $u_kv_k=\lambda$, $k=1,\ldots,s$,
and balancing $\prod_{k=1}^{s}u_k=\mu \lambda^{s/2}.$

Let us consider all possible $p$-shifts of the parameters
$u_1,\ldots, u_{s}$ and $\lambda$:
$$
u_k\to p^{n_k}u_k,\qquad \lambda\to p^{N}\lambda, \quad n_k,N\in\Z,
$$
and require that $h^{(p)}$ is invariant under the maximally
possible group of these transformations. The balancing condition
leads to the constraint
$
\sum_{k=1}^{s}n_k=sN/2.
$
For $N=0$ it is easy to check that
$h^{(p)}(\ldots,pu_a,\ldots,p^{-1}u_b,\ldots;\lambda)=
h^{(p)}(u_1,\ldots,u_{s}; \lambda)$,
i.e. $h^{(p)}$ is an elliptic function of all its parameters.
The transformations with $N\neq0$ are more complicated and
depend on  the parity of the variable $s$. For odd $s$
the integer $N$ must be even. The full symmetry group is generated then by
the transformations with $N=0$ and, say, $n_1=s,\,
n_k=0, k\neq1,$ and $N=2$ yielding
\begin{eqnarray*}
\frac{h^{(p)}(p^su_1,u_2,\ldots,u_{s};p^2\lambda)}{h^{(p)}(u_1,\ldots,u_s;\lambda)}
=\frac{\lambda^s}{\prod_{k=1}^{s}u_k^2}=1,
\end{eqnarray*}
that is the value of  $\mu$ is not fixed.

As to the even values $s=2m$, in this case there are no constraints
on $N$, and it is sufficient to consider the transformation corresponding
to the choice $n_1=m,\, N=1$, with all other $n_k=0$. Then we have
\begin{eqnarray*}
\frac{h^{(p)}(p^mu_1,u_2,\ldots,u_{2m};p\lambda)}
{h^{(p)}(u_1,\ldots,u_{2m};\lambda)}
=\frac{\lambda^m}{\prod_{k=1}^{2m}u_k}=\mu.
\end{eqnarray*}
Requiring this transformation to be a symmetry,
we fix uniquely the balancing condition, $\mu=1$.

It is not difficult to check that for any $\mu=\pm1$
the following equality is true
$$
\Delta^{(s)}(p^i q^jz)\Delta^{(s)}(z)=\Delta^{(s)}(p^i z)\Delta^{(s)}(q^jz)
$$
for all $i,j\in\Z$. Vice versa, from this condition one can
derive the balancing condition with $\mu=\pm1$.
Passing to the limits  $z\to \pm p^{-i/2}q^{-j/2}$ and using the
symmetry $\Delta^{(s)}(z)=\Delta^{(s)}(z^{-1})$, we obtain
$\Delta^{(s)}(\pm p^{i/2} q^{j/2})^2=\Delta^{(s)}(\pm p^{i/2}q^{-j/2})^2$.
The straightforward computation yields
$$
\Delta^{(s)}(\pm p^{i/2} q^{j/2})=
\left((\mp1)^s \mu\right)^{ij} \Delta^{(s)}(\pm p^{i/2}q^{-j/2}).
$$
For even values $s=2m$ and $\mu=1$, we obtain
\begin{equation}
\Delta^{(2m)}(\pm p^{i/2} q^{j/2})=\Delta^{(2m)}( \pm p^{i/2}q^{-j/2}),
\label{delta-plus}\end{equation}
provided the functions on both sides  are well defined.
The latter requirement is satisfied provided we do not hit the poles,
that is $t_r\neq \pm p^{a/2} q^{b/2}$, $a,b\in\Z$, for all $r$.

We call the integrals $\int_\T\Delta(z)dz/z$ very-well poised, if their
integration kernels have the form
$$
\Delta^{(2m+6)}_{vwp}(z)=\frac{\prod_{k=1}^{2m+6} \Gamma(t_kz^{\pm1};p,q)}
{\Gamma(z^{\pm2};p,q)}
=\theta(z^2;p)\theta(z^{-2};q)\prod_{k=1}^{2m+6} \Gamma(t_kz^{\pm1};p,q).
$$
This kernel can be obtained from $\Delta^{(2m+14)}(z)$ by
restriction of the parameters
$$
t_{2m+7},\dots,t_{2m+14}=(\pm(pq)^{1/2},\pm p^{1/2}q, \pm pq^{1/2}, \pm pq);
$$
this follows from the reflection formula and the argument duplication
formula for the function $\Gamma(z;p,q)$:
$$
\Gamma(z^2;p,q)=\Gamma(\pm z, \pm q^{1/2}z, \pm p^{1/2}z,
 \pm (pq)^{1/2}z;p,q).
$$
The balancing condition takes now the form
\begin{equation}
\prod_{k=1}^{2m+6}t_k=\mu (pq)^{m+1},
\label{balance} \end{equation}
where we count the sign choice $\mu=1$ as canonical,
since it leads to additional symmetries.
The kernel $\Delta_{vwp}^{(2m+6)}(z)$ can be obtained also
from $\Delta^{(2m+12)}(z)$ by imposing the constraints
$
t_{2m+7},\dots,t_{2m+12}=(\pm p^{1/2}q, \pm pq^{1/2}, \pm pq),
$
since the choice $t_k^2=pq$ simply removes the corresponding
 $\Gamma$-factors. However, we shall not use such a reduction,
because it changes the sign in the balancing condition:
$\prod_{k=1}^{2m+6}t_k=-\mu (pq)^{m+1}.$

In the following we shall be studying the very-well poised elliptic
hypergeometric integrals of the form
$$
I^{(m)}(t_1,\ldots,t_{2m+6})=\kappa \int_\T \frac{ \prod_{k=1}^{2m+6}
\Gamma(t_kz^{\pm 1};p,q) }{ \Gamma(z^{\pm 2};p,q)}
\frac{dz}{z}
$$
with the ``correct" balancing condition $\prod_{k=1}^{2m+6}t_k=(pq)^{m+1}$.
These integrals represent elliptic analogues of the plain hypergeometric
functions $_{m+1}F_{m}$. In particular, for $m=0$ we obtain
the elliptic beta integral. The constraints
$t_{2m+5}=(pq)^{1/2},\, t_{2m+6}=-(pq)^{1/2}$
reduce these integrals to $I^{(m-1)}(t_1,\ldots,t_{2m+4})$ with the
balancing condition $\prod_{k=1}^{2m+4}t_k=-(pq)^{m}$.
Appearance of the ``--" sign on the right-hand side simply
indicates that these integrals should be considered as some
generalizations of $_{m+1}F_{m}$-functions, and not of
$_{m}F_{m-1}$. For example, such a choice in the elliptic
beta integral yields  after taking into account of the
relation $\Gamma(-pq;p,q)=2(-p;p)_\infty (-q;q)_\infty$ the following:
\begin{eqnarray*}
&& \frac{1}{2\pi i}\int_{C}\frac{\prod_{k=1}^4\Gamma(t_kz^{\pm1}; p,q)}
{\Gamma(z^{\pm2}; p,q)}\frac{dz}{z}
\\ && \makebox[2em]{}
=4(p^2;p^2)_\infty (q^2;q^2)_\infty
\prod_{1\leq j<k\leq4}
\Gamma(t_jt_k;p,q)\prod_{k=1}^4\Gamma(pqt_k^2; p^2,q^2),
\end{eqnarray*}
where $\prod_{k=1}^4t_k=-1$ and the contour $C$ separates sequences of poles
converging to zero from those going to infinity. For $p\to 0$
we obtain a special case of the Askey--Wilson integral.

The multiplier appearing in integrals' kernel from the
very-well poisedness constraints for parameters leads to an
interesting property \cite{rai:trans}
\begin{equation}
\Delta_{vwp}^{(2m+6)}(p^{i/2}q^{j/2})=
-\mu^{ij}\Delta_{vwp}^{(2m+6)}(p^{i/2}q^{-j/2}),
\label{delta-minus}\end{equation}
that is to the ``--" sign in the right-hand side for the canonical
choice $\mu=1$, which is sharply distinct from \eqref{delta-plus}.
Indeed, for $i=0$ or $j=0$ we have zeros in both parts of the
equality, and in other cases we obtain
$$
\frac{\lim_{z\to p^{i/2}q^{j/2}}\Gamma^{-1}_{\! p,q}(z^2,z^{-2})}
{\lim_{z\to p^{i/2}q^{-j/2}}\Gamma^{-1}_{\! p,q}(z^2,z^{-2})}
=\frac{ \theta(p^{-i}q^{-j};q)\theta(p^i q^j;p) }
{\theta(p^{-i}q^j;q)\theta(p^i q^{-j};p)}
= -(pq)^{-2ij},
$$
which together with the described properties of $\Delta^{(2m+6)}(z)$
yields the presented formula. For $i,j\neq0$ there appears a
non-commutativity of two limits:
\begin{eqnarray*}
&& \frac{ \lim_{z\to p^{i/2}q^{j/2}}\Delta_{vwp}^{(2m+6)}(z) }
     { \lim_{z\to p^{i/2}q^{-j/2}}\Delta_{vwp}^{(2m+6)}(z) }=-\mu^{ij}
\\ && \makebox[1em]{}
\neq
\lim_{t_{2m+7},\dots,t_{2m+14}\to(\pm(pq)^{1/2},\pm p^{1/2}q, \pm pq^{1/2}, \pm pq)}
\frac{ \Delta^{(2m+14)}(p^{i/2}q^{j/2}) }
     { \Delta^{(2m+14)}(p^{i/2}q^{-j/2}) }=\mu^{ij},
\end{eqnarray*}
although both of them are well defined. The reason for such a contradiction
consists in the use in the very-well poisedness condition of the
``forbidden" values of parameters leading to poles and zeros of
$\Delta^{(2m+14)}(z)$ for $z=p^{i/2}q^{j/2}$.

As shown in \cite{rai:trans}, the product
$\prod_{1\leq j\leq k\leq m}(t_jt_k;p,q)_\infty I^{(m)}(\underline{t})$ is
a holomorphic function of parameters $t_j\in\C^*$. Additionally,
the ``--" sign on the right-hand side of \eqref{delta-minus} guarantees
that the integral $I^{(m)}(\underline{t})$ is holomorphic in the points
$t_k^2=p^{-a}q^{-b},$ $a,b=0,1,\ldots$. If some of the parameters
take the forbidden values $\pm p^{a/2} q^{b/2}$,
$a,b\in\Z$, then the latter property disappears.

\subsection{Series}\

According to the general definition  \cite{spi:theta1},
formal series $\sum_{n\in\Z}c_n$ are called the elliptic
hypergeometric series, if the ratio of neighbouring coefficients
$c_{n+1}/c_n$ is an elliptic function of $n\in\C$.
This definition lies in the stream of ideas of Pochhammer and Horn
which are used for building the plain and  $q$-hypergeometric
series \cite{ggr}.

As we saw already on the example of the Frenkel--Turaev sum,
the elliptic hypergeometric series appear as sums of
residues of certain sequences of poles of the elliptic
hypergeometric integrals' kernels. Indeed, let us consider the
poles of the integrand in \eqref{e-meijer} located at the
points $u=a+\omega_1 n,\, n=0,1,\ldots$, for some parameter $a$,
and denote residues of these poles as $c_n$.
For $u\to a+\omega_1 n$, we have $\Delta(u)\to c_n/(u-a-\omega_1n) +O(1)$.
Now it is not difficult to notice that the ratio
$$
\lim_{u\to a+\omega_1n}\frac{\Delta(u+\omega_1)}{\Delta(u)}=\frac{c_{n+1}}{c_n}=
\lim_{u\to a+\omega_1n} h_1(u)=h_1(a+\omega_1 n)
$$
is an elliptic function of $n$ with the periods $\omega_2/\omega_1$ and
$\omega_3/\omega_1$, which demonstrates the general connection between
the integrals and series.

An arbitrary elliptic function $h(n)$ of the order $s+1$ with the periods
$\omega_2/\omega_1$ and $\omega_3/\omega_1$ has the form
\begin{equation}
h(n)=y\;\prod_{k=1}^{s+1}\frac{\theta(t_kq^n;p)}{\theta(w_kq^u;p)},
\label{el-fn}\end{equation}
where $ \prod_{k=1}^{s+1} t_k=\prod_{k=1}^{s+1} w_k. $
We define the elliptic Pochhammer symbol $\theta(t)_n$ as the solution
of the recurrence relation $c_{n+1}=\theta(tq^n;p)\, c_n$
with the initial condition $c_0=1$:
$$
\theta(t)_n=\left\{
\begin{array}{cl}
\prod_{j=0}^{n-1}\theta(tq^j;p), &  n>0  \\ [0.5em]   \displaystyle
\prod_{j=1}^{-n}\frac{1}{\theta(tq^{-j};p)}, & n<0.
\end{array}
\right.
$$
Then it is easy to deduce the explicit form of the formal bilateral
elliptic hypergeometric series
$$
{}_{s+1}G_{s+1}\bigg({t_1,\ldots,t_{s+1}\atop w_1,\ldots, w_{s+1}};q,p;y\bigg)
=\sum_{n\in\Z}\prod_{k=1}^{s+1}\frac{\theta(t_k)_n}
{\theta(w_k)_n}\,y^n
$$
with the normalization of the zeroth coefficient $c_0=1$.
Choosing $w_{s+1}=q$ and replacing $t_{s+1}\to t_0$,
we obtain the unilateral series
\begin{equation}
{}_{s+1}E_s\bigg({t_0,t_1,\ldots  ,t_s\atop w_1,\ldots, w_s};q,p;y\bigg)
=\sum_{n=0}^\infty\frac{\theta(t_0,t_1,\ldots, t_s)_n}{
\theta(q,w_1,\ldots,w_s)_n}\, y^n.
\label{E-series}\end{equation}
For fixed $t_j$ and $w_j$, in the limit $p\to 0$ we have
$\theta(t)_n\to (t;q)_n= (1-t)(1-qt)\ldots(1-q^{n-1}t)$,
and the function ${}_{s+1}E_s$ reduces to the
$q$-hypergeometric series \cite{gas-rah:basic}
$$
{}_{s+1}\varphi_s\bigg({t_0,t_1,\ldots,t_s\atop w_1,\ldots, w_s};q;y\bigg)
=\sum_{n=0}^\infty\frac{(t_0;q)_n\ldots (t_s;q)_n}
{(q;q)_n(w_1;q)_n\ldots (w_s;q)_n}\, y^n
$$
with the condition $\prod_{j=0}^st_j=q\prod_{j=1}^sw_j$.
Parameterizing $t_j=q^{u_j}$ and $w_j=q^{v_j}$, in the limit $q\to1$
we obtain the series
$$
{}_{s+1}F_s\bigg({u_0,u_1,\ldots,u_s\atop v_1,\ldots, v_s};y\bigg)
=\sum_{n=0}^\infty\frac{(u_0)_n\ldots (u_s)_n}
{n!(v_1)_n\ldots (v_s)_n}\, y^n,
$$
where $(a)_n=a(a+1)\ldots (a+n-1)$ denotes the standard Pochhammer symbol
and  $u_0+\ldots +u_s=1+v_1+\ldots+v_s$. The latter constraint for the
parameters is not essential, since it disappears already for the
 $_sF_{s-1}$-function obtained after taking the limit $u_s\to\infty$.

Investigation of the conditions of convergence of the infinite series
\eqref{E-series} represents a serious problem and it was not solved
completely to the present moment. Therefore in the applications
of the $E$-series it is usually assumed that they terminate
because of the condition $t_k=q^{-N}p^M$ for some $k$, where
$N=0,1,\ldots,\, M\in\Z.$ It is worth of noting that the formal
 $_2E_1$-series does not represent a natural elliptic generalization
of the Gauss hypergeometric function, because it does not obey
natural analogues of many important properties of the $_2F_1$-function.

Series  \eqref{E-series} are called well poised, if the
following constraints on the parameters are satisfied
$t_0q=t_1w_1=\ldots=t_sw_s$. The balancing condition for
them takes the form
$t_1\cdots t_s=\pm  q^{(s+1)/2}t_0^{(s-1)/2}$, and the
functions $h(n)$ and $_{s+1}E_s$ become invariant with respect
to the transformations $t_j\to pt_j,\, j=1,\ldots,s-1,$
and $t_0\to p^2t_0$ (for this it is necessary to count
$t_s$ as a dependent parameter). In the same way as in the
case of integrals, for odd $s$ and the ``+" sign choice in
the balancing condition there appears an additional symmetry
--- the functions $h(n)$ and $_{s+1}E_s$ become invariant
with respect to the transformation $t_0\to pt_0$
(with the compensating transformation $t_s\to p^{(s-1)/2} t_s$).
The function $h(n)$ becomes thus an elliptic function of
all free parameters $\log t_j,\, j=0,\ldots, s-1,$
with equal periods, that is there appears some kind
of ``total ellipticity" \cite{spi:theta1,spi:thesis}.

The next structural restriction, which is needed for building the most
interesting examples of the series, looks as follows
$$
t_{s-3}=q\sqrt{t_0},\quad t_{s-2}=-q\sqrt{t_0},\quad t_{s-1}=
q\sqrt{t_0/p},\quad t_s=-q\sqrt{pt_0}
$$
and is called the very-well poisedness condition
(it is related to doubling of the argument for theta
functions). In view of the importance of the very-well
poised elliptic hypergeometric series, there is a special notation
for them \cite{spi:bailey1}:
\begin{eqnarray}\nonumber
&&{}_{s+1}V_{s}(t_0;t_1,\ldots,t_{s-4};q,p;y)= \sum_{n=0}^\infty
\frac{\theta(t_0q^{2n};p)}{\theta(t_0;p)}\prod_{m=0}^{s-4}
\frac{\theta(t_m)_n}{\theta(qt_0t_m^{-1})_n}(qy)^n
\\ \nonumber &&
={}_{s+1}E_{s}\bigg({t_0,t_1,\ldots,t_{s-4},q\sqrt{t_0},-q\sqrt{t_0},
q\sqrt{t_0/p},-q\sqrt{pt_0}  \atop
qt_0/t_1,\ldots  ,qt_0/t_{s-4},\sqrt{t_0}, -\sqrt{t_0},\sqrt{pt_0},
-\sqrt{t_0/p}};q,p;-y\bigg),
\label{V-series}\end{eqnarray}
where $\prod_{k=1}^{s-4}t_k=\pm t_0^{(s-5)/2}q^{(s-7)/2}$
(in this balancing condition for odd $s$ the ``+" sign choice
is considered as canonical). The choice $t_j=\pm\sqrt{qt_0}$
removes completely this parameter and reduces the indices of the series
by 1 (this choice can change the sign in the balancing condition).
For sufficiently large values of $s$, one can choose four
parameter values in such a way that the very-well poised part in the
series coefficients cancels out, and it becomes again only the
well poised series. For $y=1$, this argument is dropped
in the notation of ${}_{s+1}V_{s}$-series. In this scheme,
the Frenkel--Turaev formula
\eqref{ft-sum} yields a closed form expression for the terminating
$_{10}V_9(t_0;t_1,\ldots,t_5;q,p)$-series.

In order to consider modular transformations it is necessary to use the
parameterization  $t_k=q^{u_k}$,  $w_k=q^{v_k}$,
$q=e^{2\pi i\sigma}$,  $p=e^{2\pi i\tau}$ and replace
$\theta(q^a;p)$-functions by the ``elliptic numbers"
$[a]=\theta_1(\sigma a|\tau)$. Suppose that the parameter $y$ does
not depend on $\tau$. Then, it is not difficult to verify with the help of
formula \eqref{mod} that
the functions $h(n)$ and, respectively, $_{s+1}E_s$ will be
modular invariant under the restriction $\sum_{k=0}^su_k^2=1+\sum_{k=1}^s v_k^2$
defined modulo $2\tau/\sigma^2$.
For completeness of the description, we present
 an explicit expression for the most interesting $_{s+1}V_s$-series
with odd $s$ in this notation:
\begin{eqnarray*}\nonumber
&&{}_{2m}V_{2m-1}(q^{u_0};q^{u_1},\ldots,q^{u_{2m-5}}; q,p;y)
= \sum_{n=0}^\infty \frac{[2n+u_0]}{[u_0]}\prod_{k=0}^{2m-5}
\frac{[u_k]_n}{[1+u_0-u_k]_n}y^n,
\end{eqnarray*}
where $[a]_n=[a][a+1]\ldots [a_n-1]$ and
the balancing condition $\sum_{k=1}^{2m-5}u_k=m-4+(m-3)u_0$.
This series is automatically modular invariant.
In \cite{gas-rah:basic}, the symbols $_{s+1}e_s$ and $_{s+1}v_s$
were suggested for the additive system of notation
for series, but we do not use them here.

Importance of the balancing condition for the plain and
$q$-hypergeometric functions was known for a long time,
because in its presence there appear some additional identities
\cite{gas-rah:basic}. The same is true for the notions of
well poisedness and very-well poisedness. However, the corresponding
constraints on the parameters were appearing in a rather ad hoc manner,
and their deep meaning was missing. The elliptic hypergeometric
functions clarify the origin of these old concepts. Namely,
the balancing condition is connected with the condition of double
periodicity of the main elliptic function used in the
construction of series or integrals. The condition of well poisedness
is connected with the condition of ellipticity in all parameters
determining the divisor of this elliptic function. The condition
of very-well poisedness is tied to the rule of doubling of the argument
of theta functions. Strictly speaking all these notions
are defined in fact in a self-contained manner only at the elliptic level.
Indeed, there are limiting transitions from the elliptic hypergeometric
identities involving $_{s+1}V_s$-series
to the $q$-hypergeometric relations such that there appear basic
$_{s+1}\varphi_s$-functions which do not obey one of the mentioned
properties \cite{rai:abelian,spi:bailey1}.

The degeneration of  $_{s+1}E_s$-series
described above at $p\to 0$ leads to $_{s+1}\varphi_s$-series with
the constraint for parameters resembling
the old $q$-balancing condition \cite{aar,gas-rah:basic}, but not
coinciding with it. The limit $p\to 0$ for $_{s+1}V_s$-series
with fixed parameters leads to the very-well poised balanced
 $_{s-1}\varphi_{s-2}$-series having their own notation
\cite{gas-rah:basic}:
\begin{eqnarray*}
&&\lim_{p\to 0}\; {}_{s+1}V_s(t_0;t_1,\ldots,t_{s-4};q,p;z)
=\sum_{n=0}^\infty\frac{1-t_0q^{2n}}{1-t_0}\prod_{k=0}^{s-4}\frac{(t_k;q)_n}
{(qt_0/t_k;q)_n}(qz)^n
\\ && \makebox[4em]{}
=:{}_{s-1}W_{s-2}(t_0;t_1,\ldots,t_{s-4}; q,qz).
\end{eqnarray*}
In a remarkable way, the balancing condition for the
$_{s+1}V_s$-series coincides in this case with the usual
balancing condition for the $_{s-1}W_{s-2}$-series
\cite{gas-rah:basic,spi:theta1}.
The possibility of fixing the sign in the balancing condition
for odd $s$ from the requirement of existence of
an additional symmetry strengthens the ``elliptic" point of view
on the functions of hypergeometric type and indicates
on an indispensable connection of these two classes of functions.

In conclusion of this section, let us mention that the quadratic
transformations for the $_{s+1}V_s$-series
and related summation formulae were considered in
\cite{chu,spi:bailey1,war:summation,war:summation2}.
Other specific elliptic hypergeometric series were investigated
in the papers \cite{gs:theta,ros-sch2}. An interesting application
of a $_{s+1}E_s$-series with the nontrivial power variable $y$
appeared recently in \cite{zhe:ell}.

\section{An elliptic analogue of the Gauss hypergeometric function}

\subsection{Definition of the $V$-function and a connection with
the root system $E_7$}\

The Euler integral representation  for the
$_2F_1(a,b;c;x)$-function differs from the beta integral by the
presence in the integrand of an additional term depending on two
new parameters \cite{aar}. An elliptic analogue of the
Gauss hypergeometric function, which we shall be denoting by the
symbol $V(t_1,\ldots,t_8;p,q)$, is also given by a two-parameter
extension of the elliptic beta integral \cite{spi:thesis}:
\be
V(t_1,\ldots,t_8;p,q)=\kappa\int_\T\frac{\prod_{j=1}^8
\eg(t_jz^{{\pm 1}};p,q)}{\eg(z^{\pm 2};p,q)}\frac{dz}{z},
\lab{ehf}\ee
where eight parameters $t_1,\ldots, t_8\in\C $ and two basic variables
$p,q\in\C$ satisfy the constraints $|t_j|,|p|,|q|<1$ and
the balancing condition $\prod_{j=1}^8 t_j=p^2q^2$.
For other values of the parameters $t_j$ the $V$-function is
defined by analytical continuation of the integral \eqref{ehf}.
This continuation is build by the replacement of the contour of integration
$\T$ by a contour $C$ separating the sequences of poles
of the integrand at $z=t_jp^{a}p^{b},\, a,b=0,1,\ldots$,
converging to zero from the poles at $z=t_j^{-1}p^{-a}p^{-b}$,
diverging to infinity; this does not assume now the restrictions
 $|t_j|<1$. Shrinking the contour $C$, one can represent the resulting
function as an integral over $\T$ and the sum of residues
crossed by the contour during this deformation.

Suppose that a pair of parameters satisfies the condition
 $t_jt_k=pq$. Then the $V$-function is reduced to the elliptic beta integral
\eqref{ell-int}, which follows from the reflection formula
for the elliptic gamma function $\Gamma(z;p,q)$.

Let us consider symmetries of the $V$-function.
Evidently, it is invariant with respect to the permutation of
$p$ and $q$ and $S_8$-group of permutations of $t_j$
related to the root system $A_7$. It appears that there
exists a transformation extending $S_8$ to the Weyl group
for the exceptional root system $E_7$. It is derived with
the help of the double integral
$$
\kappa \int_{\T^2}
\frac{\prod_{j=1}^4\eg(a_jz^{\pm 1}, b_jw^{\pm 1};p,q)\;\eg(cz^{\pm 1} w^{\pm 1};p,q)}
{\eg(z^{\pm 2},w^{\pm 2};p,q)} \frac{dz}{z}\frac{dw}{w},
$$
where $a_j,b_j,c\in \C $, $|a_j|,|b_j|,|c|<1,$ and
$c^2\prod_{j=1}^4a_j= c^2\prod_{j=1}^4b_j=pq.$
Computation of the integrals over $z$ or over $w$
in different orders with the help of formula
\eqref{ell-int} yields the fundamentally important relation \cite{spi:theta2}
\be
V(\underline{t})=\prod_{1\le j<k\le 4}\eg(t_jt_k,t_{j+4}t_{k+4};p,q)\,
V(\underline{s}),
\lab{E7-1}\ee
where $V(\underline{t})=V(t_1,\ldots, t_8;p,q)$ and
$$
\left\{
\begin{array}{cl}
s_j =\ve t_j,&   j=1,2,3,4  \\
s_j = \ve^{-1} t_j, &    j=5,6,7,8
\end{array}
\right.;
\quad \ve=\sqrt{\frac{pq}{t_1t_2t_3t_4}}=\sqrt{\frac{t_5t_6t_7t_8}{pq}}
$$
and $|t_j|, |s_j|<1$.

Repetition of transformation \eqref{E7-1} with the parameters $s_{3,4,5,6}$
playing the role of $t_{1,2,3,4}$ and subsequent permutation of
the parameters $t_3,t_4$ with $t_5,t_6$ result in the equality
\be
V(\underline{t})=\prod_{j,k=1}^4
\eg(t_jt_{k+4};p,q)\ V(T^{1\over 2}\!/t_1,\ldots,T^{1\over 2}\!/t_4,
U^{1\over 2}\!/t_5,\ldots,U^{1\over 2}\!/t_8),
\lab{E7-2}\ee
where $ T=t_1t_2t_3t_4$, $ U=t_5t_6t_7t_8$ and
$|T|^{1/2}<|t_j|<1,$ $|U|^{1/2}<|t_{j+4}|<1,\, j=1,2,3,4$.
Equating the right-hand sides of equalities
 \eqref{E7-1} and \eqref{E7-2} and expressing parameters
$t_j$ in terms of $s_j$, we obtain the third transformation
\be
V(\underline{s})=\prod_{1\le j<k\le 8}\eg(s_js_k;p,q)\,
V(\sqrt{pq}/s_1,\ldots,\sqrt{pq}/s_8),
\lab{E7-3}\ee
where $|pq|^{1/2}<|s_j|<1$ for all $j$.

Let us connect parameters of the function $V(\underline{t})$
with the coordinates of an Euclidean space $x_j\in\R^8 $
by the relations $t_j=e^{2\pi i x_j}(pq)^{1/4}$.
Denote as $\langle x, y\rangle$ the scalar product in
$\R^8 $ and as $e_i$ -- an orthonormal basis,
$\langle e_i, e_j\rangle=\delta_{ij}.$ The root system $A_7$
consists of the vectors $v=\{e_i-e_j,\, i\neq j\},$
and its Weyl group $S_8$ -- from the reflections
$x\to S_v(x)=x-2v\langle v, x\rangle/\langle v, v\rangle$,
acting in the hyperplane orthogonal to the vector
 $\sum_{i=1}^8e_i$.  In this case the coordinates of the vectors
 $x=\sum_{i=1}^8x_ie_i$ satisfy the relation $\sum_{i=1}^8x_i=0$,
which is guaranteed by the balancing condition.

Transformation of the coordinates in \eqref{E7-1}
corresponds to the reflection $S_v(x)$ with respect to the
vector $v=(\sum_{i=5}^8e_i-\sum_{i=1}^4 e_i)/2$ of the length
$\langle v, v\rangle=2$. It extends the group
 $S_8$ to the Weyl group for the exceptional root system $E_7$.
A different proof of  equalities \eqref{E7-1}--\eqref{E7-3}
is given by Rains in \cite{rai:trans}, where a relation to the
$E_7$-group is indicated for the first time.

In \cite{spi:theta2} it was shown that the $V$-function is reduced to
the product of two $_{12}V_{11}$-series for special discrete values of
one of the parameters. Let us denote
$A=t_0\ldots t_4$; then
\begin{eqnarray}\nonumber
&& V(t_0,\ldots,t_5,t_5^{-1}p^{n+1}q^{m+1},A^{-1}p^{1-n}q^{1-m};p,q)
\\ && \makebox[1em]{}
=\frac{\prod_{0\leq j<k\leq 4}\Gamma(t_jt_k;p,q)}
{\prod_{j=0}^4\Gamma(t_j^{-1}A;p,q)}
\frac{\Gamma(p^{n+1}q^{m+1}t_5^{-1}t_0^{\pm1},t_5t_0^{\pm1},At_0^{\pm1};p,q)}
{\Gamma(p^nq^m At_0^{\pm1};p,q)}
  \nonumber \\ && \makebox[1em]{}\times
{_{12}V_{11}}\left(\frac{At_0}{q};\frac{At_5}{q},t_0t_1,t_0t_2,t_0t_3,t_0t_4,
q^{-m},\frac{Aq^{m}}{t_5};q,p\right) \nonumber \\
&& \makebox[1em]{} \times {_{12}V_{11}}
\left(\frac{At_0}{p};\frac{At_5}{p},t_0t_1,t_0t_2,t_0t_3,t_0t_4,
p^{-n},\frac{Ap^{n}}{t_5};p,q\right),
\label{intrep}\end{eqnarray}
where the contour of integration in the definition
of $V$-function is chosen in such a way that it separates sequences
of the integrand poles converging to zero and diverging to infinity.
For $m=0$, we obtain an integral representation of a separate terminating
$_{12}V_{11}$-series. Since the left-hand side of \eqref{intrep}
is symmetric in $t_0,\ldots, t_4$, the same should hold for the
right-hand side as well, which follows from the symmetry transformation
for  series to be considered below.

\subsection{The elliptic hypergeometric equation}\

The addition formula for theta functions \eqref{ident}, being written in the form
$$
t_8\theta(t_7t_8^{\pm1},t_6z^{\pm1};p)+t_6\theta(t_8t_6^{\pm1},t_7z^{\pm1};p)
+t_7\theta(t_6t_7^{\pm1},t_8z^{\pm1};p)=0,
$$
leads to the connection formula
\begin{equation}
\frac{t_6V(qt_6)}{\theta(t_6t_7^{\pm1},t_6t_8^{\pm 1};p)}
+\frac{t_7V(qt_7)}{\theta(t_7t_6^{\pm1},t_7t_8^{\pm 1};p)}
+\frac{t_8V(qt_8)}{\theta(t_8t_6^{\pm1},t_8t_7^{\pm 1};p)} = 0,
\label{cont-1}\end{equation}
where $V(qt_j)$ denotes the  $V$-function with the parameter
$t_j$ replaced by  $qt_j$ (this leads to the balancing condition
$\prod_{j=1}^8t_j=p^2q$). Indeed, it is easy to check that
the same equality holds for the $V$-function kernel,
after integration of which over $z$ one obtains \eqref{cont-1}.

A substitution of transformation \eqref{E7-3} in \eqref{cont-1} yields
\begin{eqnarray}\nonumber
&& \frac{\prod_{j=1}^5\theta\left(t_6t_j/q;p\right)}
{t_6\theta(t_7/t_6,t_8/t_6;p)}V(q^{-1}t_6)
+\frac{\prod_{j=1}^5\theta\left(t_7t_j/q;p\right)}
{t_7\theta(t_6/t_7,t_8/t_7;p)}V(q^{-1}t_7)
\\ &&\makebox[2em]{}
+\frac{\prod_{j=1}^5\theta\left(t_8t_j/q;p\right)}
{t_8\theta(t_6/t_8,t_7/t_8;p)}V(q^{-1}t_8) = 0,
\label{cont-3}\end{eqnarray}
where $\prod_{j=1}^8t_j=p^2q^3$.

Let us consider three equations appearing after the replacements
 $t_8\to qt_8$ in \eqref{cont-3} and $t_6\to q^{-1}t_6$ or
$t_7\to q^{-1}t_7$ in \eqref{cont-1}. Excluding from them the
functions $V(q^{-1}t_6,qt_8)$ and $V(q^{-1}t_7,qt_8)$,
we obtain the elliptic hypergeometric equation \cite{spi:thesis,spi:cs}:
\begin{eqnarray}\label{eheq}
&& \makebox[-2em]{}
\mathcal{A}(t_1,\ldots,t_6,t_7,t_8,q;p)\Big(U(qt_6,q^{-1}t_7)
-U(\underline{t})\Big)
\\ &&
+\mathcal{A}(t_1,\ldots,t_7,t_6,t_8,q;p)\Big(U(q^{-1}t_6,qt_7)
-U(\underline{t})\Big) + U(\underline{t})=0,
\nonumber\end{eqnarray}
where
\begin{equation}
 \mathcal{A}(t_1,\ldots, t_8,q;p):=\frac{\theta(t_6/qt_8,t_8t_6,t_8/t_6;p)}
                 {\theta(t_6/t_7,t_7/qt_6,t_6t_7/q;p)}
\prod_{k=1}^5\frac{\theta(t_7t_k/q;p)}{\theta(t_8t_k;p)}
\end{equation}
and
$$
U(\underline{t}):=U(\underline{t};p,q)=\frac{V(\underline{t};p,q)}
{\prod_{k=6}^7 \eg(t_kt_8^{\pm 1};p,q)}.
$$
The potential $\mathcal{A}(t_1,\ldots, t_8,q;p)$ is an elliptic
function of all its parameters, i.e. it does not change under the
transformations $t_j\to p t_j$, $j=1,\ldots,7$, provided we count $t_8$
as a parameter depending on others through the balancing condition.
Because of the symmetry of the $U$-function in $p$ and $q$, it satisfies
another elliptic hypergeometric equation, which is obtained from
\eqref{eheq} by permutation of $p$ and $q$.
A change of variables
brings the function $\mathcal{A}(t_1,\ldots, t_8,q;p)$
to the completely $S_8$-symmetric form \eqref{wp-ef}
with arbitrary $u_1,\ldots, u_8$ and  $\mu=1$
(a remark from D. Zagier to the author).
However, the elliptic hypergeometric
equation itself is maximally $S_6$-symmetric \cite{spi:thesis,spi:cs}.

Let us denote $t_6=cx, t_7=c/x$ and pass to a new set of parameters
$$
\ve_k=\frac{q}{ct_k},\; k=1,\ldots,5,\quad \ve_8=\frac{c}{t_8},
\quad \ve_7=\frac{\ve_8}{q}.
$$
We fix the parameter $\ve_6$ from the requirement that
the balancing condition takes the form $\prod_{k=1}^8\ve_k=p^2q^2$,
which yields $c=\sqrt{\ve_6\ve_8}/p^2$. Now, after the replacement
of $U(\underline{t})$ by some unknown function
$f(x)$ in \eqref{eheq}, we can rewrite this elliptic hypergeometric
equation in the form of a $q$-difference equation of the second order:
\begin{eqnarray}
&& A(x)\left( f(qx)-f(x)\right)
+ A(x^{-1})\left( f(q^{-1}x)-f(x)\right) + \nu f(x)=0,
\label{eheq2}
\\ && \qquad
A(x)=\frac{\prod_{k=1}^8 \theta(\ve_kx;p)}{\theta(x^2,qx^2;p)},
\qquad
\nu=\prod_{k=1}^6\theta\left(\frac{\ve_k\ve_8}{q};p\right),
\lab{pot}\end{eqnarray}
where one can explicitly see the  $S_6$-group of symmetries in
parameters. We have already one functional solution of this equation:
\be
f_1(x;\underline{\ve};p,q)=\frac{ V(q/c\ve_1,\ldots,q/c\ve_5,cx,c/x,c/\ve_8;p,q)}
{\eg(c^2x^{\pm 1}/\ve_8,x^{\pm 1} \ve_8;p,q)}.
\lab{sol1}\ee
In order to build other linearly independent solutions, one can use
symmetries of the equation \eqref{eheq2} which do not represent symmetries
of the function \eqref{sol1}. For instance, the second solution can be obtained
by multiplication of one of the parameters $\ve_1,\ldots,\ve_5$ or
$x$ by the powers of $p$ or by permutations of $\ve_1,\ldots,\ve_5$ with $\ve_6$.

Since the function $\mathcal{A}(\underline{t},q;p)$ in equation
\eqref{eheq} does not change after the replacements
$t_6\to p^{-1}t_6, t_7\to pt_7$, the function $U(p^{-1}t_6,pt_7)$
also represents a solution of this equation.
The Casoratian (a discrete Wronskian) of these two solutions was computed
in the paper \cite{rs:det}. Let us multiply equation
\eqref{eheq} by $U(p^{-1}t_6,pt_7)$, and the equation for
$U(p^{-1}t_6,pt_7)$ by  $U(\underline{t})$ and subtract
one of the resulting relations from another.
As a result, we obtain the equality
\begin{equation}\label{cas-eqn}
\mathcal{A}(t_1,\ldots,t_6,t_7,t_8,q;p){D}(p^{-1}t_6,q^{-1}t_7)
=\mathcal{A}(t_1,\ldots,t_7,t_6,t_8,q;p){D}(p^{-1}q^{-1}t_6,t_7),
\end{equation}
where $D$ denotes the Casoratian
$$
{D}(t_6,t_7)=U(pqt_6,t_7)U(t_6,pqt_7)-U(qt_6,pt_7)U(pt_6,qt_7).
$$
After substitutions $t_6\to pt_6,\, t_7\to qt_7$ the balancing condition
takes the form $t_6t_7=pq/t_1\ldots t_5t_8$, and \eqref{cas-eqn} can be considered
as a $q$-difference equation of the first order in $t_7$.
It is easily solved, and the solution is defined up to the multiplication by
an arbitrary $q$-periodic function $\varphi(qt_7)=\varphi(t_7)$.

Repeating the same procedure with the equations appearing after the
permutation of the parameters $p$ and $q$, we find
$\varphi(pt_7)=\varphi(t_7)$,
so that $\varphi(t_7)$ does not depend on $t_7$.
Further investigation of the structure of pole residues of the
$V$-functions in ${D}(t_6,t_7)$, which are crossed by the integration
contour in the limit $t_7\to t_8^{-1}$, fixes completely the form of
$\varphi$ and leads to the identity \cite{rs:det}
\begin{equation}
V(pqt_6,t_7)V(t_6,pqt_7)-t_6^{-2}t_7^{-2} V(qt_6,pt_7)V(pt_6,qt_7)
=\frac{\prod_{1\leq j<k\leq8}\eg(t_jt_k;p,q)}{\eg(t_6^{\pm1}t_7^{\pm1};p,q)},
\label{V-det}\end{equation}
where $t_6$ and $t_7$ can be replaced by any other pair of parameters.

The described solutions of the elliptic hypergeometric equation exist
for $|q|<1$, whereas the equation itself \eqref{eheq} does not
demand such a restriction.  Because of the symmetry
$$
\mathcal{A}\left(\frac{p^{1/2}}{t_1}, \ldots,\frac{p^{1/2}}{t_8},q;p\right)
=\mathcal{A}\left(t_1,\ldots,t_8,q^{-1};p\right),
$$
the transformation $t_j\to p^{1/2}/t_j$, $j=1,\ldots, 8,$ leads
to the change of the base variable $q\to q^{-1}$ in \eqref{eheq}.
This yields the following solution of the elliptic hypergeometric
equation for  $|q|>1$ \cite{rs:det}
\begin{eqnarray}
U(\underline{t};p,q)=\frac{V(p^{1/2}/t_1,\ldots,p^{1/2}/t_8;p,q^{-1})}
{\prod_{k=6}^7\Gamma(p/t_kt_8,t_8/t_k;p,q^{-1}) }.
\label{q>1}\end{eqnarray}

In order to build solutions of equation \eqref{eheq} (or \eqref{eheq2})
for $|q|=1$, one can use the parameterization of the base
variables \eqref{bases}, \eqref{mod-bases}, make substitutions
$t_k=e^{2\pi i g_k/\omega_2}$ and repeat the whole chain of arguments
given above with the replacement of $\Gamma(z;p,q)$ by the
modified elliptic gamma function $G(u;\mathbf{\omega})$.
In the same way as in the case of modified elliptic beta integral,
we obtain again the $V$-function, but with a different parameterization.
This procedure appears to be equivalent to the application of
the modular transformation  $(\omega_2,\omega_3)\to (-\omega_3,\omega_2)$
to solutions described above.
Indeed, the potential $\mathcal{A}(e^{2\pi i g_1/\omega_2},
\ldots, e^{2\pi i g_8/\omega_2},q;p)$ is invariant under this
transformation, but $U(t_1,\ldots, t_8; p, q)$ gets transformed to
$U(e^{-2\pi i g_1/\omega_3},\ldots, e^{-2\pi i g_8/\omega_3};\tilde p,\tilde r)$.
The latter function provides thus a new solution of the elliptic hypergeometric
equation well defined for $|q|=1$. It is evident that this function satisfies
also a partner of equation \eqref{eheq} obtained from it
by the permutation of $\omega_1$ and $\omega_2$.
An example of a similar situation with two equations for
one function at the $q$-hypergeometric level is given in \cite{rui:generalized}.

Different formal degenerations of the
$V$-function to $q$-hypergeometric integrals of the Mellin--Barnes
or Euler type are briefly considered in
\cite{spi:thesis,spi:cs}. A detailed and rigorous
analysis of the degeneration procedure is performed
in \cite{rai:limits,brs}. It is necessary to note \cite{spi:cs} that
the elliptic hypergeometric equation emerges as a particular case
of the eigenvalue problem equation for the one particle Hamiltonian
of the quantum model proposed by van Diejen \cite{die:integrability}
and investigated in detail by Komori and Hikami \cite{kom-hik:quantum}.
This Calogero--Sutherland type model represents a generalization
of the Ruijsenaars \cite{rui:complete} and Inozemtsev \cite{ino:lax}
systems.

\section{Chains of symmetry transformations for functions}

Symmetry transformations for the plain and $q$-hypergeometric series are built
with the help of the Bailey chains \cite{aar,and:bailey}.
This technique was generalized to the elliptic level in
 \cite{spi:bailey1,war:extensions}. Let us sketch
this generalization.

Two sequences of numbers $\alpha_n(a,k)$ and $\beta_n(a,k)$
by definition form an elliptic Bailey pair with respect to the parameters
 $a$ and $k$, if
\begin{equation}
\beta_n(a,k)=\sum_{0\leq m\leq n} M_{nm}(a,k)\alpha_m(a,k),
\label{pair}\end{equation}
where
\begin{equation}
M_{nm}(a,k)=\frac{\theta(k/a)_{n-m}\theta(k)_{n+m}}
{\theta(q)_{n-m}\theta(aq)_{n+m}}\frac{\theta(aq^{2m};p)}
{\theta(a;p)}a^{n-m}.
\end{equation}
In the matrix form $\beta(a,k)=M(a,k)\alpha(a,k)$, where $\alpha$ and $\beta$
denote the columns formed by $\alpha_n$ and $\beta_n$.

Let us introduce the diagonal matrix
\begin{equation}
D_{nm}(a;b,c)=D_m(a;b,c)\delta_{nm},\qquad
D_m(a;b,c)=\frac{\theta(b,c)_m}{\theta(aq/b,aq/c)_m}
\left(\frac{aq}{bc}\right)^m.
\end{equation}

\begin{theorem}
Let $\alpha(a,t)$ and $\beta(a,t)$ form an elliptic Bailey pair with respect
to the parameters $a$ and $t$. Then the quantities
\begin{eqnarray}\label{rel1}
&& \alpha'(a,k)=D(a;b,c)\alpha(a,t),\qquad
\beta'(a,k)=K(t,k,b,c)\beta(a,t),
\\ \nonumber  && \makebox[2em]{}
K(t,k,b,c):= D(k;qt/b,qt/c)M(t,k)D(t;b,c),
\end{eqnarray}
where $qat=kbc$ and $b,c$ are two arbitrary new parameters,
form a new elliptic Bailey pair with respect to $a$ and $k$.
\end{theorem}
\begin{proof}
Substitution of \eqref{rel1} and of the relation $\beta=M\alpha$
in the required equality $\beta'=M\alpha'$ leads to the
matrix identity
\begin{equation}
M(a,k)D(a;b,c)M(t,a)=D(k;qt/b,qt/c)M(t,k)D(t;b,c).
\label{mat-id}\end{equation}
After substitution of the explicit expressions for matrices,
one can see that it is equivalent to the Frenkel--Turaev summation formula
\eqref{ft-sum}.
\end{proof}

Because $M_{nm}(a,a)=\delta_{nm}$ and $D_{nm}(bc/q;b,c)=\delta_{nm}$,
we find after setting $t=k$ in \eqref{mat-id} that $M(a,k)M(k,a)=1$.
The inversion of the matrix $M$ is reached thus by the permutation
of parameters $a$ and $k$ (in the $p=0$ case this fact was established
in \cite{bre}; a more detailed discussion of such matrix inversions
is given in \cite{war:summation,ros-sch2,ma}).
Therefore, $\tilde \alpha(a,k)=\beta(k,a)$
and $\tilde \beta(a,k)=\alpha(k,a)$ define new Bailey pairs to which
one can apply the transformation \eqref{rel1}. The described rules
of composition of new Bailey pairs generate a binary tree of identities for
different products of matrices $M$ and $D$, which are equivalent to
some nontrivial identities for (multiple) elliptic hypergeometric series.

From the relation $M(a,f)M(f,a)=1$ we find the simplest Bailey pairs
$\alpha_n^{(i)}(a,f)=M_{ni}(f,a)$ and $\beta_n^{(i)}(a,f)$ $=\delta_{ni}$.
Let us set $\alpha'(a,t)=D(a;d,e)M(f,a)$ and $\beta'(a,t)=K(f,t,d,e)$,
where $qaf=tde$ and it is assumed that $\alpha'$ and $\beta'$
are the matrices whose columns form the Bailey pairs.
Then the equality $\beta'(a,t)=M(a,t)\alpha'(a,t)$
is equivalent to \eqref{mat-id}. The relation
$\beta''(a,k)=M(a,k)\alpha''(a,k)$ with $\alpha''(a,k)=D(a;b,c)\alpha'(a,t)$ and
$\beta''(a,k)=K(t,k,b,c)\beta'(a,t)$, where $qat=kbc$, leads to the
identity
\begin{eqnarray}\nonumber
&& _{12}V_{11}(a;b,c,d,e,kq^n,q^{-n},af^{-1};q,p)
=\frac{\theta(qa,k/t,qt/b,qt/c)_n}{\theta(k/a,qt,kb/t,kc/t)_n}
\\ && \makebox[8em]{} \times
{_{12}V_{11}}(t;b,c,td/a,te/a,kq^n,q^{-n},tf^{-1};q,p),
\label{bailey}\end{eqnarray}
where $f=kbcde/(aq)^2$. This relation represents an elliptic analogue
of the Bailey transformation for terminating
$_{10}\varphi_9$-series \cite{gas-rah:basic}, which was proved for the
first time in \cite{ft} by a different method.
There exists also a four term Bailey transformation for non-terminating
$_{10}\varphi_9$-series \cite{gas-rah:basic}. Its elliptic generalization
is given by the $V$-function transformation \eqref{E7-1} (written in the
integral form since its infinite series version is not well defined).

Integral analogues of the Bailey chains
were discovered in the paper \cite{spi:bailey2},
where a number of symmetry transformations for the elliptic hypergeometric
integrals has been built with the help of this technique.
The functions $\alpha(z,t)$ and $\beta(z,t)$ form by definition
an integral elliptic Bailey pair with respect to the parameter
 $t$, if they are connected to each other by the relation
\be
\beta(w,t)=\kappa\int_{\T}\eg(tw^{\pm 1}z^{{\pm 1}};p,q)\,\alpha(z,t)\frac{dz}{z}.
\lab{int-bp}\ee
\begin{theorem}\label{bpairs}
Let $\alpha(z,t)$ and $\beta(z,t)$ form an integral elliptic Bailey pair
with respect to the parameter $t$, $|t|<1$. We take the parameters $w, u,s$
satisfying the conditions $w\in \T$ and $|s|,|u|<1, |pq|<|t^2s^2u|$.
Then the functions
\begin{eqnarray*}
&& \alpha'(w,st)=\frac{\eg(tuw^{{\pm 1}};p,q)}
{\eg(ts^2uw^{{\pm 1}};p,q)}\,\alpha(w,t),\\
&& \beta'(w,st)=\kappa \frac{\eg(t^2s^2,t^2suw^{{\pm 1}};p,q)}
{\eg(s^2,t^2,suw^{{\pm 1}};p,q)}
\int_{{\mathbb T}} \frac{\eg(sw^{{\pm 1}}x^{{\pm 1}},ux^{{\pm 1}};p,q)}
{\eg(x^{\pm 2},t^2s^2ux^{{\pm 1}};p,q)}\,\beta(x,t)\frac{dx}{x}
\end{eqnarray*}
define a new integral Bailey pair with respect to
the parameter  $st$, and the functions
\begin{eqnarray*}
&& \alpha'(w,t)=\kappa \frac{\eg(s^2t^2,uw^{{\pm 1}};p,q)}
{\eg(s^2,t^2,w^{\pm 2},t^2s^2uw^{\pm1};p,q)}
\int_{\T}\frac{\eg(t^2sux^{\pm1},sw^{\pm1}x^{\pm1}p,q)}
{\eg(sux^{\pm1};p,q)}\alpha(x,st)\frac{dx}{x},\\
&& \beta'(w,t)=\frac{\eg(tuw^{\pm1};p,q)}{\eg(ts^2uw^{\pm1};p,q)}\beta(w,st)
\end{eqnarray*}
define a new integral Bailey pair with respect to the parameter $t$.
\end{theorem}
In order to prove the first statement, it is sufficient to
substitute relation \eqref{int-bp}
in the definition of $\beta'(w,st)$, to change the order of
integrations, and to use the elliptic beta integral
\eqref{ell-int}. The second statement is proved in an analogous way.

In the paper \cite{spi-war:inversions} it is  shown that
under certain restrictions the integral transformation
$\alpha(x,t)\to \beta(x,t)$ has a very simple inversion.
Let $p,q,t\in\C$ are such that  $|p|, |q|<|t|^2<1$.
For a fixed $w\in\T$, we denote as $C_w$ a contour inside the
annulus $A=\{z\in\C|~|t|-\epsilon<|z|<|t|^{-1}+\epsilon\}$
for some infinitesimally small positive $\epsilon$, such that
the points $t^{-1}w^{\pm1}$ are lying inside $C_w$. Let $f(z,t)$ be
a holomorphic function in $A$ satisfying $f(z,t)=f(z^{-1},t)$.
We define an integral transformation
\begin{equation}\label{g}
g(w,t)=\kappa \int_{C_w} \delta(z,w;t^{-1}) f(z,t)\frac{dz}{z},
\end{equation}
where
\begin{equation*}
\delta(z,w;t^{-1})=\frac{\Gamma(t^{-1}w^{\pm1}z^{\pm1};p,q)}
{\Gamma(t^2,z^{\pm 2};p,q)}.
\end{equation*}
Then for $|t|<|z|<|t|^{-1}$, we have
\begin{equation}\label{f}
f(z,t)=\kappa \int_{\T}\delta(w,z;t) g(w,t)\frac{dw}{w}.
\end{equation}
This relation coincides in essence with the definition of
elliptic integral Bailey pairs. It can be shown that the two ways of building
chains of integral Bailey pairs  indicated above are
related to each other by the described inversion of
integral transformation \eqref{g}.

With the help of Theorem \ref{bpairs} one can build
infinite sequences of Bailey pairs starting from a
given initial pair. The simplest pair can be built
with the help of the elliptic beta integral \eqref{ell-int}.
Each new application of the substitutions indicated above
brings in two new parameters. Equality \eqref{int-bp},
being applied to the appearing new Bailey pairs, leads
to a binary tree of identities for multiple elliptic hypergeometric
integrals with many parameters.

As an illustration, we describe a chain of nontrivial
relations for the integrals
\begin{equation}
 I^{(m)}(t_1,\ldots,t_{2m+6})
=\kappa\int_\T\frac{\prod_{j=1}^{2m+6}\Gamma(t_jz^{\pm 1};p,q)}
{\Gamma(z^{\pm2};p,q)}
\frac{dz}{z}, \quad
 \prod_{j=1}^{2m+6}t_j=(pq)^{m+1},
\label{ell-hyp-int}\end{equation}
where $|t_j|<1$. With the help of the elliptic beta integral
\eqref{ell-int} it is easy to verify validity of the recursion
\ba \label{rec-int}
&& I^{(m+1)}(t_1,\ldots,t_{2m+8})=\frac{\prod_{2m+5\leq k<l\leq 2m+8}
\Gamma(t_kt_l;p,q)}{\Gamma(\rho_m^2;p,q)}
\\ && \makebox[1em]{} \times
\kappa\int_\T\frac{\prod_{k=2m+5}^{2m+8}\Gamma(\rho_m^{-1}t_kw^{\pm 1};p,q)}
{\Gamma(w^{\pm2};p,q)}
I^{(m)}(t_1,\ldots,t_{2m+4},\rho_m w,\rho_m w^{-1})\frac{dw}{w},
\nonumber\ea
where $\rho_m^2=\prod_{k=2m+5}^{2m+8}t_k/pq$.
This equality gives a concrete realization of the Bailey pairs
$\alpha\sim I^{(m)}$ and $\beta\sim I^{(m+1)}$
after a change of parameters. For $m=0$, substitution of the
explicit  expression for  $I^{(0)}$ \eqref{ell-int} to the
right-hand side of \eqref{rec-int} leads to the identity \eqref{E7-1}.
Another important consequence of recursion \eqref{rec-int}
is considered in the next section. In general equality \eqref{rec-int} yields
an $m$-tuple integral representation for $I^{(m)}$ analogous to the
Euler representation for the $_{m+1}F_m$-function.

\section{Biorthogonal functions of the hypergeometric type}

\subsection{Discrete biorthogonal functions with the continuous measure}\

Let us denote
$ \mathcal{E}(\mathbf{t}) :=
{_{12}}V_{11}(t_0^2;t_0t_1,\ldots,t_0t_7;q,p)$
with the balancing condition $\prod_{m=0}^7t_m=q^2$ and the termination
condition $t_0t_k=q^{-n}$ for some $k$. The contiguous with
$\mathcal{E}(\mathbf{t})$ functions, which are obtained by the
change of parameters $t_i$ and $t_k$ to $t_iq$ and $t_kq^{-1}$,
will be denoted as $\mathcal{E}(t_i^+,t_k^-)$. Then, using the
addition formula for theta functions, it is not difficult to
check the equality \cite{spi-zhe:spectral,spi-zhe:rims}
\begin{eqnarray}\nonumber
\lefteqn{\mathcal{E}(\mathbf{t}) - \mathcal{E}(t_6^-, t_7^+) } &&  \\
&&
=\frac{\theta(qt_0^2,q^2t_0^2,qt_7/t_6,t_6t_7/q;p)}
{\theta(qt_0/t_6,q^2t_0/t_6,t_0/t_7,t_7/qt_0;p)}
\prod_{r=1}^5\frac{\theta(t_0t_r;p)}{\theta(qt_0/t_r;p)}\,
\mathcal{E}(t_0^+, t_6^-).
\label{1_con} \end{eqnarray}
Substitution of an elliptic analogue of the Bailey transformation
and of its iterations in \eqref{1_con} allows one to build many
formulae of such a type. One of them has the form
\begin{eqnarray}\nonumber
\lefteqn{\frac{\theta(t_0t_7;p)}{\theta(t_6/qt_0,t_6/q^2t_0,t_6/t_7;p)}
\prod_{r=1}^5 \theta(t_rt_6/q;p)\,
\mathcal{E}(t_0^+, t_6^-)  } && \\ && \makebox[4em]{}
+\frac{\theta(t_0t_6;p)}{\theta(t_7/qt_0,t_7/q^2t_0,t_7/t_6;p)}
\prod_{r=1}^5 \theta(t_rt_7/q;p)\, \mathcal{E}(t_0^+, t_7^-)
\nonumber \\ && \makebox[8em]{}
=\frac{1}{\theta(qt_0^2,q^2t_0^2;p)}
\prod_{r=1}^5\theta(qt_0/t_r;p)\, \mathcal{E}(\mathbf{t}).
\label{2_con} \end{eqnarray}

Let us replace in (\ref{1_con}) the parameter $t_6$ by $qt_6$,
and $t_7$ by $t_7/q$, and substitute $\mathcal{E}(t_0^+,t_7^-)$
from the resulting equality and $\mathcal{E}(t_0^+,t_6^-)$ from (\ref{1_con})
in (\ref{2_con}). This leads to the equation
\begin{eqnarray}\nonumber
\lefteqn{
\frac{\theta(t_0t_7,t_0/t_7,qt_0/t_7;p)}{\theta(qt_7/t_6,t_7/t_6;p)}
\prod_{r=1}^5\theta(q/t_6t_r;p)\left(
\mathcal{E}(t_6^-,t_7^+) - \mathcal{E}(\mathbf{t}) \right) } &&
\\  \nonumber
&& +\frac{\theta(t_0t_6,t_0/t_6,qt_0/t_6;p)}
{\theta(qt_6/t_7,t_6/t_7;p)}\prod_{r=1}^5\theta(q/t_7t_r;p)
\left(\mathcal{E}(t_6^+,t_7^-)-\mathcal{E}(\mathbf{t})\right) \\
&& \makebox[4em]{}
+\theta(q/t_6t_7;p)\prod_{r=1}^5\theta(t_0t_r;p)\,
\mathcal{E}(\mathbf{t})=0.
\label{3_con}\end{eqnarray}
These relations are analogues of equations
 \eqref{cont-1}, \eqref{cont-3}, and  \eqref{eheq} for
the $V$-function, and they can be obtained from them by
the residue analysis for some sequences of poles of the integrand.
Let us replace in \eqref{3_con} the $\mathcal{E}$-function by
\begin{eqnarray}
R_n(x;q,p)= {_{12}}V_{11}\Biggl(\frac{\ve_6}{\ve_8};
\frac{q}{\ve_1\ve_8}, \frac{q}{\ve_2\ve_8},\frac{q}{\ve_3\ve_8},
q^{-n},\frac{Aq^{n-1}}{\ve_8},\frac{\ve_6}{x},\ve_6 x;q,p\Biggr),
\label{R_n} \end{eqnarray}
where $n=0,1,\ldots$ and $A=\ve_1\ve_2\ve_3\ve_6\ve_8$,
and denote $q^{-n}=pq/\ve_4\ve_8$, $Aq^{n-1}/\ve_8=pq/\ve_5\ve_8.$
Then, after the change of notation for parameters
$t_0^2=\ve_6/\ve_8,\ldots, t_0t_7=\ve_6 x$,
we obtain equation \eqref{eheq2} with $\ve_7=\ve_8/q$ and
discrete values of one of the parameters.
The function $R_n(x;q,p)$ defines thus a particular solution of
the elliptic hypergeometric equation obeying the property
$R_n(px;q,p)=R_n(x;q,p)$.
(Notations of the paper \cite{spi:theta2}, where this function was
investigated, pass to ours after the changes
$t_{0,1,2}\to \ve_{1,2,3},\,
t_3\to\ve_6,\,t_4\to\ve_8,\, \mu\to\ve_4\ve_8/pq$ and $A\mu/qt_4\to pq/\ve_5\ve_8$.)

The parameters $\ve_1,\ldots,\ve_6$ enter \eqref{eheq2} symmetrically.
Because of the balancing condition $\prod_{k=1}^8\ve_k=p^2q^2$,
the function $R_n(x;q,p)$ is invariant with respect to
changes $\ve_k\to p\ve_k,\, k=1,\ldots, 5$. This  guarantees
the symmetry of \eqref{R_n} in  $\ve_1,\ldots,\ve_5$, and any
of these parameters can be used for the series termination.
Permutation of one of the parameters $\ve_{1,2,3,5}$ with $\ve_6$
and application of the elliptic Bailey transformation for terminating
${}_{12}V_{11}$-series leads to $R_n(x;q,p)$
up to some multiplier which does not depend on $x$.

We identify now the parameters in  equation
\eqref{3_con} in a different way:
$$
t_0^2=\frac{\ve_6}{\ve_8},\; t_0t_{1,2,3}=\frac{q}{\ve_{1,2,3}\ve_8},\;
t_0t_4=\frac{\ve_6}{x},\; t_0t_5=\ve_6x,\; t_0t_6=q^{-n},\;
t_0t_7=\frac{A}{\ve_8}q^{n-1}.
$$
This leads to a three term recurrence relation in the index $n$:
\begin{eqnarray} \nonumber
&& \makebox[-1em]{}
\frac{\theta\left(\frac{Aq^{n-1}}{\ve_8},\frac{\ve_6q^{2-n}}{A},\frac{\ve_6q^{1-n}}{A};p
\right)} {\theta\left(\frac{Aq^{2n-1}}{\ve_8},\frac{Aq^{2n}}{\ve_8};p\right)}
\prod_{r=1}^3\theta\left(\ve_r\ve_6q^n;p\right)\,
\theta\left(\frac{q^{n+1}x^{\pm 1}}{\ve_8};p\right)\left(R_{n+1}-R_n\right)
\\ &&
+\frac{\theta\left(q^{-n},\frac{\ve_6q^{n}}{\ve_8},\frac{\ve_6q^{1+n}}{\ve_8};p
\right)}
{\theta\left(\frac{q^{1-2n}\ve_8}{A},\frac{q^{2-2n}\ve_8}{A};p\right)}
\prod_{r=1}^3\theta\left(\frac{\ve_r\ve_6\ve_8q^{1-n}}{A};p\right)\,
\theta\left(\frac{q^{2-n}x^{\pm 1}}{A};p\right)\left(R_{n-1}-R_n\right)
\nonumber\\ && \makebox[4em]{}
+\theta\left(\frac{\ve_6q^n}{A},\ve_6x^{\pm 1};p\right)\prod_{r=1}^3
\theta\left(\frac{q}{\ve_r\ve_8};p\right) R_n=0
\lab{3ttr}\end{eqnarray}
with the initial conditions $R_{-1}=0$ and $ R_0=1$.

Let us introduce the functions
$$
z(x)= \frac{\theta(x\xi^{\pm 1};p)}{\theta(x\eta^{\pm 1};p)}, \qquad
\alpha_n= z(q^n/\ve_8),\qquad \beta_n= z(Aq^{n-1}),
$$
where $\xi$ and $\eta$ are arbitrary gauge parameters,
$\xi\neq \eta^{\pm1} p^k,\, k\in\Z$. Then relation
\eqref{3ttr} can be rewritten in a more structured form
\begin{eqnarray}\nonumber
&& \makebox[-2em]{}
(z(x)-\alpha_{n+1})\rho(Aq^{n-1}/\ve_8)\left(R_{n+1}(x;q,p)-R_n(x;q,p)\right)
\\ &&
+(z(x)-\beta_{n-1})\rho(q^{-n})\left(R_{n-1}(x;q,p)-R_n(x;q,p)\right)
\nonumber \\ && \makebox[4em]{}
+\delta (z(x)-z(\ve_6)) R_n(x;q,p)=0,
 \lab{ttr}\end{eqnarray}
where
\begin{eqnarray*}
&& \rho(t)=\frac{\theta\left(t,\frac{\ve_6}{\ve_8t},\frac{q\ve_6}{\ve_8t},
\frac{qt}{\ve_1\ve_2},\frac{qt}{\ve_2\ve_3},\frac{qt}{\ve_1\ve_3},
\frac{q^2t\eta^{\pm 1}}{A};p\right)}
{\theta\left(\frac{qt^2\ve_8}{A},\frac{q^2t^2\ve_8}{A};p\right)},
\\ &&
\delta=\theta\left(\frac{q^2\ve_6}{A},\frac{q}{\ve_1\ve_8},
\frac{q}{\ve_2\ve_8},\frac{q}{\ve_3\ve_8},\ve_6\eta^{\pm 1};p\right).
\nonumber
\end{eqnarray*}
The initial conditions $R_{-1}=0$ and $ R_0=1$ guarantee that $R_n(x;q,p)$
are rational functions of  $z(x)$ with the poles at the points
$\alpha_1,\ldots,\alpha_n$ (i.e., all the dependence
on $x$ enters $R_n$ only through the variable $z(x)$).

Suppose that  $\phi_\lambda$ is a solution of an abstract
generalized eigenvalue problem
$\mathcal{D}_1\phi_\lambda= \lambda \mathcal{D}_2\phi_\lambda$
for some operators $\mathcal{D}_{1,2}$. Let also a scalar product
$\langle \psi|\phi\rangle$ is given, which defines the formal
conjugated operators $\mathcal{D}_{1,2}^T$
by the standard rule $\langle \mathcal{D}_{1,2}^T\psi|\phi\rangle
=\langle \psi|\mathcal{D}_{1,2}\phi\rangle.$
Let $\psi_\lambda$ denote solutions of the dual
generalized eigenvalue problem
 $\mathcal{D}_1^T\psi_\lambda= \lambda \mathcal{D}_2^T\psi_\lambda$.
Then
$0=\langle \psi_\mu|(\mathcal{D}_1-\lambda \mathcal{D}_2)\phi_\lambda\rangle
=(\mu-\lambda)\langle \mathcal{D}_2^T\psi_\mu|\phi_\lambda\rangle$,
that is the function $\mathcal{D}_2^T\psi_\mu$ is orthogonal to $\phi_\lambda$
for $\mu\neq\lambda$. Consequences of this well known fact of the
linear algebra were investigated in detail by Zhedanov \cite{zhe:gevp}
in the case when $\mathcal{D}_{1,2}$ are the Jacobi matrices
(i.e., the tridiagonal matrices). In particular, such generalized
eigenvalue problems were shown to be equivalent to the
theory of biorthogonal rational functions generalizing the orthogonal
polynomials. They are connected also to the recurrence relation
of the $R_{II}$-type investigated in \cite{ism-mas:general} and to the
orthogonality relations appearing within the theory of multipoint Pad\'e
approximation \cite{gon:speed,GL}. Recurrence relation \eqref{ttr} belongs
to this class of problems and, therefore, there exists
a linear functional $\mathcal{L}$ with the condition
$\mathcal{L}\{T_m(x;q,p)R_n(x;q,p)\}=h_n\delta_{nm}$
for some rational functions $T_m$ and normalization constants $h_n$.

The elliptic hypergeometric equation for $R_n$-functions can be rewritten
in the form of a generalized eigenvalues problem  \cite{spi:theta2}:
$$
\mathcal{D}(\ve_4',\ve_5')R_n=\lambda_n\mathcal{D}(\ve_4'',\ve_5'')R_n,
$$
where
$$
\mathcal{D}(\ve_4,\ve_5)=A(x)(T_{q,x}-1)+ A(x^{-1})(T_{q,x}^{-1}-1) + \nu,
\quad T_{q,x}f(x)=f(qx),
$$
denotes the operator permitting to rewrite equation \eqref{eheq2}
as $\mathcal{D}(\ve_4,\ve_5)f(x)=0$. The primed parameters are
arbitrary under the restriction $\ve_4'\ve_5'= \ve_4''\ve_5''=\ve_4\ve_5$,
and other parameters remain untouched (that is why the dependence on them
is not indicated). The spectral variable
$$
\lambda_n=\frac{\theta(\ve_4/\ve_4',\ve_4/\ve_5';p)}
{\theta(\ve_4/\ve_4'',\ve_4/\ve_5'';p)}
$$
is discrete because $\ve_4=pq^{n+1}/\ve_8$.
If we take as the functional $\mathcal{L}$ the integral whose
kernel coincides with the kernel of elliptic beta integral,
$$
\mathcal{L}\{\phi(x)\psi(x)\}=\kappa\int_{\T}
\frac{\prod_{j=1,2,3,6,8}\eg(\ve_j x^{\pm 1};p,q)}
{\eg(x^{\pm2}, Ax^{\pm 1};p,q)}\phi(x)\psi(x)\frac{dx}{x},
$$
then the functions
\be
T_n(x;q,p)=
{}_{12}V_{11}\left(\frac{A\ve_6}{q};\frac{A}{\ve_1},\frac{A}{\ve_2},\frac{A}{\ve_3},
\ve_6x,\frac{\ve_6}{x},q^{-n},\frac{Aq^{n-1}}{\ve_8};q,p\right),
\lab{T_n}\ee
serve as an analogue of $\mathcal{D}_2^T\psi_\mu$ for $R_n(z;q,p)$.
$T_n(x;q,p)$ are the rational functions of $z(x)$ with the poles
at the points $\beta_1,\ldots,\beta_n$ which are obtained
from $R_n(x;q,p)$ after the parameter change
$\ve_8\to pq/A$ (dependence on  $p$
in the parameters disappears because of the ellipticity of
the $_{12}V_{11}$-series in them).

We denote $R_{nm}(x):= R_n(x;q,p)R_m(x;p,q)$ and
$T_{nm}(x):= T_n(x;q,p)T_m(x;p,q)$, where all $_{12}V_{11}$-series
terminate simultaneously because of the modified termination
condition $\ve_4\ve_8=p^{m+1}q^{n+1},\, n,m=0,1,\ldots$.
Since $R_m(qx;p,q)=R_m(x;p,q)$, the functions
$R_{nm}$ represent now solutions of not one, but two generalized
eigenvalue problems differing from each other by
permutation of $p$ and $q$.
Therefore the orthogonality relations in our case
appear to be more complicated than for the biorthogonal rational functions.

\begin{theorem}
The functions $R_{nm}(x)$ and $T_{nm}(x)$
satisfy the following two-index biorthogonality relations:
\be
\kappa\int_{C_{mn,kl}}T_{nl}(x)R_{mk}(x)
\frac{\prod_{j\in S}\eg(\ve_j x^{\pm 1};p,q)}
{\eg(x^{\pm2}, Ax^{\pm 1};p,q)}\frac{dx}{x}
=h_{nl}\: \delta_{mn}\: \delta_{kl},
\lab{2ib}\ee
where $S=\{1,2,3,6,8\}$, $C_{mn,kl}$ denotes a contour separating the
sequences of points
$
\ve_jp^aq^b \,(j=1,2,3,6),\; \ve_8 p^{a-k}q^{b-m}, p^{a+1-l}q^{b+1-n}/A,\;
a,b=0,1,\ldots,
$
from their $x\to x^{-1}$ reciprocals,
and the normalization constants have the form
\begin{eqnarray*}
h_{nl}&=&
\frac{\prod_{j< k,\, j,k\in S} \eg(\ve_j\ve_k;p,q)}
{\prod_{j\in S} \eg(A\ve_j^{-1};p,q)}\,h_n(q,p)\cdot h_l(p,q), \\
h_n(q,p)&=&\frac{\theta(A/q\ve_8;p)
\theta(q,q\ve_6/\ve_8,\ve_1\ve_2,\ve_1\ve_3,\ve_2\ve_3,A\ve_6)_n\,q^{-n}}
{\theta(Aq^{2n}/q\ve_8;p) \theta(1/\ve_6\ve_8,\ve_1\ve_6,\ve_2\ve_6,\ve_3\ve_6,
A/q\ve_6,A/q\ve_8)_n}.
\end{eqnarray*}
\end{theorem}

A direct proof of this statement by a straightforward computation
of the integral on the left-hand side with the help
of formula \eqref{ell-int} and the Frenkel--Turaev sum
is given in \cite{spi:theta2}. Appearance of the two-index orthogonality
relations for univariate functions is a new phenomenon in the
theory of special functions. It is worth of noting that
$R_{nm}(x)$ and $T_{nm}(x)$ are meromorphic functions
of $x\in\C^*$ with essential singularities at
$x=0, \infty$; only for $k=l=0$ or $n=m=0$
they become rational functions of some argument depending on $x$.
For $k=l=0$, one can take the limit $p\to 0$ with fixed
parameters and obtain the functions $R_n(x;q,0)$ and $T_n(x;q,0)$,
which coincide with the family of continuous
 ${}_{10}\varphi_9$ biorthogonal rational functions of Rahman
\cite{rah:integral}. A further degeneration of these functions leads
to the Askey--Wilson polynomials \cite{aw}. Additional restrictions
on one of the parameters in $R_n(x;q,p)$ and $T_n(x;q,p)$
leads to a finite-dimensional systems of biorthogonal rational
functions of a discrete argument \cite{spi-zhe:spectral,spi-zhe:classical},
generalizing the Wilson functions \cite{wil:orthogonal}.
An elementary approach to the analysis of these functions,
related to the elliptic $6j$-symbols \cite{ft},
was suggested by Rosengren in \cite{ros:elementary}.
Some properties of the functions $R_n(x;q,p)$ are investigated
in the recent paper \cite{e-taylor}.

One can build a relation analogous to \eqref{2ib} on the basis
of the modified elliptic beta integral \eqref{circle-int} \cite{spi:thesis}.
For this it is necessary to use the parameterization of base variables
in terms of the quasiperiods $\omega_{1,2,3}$ and pass to the functions
$r_n(u;\omega_1,\omega_2,\omega_3)=R_n(e^{2\pi i u/\omega_2};
e^{2\pi i\omega_1/\omega_2},e^{2\pi i\omega_3/\omega_2})$, where
we have also substituted $\ve_j=e^{2\pi i g_j/\omega_2}$.
Analogously, it is necessary to redenote $T_n(x;q,p)$ as
$s_n(u;\omega_1,\omega_2,\omega_3)$ and $h_n(q,p)$ as
$h_n(\omega_1,\omega_2,\omega_3)$. The products
$r_{nm}(u)= r_n(u;\omega_1,\omega_2,\omega_3)$
$r_m(u;\omega_2,\omega_1,\omega_3)$ and $s_{nm}(u)=
s_n(u;\omega_1,\omega_2,\omega_3)
s_m(u;\omega_2,\omega_1,\omega_3)$ are invariant with respect to the
permutations
 $\omega_1\leftrightarrow\omega_2$, $n\leftrightarrow m$.
Then, for a specially chosen  contour  $\tilde C_{mn,kl}$, we have
\be
\tilde \kappa\int_{\tilde C_{mn,kl}}s_{nl}(u)r_{mk}(u)
\frac{\prod_{j\in S}G(g_j\pm u;\mathbf{\omega})}
{G(\pm2 u, \mathcal{A}\pm u;\mathbf{\omega})}\frac{du}{\omega_2}
=\tilde h_{nl}\: \delta_{mn}\: \delta_{kl},
\lab{2ib-mod}\ee
where $\mathcal{A}=\sum_{j\in S}g_j$ and
$$
\tilde h_{nl}=\frac{\prod_{j<m,\, j,m\in S}G(g_j+g_m;{\bf \omega})}
{\prod_{j\in S}G(\mathcal{A}-g_j;\mathbf{\omega})}
h_n(\omega_1,\omega_2,\omega_3)h_l(\omega_2,\omega_1,\omega_3).
$$
In distinction from the previous case, the limiting transition to the
 $q$-hypergeometric level
$\text{Im}(\omega_3/\omega_2), \text{Im}(\omega_3/\omega_1)\to+\infty$
(that is, $p,r\to 0$) is well defined and preserves the two-index
structure of biorthogonality relations. In particular,
the $r_{nm}(u)$-function degenerates now to the product of two
$q$-hypergeometric series
\begin{eqnarray*}
&&r_{nm}(u;\omega_1,\omega_2)={}_{10}W_9\Big(e^{2\pi i(g_6-g_8)/\omega_2};
e^{2\pi i(\omega_1-g_1-g_8)/\omega_2},\ldots,
e^{2\pi i(g_6+u)/\omega_2};q,q\Big)\\
&&\makebox[4em]{}\times
{}_{10}W_9\Big(e^{2\pi i(g_6-g_8)/\omega_1};
e^{2\pi i(\omega_2-g_1-g_8)/\omega_1},\ldots,
e^{2\pi i(g_6+u)/\omega_1}; \tilde{q}^{-1}, \tilde{q}^{-1}\Big),
\end{eqnarray*}
whose basic variables are related by a modular transformation and the
normalization of the measure is given by the integral \eqref{unit-rah}.

\subsection{A terminating continued fraction}\

A terminating continued
fraction related to the rational functions $R_n(x;q,p)$ is computed
in the paper \cite{spi-zhe:rims}. Let $U_n$ and $V_n$ denote two
sequences of numbers satisfying the three term recurrence relation
\begin{equation}
\psi_{n+1}= \xi_n \psi_n + \eta_n \psi_{n-1}, \quad n=1,2,\ldots
\label{ttr-cf}\end{equation}
with some coefficients $\xi_n$  and $\eta_n$ and the initial conditions
$U_0=0,\; U_1 = 1$ and $V_0 = 1, \; V_1 = \xi_0.$
It is well known that their ratio is related to the finite continued fraction
\begin{equation}
\frac{U_n}{V_n}={1\over\displaystyle \xi_0 + {\strut \eta_1
\over\displaystyle \xi_1 +{\strut  \dots
{\atop \displaystyle +\frac{\eta_{n-1}}{\xi_{n-1}} }}}},
\qquad n=1,2,\ldots.
\label{cf-chain} \end{equation}
Let us define polynomials of $z(x)$ of the $n$-th degree:
$$
P_n(z(x))=\kappa_n \prod_{k=1}^n(z-\alpha_k)\, R_n(x;q,p),
\quad \kappa_n=\prod_{j=0}^{n-1}\rho(aq^{j-1})
$$
and set $P_0(z(x))=1$. Replacing $R_n(x;q,p)$ by $P_n(z(x))$ in \eqref{ttr},
we obtain the following recurrence relation
\begin{equation}
P_{n+1}(z) + (v_n -\rho_n z) P_{n}(z) +
u_n(z-\alpha_n)(z-\beta_{n-1}) P_{n-1}(z) =0
\label{rec_P} \end{equation}
with the initial conditions $P_{-1}=0,\, P_0=1$ and the recurrence
coefficients
\begin{eqnarray}\nonumber
&& u_n=\rho(q^{-n})\rho(Aq^{n-2}/\ve_8), \quad
\rho_n=\rho(Aq^{n-1}/\ve_8)+\rho(q^{-n})-\delta,
\\ &&
v_n=\alpha_{n+1}\rho(Aq^{n-1}/\ve_8)+\beta_{n-1}\rho(q^{-n})-\delta z(\ve_6).
\lab{rc3}\end{eqnarray}
In this case $V_n$ $=P_n(z)$, and $U_n = P^{(1)}_{n-1}(z)$ are the
associated polynomials of the degree $n-1$ in $z$.

Let us suppose that the polynomial $P_{N+1}(z)$ has only simple
zeros, that is $P_{N+1}(z_s)=0$, $z_s\equiv z_s^{(N+1)}, \: s=0,1,\dots,N,$
$z_s\ne z_{s'}$ for $s \ne {s'}$. Then the corresponding continued fraction
 can be expanded into the partial fraction (as a rational function of $z$)
\begin{equation}
\frac{P^{(1)}_{N}(z)}{P_{N+1}(z)}
=\sum_{s=0}^N \frac{g_s}{z-z_s}, \qquad
g_s = \frac{P^{(1)}_N(z_s)}{P_{N+1}'(z_s)}.
\label{par_frac} \end{equation}

The Casoratian of any two solutions $U_n$ and $V_n$ of \eqref{ttr-cf}
satisfies the relation
$$
U_{n+1}V_n-U_nV_{n+1}=(-1)^n\eta_1\cdots\eta_n (U_1V_0-U_0V_1),
$$
which yields
\begin{equation}
P_n(z) P^{(1)}_n(z) - P_{n+1}(z)P^{(1)}_{n-1}(z) = h_n A_n(z) {} B_n(z),
\label{Wron} \end{equation}
where $h_n = u_1 u_2 \cdots u_n$ and
$$
A_n(z)=\prod_{i=1}^n(z-\alpha_i), \qquad
{} B_n(z)=\prod_{i=1}^n(z-\beta_{i-1}).
$$
Fixing  $n=N$ and setting $z=z_s, \; s=0,1,\dots,N,$ in \eqref{Wron},
we can express $P^{(1)}_N(z_s)$ in terms of $P_N(z_s),$ $h_N,$  $A_N(z_s)$ and
${}B_N(z_s)$. This leads to the following convenient for computations
expression for the residues of the poles  $g_s$ in \eqref{par_frac}:
\begin{equation}
g_s = \frac{h_N A_N(z_s){} B_N(z_s)}{P_{N+1}'(z_s) P_N(z_s)}.
\label{ex_g} \end{equation}

In the limit
\begin{equation}
\ve_3\ve_6=q^{-N+\epsilon}, \quad N=0,1,\ldots,\quad \epsilon\to 0,
\label{obryv}\end{equation}
we find that  $u_{N+1}\to 0$, and the continued fraction terminates
automatically. It appears that this fraction can be computed in the
closed form using formula \eqref{par_frac}.
For $\epsilon\to 0$, the rational function $R_{N+1}(x;q,p)$ diverges,
since the elliptic Pochhammer symbol $\theta(\ve_3\ve_6)_{N+1}\to 0$
in the denominator of the last term of the $_{12}V_{11}$-series.
However, we have simultaneously $\kappa_{N+1}\to 0$, so that the polynomial
$P_{N+1}(z)$ takes the finite value, and its zeros are found explicitly:
$z_s=z(\ve_6q^s),\, s=0,\ldots, N.$ It appears also that the polynomial
 $P_N(z_s)$ is computable in the closed form owing to the Frenkel--Turaev
summation formula. The other quantities defining the residues
$g_s$ are found sufficiently easily, although they are given by
rather cumbersome expressions.

Suppose that the conditions of the simplicity of zeros
$z_s$ are satisfied as well as other restrictions on parameters
guaranteeing that $u_k\neq 0$ for $k=1,\ldots, N,$
the descriptions of which we skip. Then the terminating
elliptic hypergeometric continued fraction has the following
explicit representation:
\begin{eqnarray}\lab{cf-fin}
&& \makebox[-1em]{}
{1\over\displaystyle \rho_0 z-v_0 -
{\strut u_1(z-\alpha_1)(z-\beta_0) \over\displaystyle \rho_1  z-v_1 -
\dots {\atop \strut \displaystyle -\frac{u_N(z-\alpha_N)(z-\beta_{N-1})}
{\rho_N z-v_N} }}}
\\ && \makebox[1em]{}
=\frac{1}{(z(\ve_6)-z(x))\delta}\,
{}_{12}V_{11}\left(\frac{q\ve_6}{\ve_8}; q,\ve_6\ve_1,\ve_6\ve_2,\frac{qx}{\ve_8},
\frac{q}{\ve_8 x},q^{-N},\frac{\ve_6 q^{N+2}}{\ve_1\ve_2\ve_8};q,p\right),
\nonumber\end{eqnarray}
where in the expressions for all recurrence coefficients,
including $\delta$, it is necessary to substitute
 $\ve_3=q^{-N}/\ve_6$ and $z(x)=\theta(x\xi^{\pm 1};p)/\theta(x\eta^{\pm 1};p)$.
In \cite{spi-zhe:rims} this result was presented in the different (additive)
system of notation.

This formula describes the most general terminating continued fraction
of the hypergeometric type, which was found to the present moment.
For fixed values of parameters, in the limit $p\to0$ one obtains
the terminating continued fraction of Gupta and Masson
\cite{gup-mas:contiguous} (see Corollary 3.3) described by
a very-well poised balanced $_{10}\varphi_9$-series.
Further specification of parameters leads to the continued fraction
of Watson which, in its turn, is a $q$-analogue of the famous
Ramanujan continued fraction (see the details in \cite{gup-mas:watson}).

\subsection{Continuous biorthogonality of the $V$-function}\

The relations described above \eqref{2ib} correspond to discrete
values of one of the parameters in the elliptic hypergeometric
equation. In paper \cite{spi:aa2}, it was shown that the
$V$-function with general set of continuous parameters obeys also
some biorthogonality relations characteristic to the continuous spectra.

Let us consider the $m=2$ case in recursion \eqref{rec-int}.
After imposing the constraints on parameters
 $t_5t_7=t_6t_{8}=pq$ the integral on the left-hand side
becomes explicitly computable. Then, after a number of notational changes
and application of transformation \eqref{E7-1}, there appears the equality
\be
\phi(x;c,d|\xi;s)=\kappa\int_\T R(c,d,a,b;x,w|s)\phi(w;a,b|\xi;s)\frac{dw}{w},
\lab{key-rel-fin}\ee
where the basis vectors have the form
\be
\phi(w;a,b|\xi;s)=\Gamma(sa\xi^{\pm 1},sb\xi^{\pm 1},
\sqrt{\frac{pq}{ab}}w^{\pm 1}\xi^{\pm 1};p,q)
\lab{phi}\end{equation}
and
\begin{eqnarray}\nonumber
&& R(c,d,a,b;x,w|s)=
\frac{1}{\Gamma(\frac{pq}{ab},\frac{ab}{pq},w^{\pm2};p,q)}
\\ && \makebox[4em]{} \times
V\left(s c,s d,\sqrt{\frac{pq}{cd}}x,
\ve\sqrt{\frac{pq}{cd}}x^{-1},\frac{pq}{as},\frac{pq}{bs},
\sqrt{\frac{ab}{pq}}w,\sqrt{\frac{ab}{pq}}w^{-1}\right).
\label{R}\end{eqnarray}
Here the parameters $a,b,c,d, s$ and $x,\xi$ are arbitrary,
but their choice should match with the condition that
taken contours of integration separate converging
to zero and diverging to infinity sequences of poles
of the integrands. The variable $\xi$ enters only the
basis vectors $\phi(w;a,b|\xi;s)$ and the kernel $R$
can be considered as a ``rotation matrix"
with continuous indices $x$ and $w$, which
permits to change arbitrarily the parameters
$a$ and $b$. Equality \eqref{R} represents an integral generalization
of the relation which was used by Rosengren in \cite{ros:elementary}
for derivation of properties of the elliptic $6j$-symbols.

Denoting $w=e^{i\theta}$,  $y=\cos\theta$, and using the equality
$$
\int_\T f(\cos\theta)\frac{dw}{iw}=\int_0^{2\pi}f(\cos\theta)d\theta
=2\int_{-1}^1f(y)\frac{dy}{\sqrt{1-y^2}},
$$
we can write
\be
\phi(x;c,d|\xi;s)=2i\kappa
\int_{-1}^1 R(c,d,a,b;x,e^{i\theta}|s)
\phi(e^{i\theta};a,b|\xi;s) \frac{dy}{\sqrt{1-y^2}}.
\lab{key-real}\ee
In the limit $c\to a$ and $d\to b$, there appears the following
relation in the distributional sense
\be
\lim_{c\to a,d\to b}
R(c,d,a,b;e^{i\varphi},e^{i\theta})
=\frac{2\pi\sqrt{1-y^2}}{(p;p)_\infty(q;q)_\infty}\,\delta(v-y),
\lab{phi-ort}\ee
where $v=\cos\varphi$.

A double application of relation
\eqref{key-rel-fin} with different parameters leads in an evident way
to the self-reproducing property for the kernel
\begin{equation}
\kappa\int_\T R(a,b,c,d;x,w|s)R(c,d,e,f;w,z|s)\frac{dw}{w}
=R(a,b,e,f;x,z|s)
\label{repr}\end{equation}
and the biorthogonality relation
\begin{equation}\label{biort}
\int_{-1}^1 R(a,b,c,d;e^{i\varphi},e^{i\theta}|s)
R(c,d,a,b;e^{i\theta},e^{i\varphi'}|s)
\frac{dy}{\sqrt{1-y^2}}=\frac{\sqrt{1-v^2}}{(2i\kappa)^2}\,\delta(v-v'),
\end{equation}
where $v=\cos\varphi$ and $v'=\cos\varphi'$.
Substitution of the expression for $R$-function \eqref{R}
in \eqref{biort} results in the equality
\begin{eqnarray}\nonumber
&&
\int_{-1}^1 \frac{1}{\Gamma(e^{\pm 2i\theta};p,q)}
V\left(sa,sb,\sqrt{\frac{pq}{ab}}e^{i\varphi},
\sqrt{\frac{pq}{ab}}e^{-i\varphi},\frac{pq}{c},\frac{pq}{d},
\sqrt{\frac{cd}{\rho}}e^{i\theta},\sqrt{\frac{cd}{\rho}}e^{-i\theta}\right)
\\ && \makebox[2em]{} \times
V\left(sc,sd,\sqrt{\frac{pq}{cd}}e^{i\theta},
\sqrt{\frac{pq}{cd}}e^{-i\theta},\frac{pq}{as},\frac{pq}{bs},
\sqrt{\frac{ab}{pq}}e^{i\varphi'},\sqrt{\frac{ab}{pq}}e^{-i\varphi'}\right)
\frac{dy}{\sqrt{1-y^2}}
\nonumber \\ && \makebox[4em]{}
=\Gamma\left(\frac{ab}{pq},\frac{pq}{ab},
\frac{cd}{pq},\frac{pq}{cd},e^{\pm2i\varphi};p,q\right)
\frac{\sqrt{1-v^2}}{(2i\kappa)^2}\,\delta(v-v').
\label{V-ort}\end{eqnarray}
The parameters $v$ and $v'$ can be considered as continuous
spectral variables in the operator formulation of the elliptic
hypergeometric equation. Therefore relation \eqref{V-ort}
should follow from the latter equation, but the precise
connection between them is not established yet.

\section{Connection with the Sklyanin algebra}

In \cite{rai:trans}, Rains introduced an interesting
finite difference operator connected with the root system
$BC_n$. For $n=1$, it can be represented in the form
\begin{equation}
D(a,b,c,d;p;q)=
\frac{\theta(az,bz,cz,dz;p)}{z\theta(z^2;p)}T_{z,q}^{1/2}
+\frac{\theta(az^{-1},bz^{-1},cz^{-1},dz^{-1};p)}
{z^{-1}\theta(z^{-2};p)}T_{z,q}^{-1/2},
\label{D-oper}\end{equation}
where $T_{z,q}^{\pm1/2}f(z)=f(q^{\pm1/2}z)$ is the $q$-shift operator
and $a,b,c,d$ are arbitrary parameters. Later on Rains noticed also
\cite{rai:abelian,ros:sklyanin} that this operator is equivalent to
the general linear combination of four generators of the Sklyanin algebra
$S_0,\ldots,S_3$ \cite{skl:alg1,skl:alg2}.

Defining relations of the Sklyanin algebra have the form
\begin{eqnarray}\nonumber
&& S_\alpha S_\beta -S_\beta S_\alpha =i(S_0S_\gamma+S_\gamma S_0),
\\ &&
S_0S_\alpha - S_\alpha S_0 = iJ_{\beta\gamma}(S_\beta S_\gamma+S_\gamma S_\beta),
\label{s-rel2}\end{eqnarray}
where $J_{\beta\gamma}$ are the structure constants of the algebra
and $(\alpha,\beta,\gamma)$ is an arbitrary cyclic permutation of
the triple $(1,2,3)$.
A representation of $S_a$ as finite difference operators has been found
in  \cite{skl:alg2}:
$$
S_a=i^{\delta_{a,2}}\frac{\theta_{a+1}(\eta|\tau)}{\theta_1(2u|\tau)}
\left(\theta_{a+1}(2u-2g|\tau)e^{\eta\partial_u}
-\theta_{a+1}(-2u-2g|\tau)e^{-\eta\partial_u}\right),
$$
where $e^{\pm\eta \partial_u}f(u)=f(u\pm\eta)$,
and under the quantization condition $g=\ell \eta$, $\ell=0,1/2,1,\ldots,$
their action in the space of theta functions of the order $4\ell$
was described. The combination of the generators
\begin{eqnarray*}
&& 2\Delta(a_1,a_2,a_3,a_4):=\frac{\prod_{j=1}^3\theta_1(a_j+a_4+2g)}
{\theta_1(\eta)}S_0
-\frac{\prod_{j=1}^3\theta_1(a_j+a_4+2g+\frac{1}{2})}
{\theta_1(\eta+\frac{1}{2})}S_1
\\ && \qquad
-ie^{\pi i(\frac{\tau}{2}+2a_4+2g-\eta)}
\frac{\prod_{j=1}^3\theta_1(a_j+a_4+2g+\frac{1+\tau}{2})}
{\theta_1(\eta+\frac{1+\tau}{2})}S_2
\\ && \qquad
+e^{\pi i(\frac{\tau}{2}+2a_4+2g-\eta)}
\frac{\prod_{j=1}^3\theta_1(a_j+a_4+2g+\frac{\tau}{2})}{\theta_1(\eta+\frac{\tau}{2})}
S_3,
\end{eqnarray*}
with the normalization $\sum_{j=1}^4a_j=-4g$, can be represented
in the form \cite{ros:sklyanin}
$$
\Delta(a_1,a_2,a_3,a_4)=\frac{\prod_{j=1}^4\theta_1(a_j+u)}
{\theta_1(2u)}e^{\eta\partial_u}+\frac{\prod_{j=1}^4\theta_1(a_j-u)}
{\theta_1(-2u)}e^{-\eta\partial_u}.
$$
After the transition to multiplicative system of notation
$$
(a,b,c,d):=e^{2\pi ia_{1,2,3,4}}, \quad \rho:=abcd= e^{-8\pi i g},\quad
z:=e^{2\pi i u},\quad q:=e^{4\pi i\eta}
$$
there appears the operator described above \eqref{D-oper}:
$$
\Delta(a_1,a_2,a_3,a_4)=\left(ip^{1/8}(p;p)_\infty\right)^3e^{4\pi i g}
D(a,b,c,d;p;q).
$$

The standard eigenvalue problem $D\psi=\lambda\psi$
appears to be very complicated, since it represents a difference
analogue of the Heun equation \cite{spi:aa2}. On the one hand, it is
known that the eigenvalue problem for the one particle Hamiltonian
of the Inozemtsev model \cite{ino:lax} is equivalent to the Heun equation.
On the other hand, the classical equations of motion for this Hamiltonian
with the modular parameter $\tau$ considered as a time variable
leads to the Painlev\'e VI equation \cite{man:painleve}.
The $D$-operator itself appears to be related to the van Diejen model
\cite{die:integrability,spi:aa2}. All this and the connection with
elliptic hypergeometric functions described below
demonstrate  some mathematical universality of the operator \eqref{D-oper}.

Consider the generalized eigenvalue problem of the form
\be
D(a,b,c,d;p;q)f(z;w;q^{1/2}a,q^{1/2}b;\rho)=
\lambda(w) D(a,b,c',d';p;q)f(z;w;q^{1/2}a,q^{1/2}b;\rho),
\label{gevp}\ee
where $cd=c'd'$. This equation is solved explicitly. Using the parameterization
$$
\lambda(w)=\frac{\theta(w\sqrt{c/d},w\sqrt{d/c};p)}
{\theta(w\sqrt{c'/d'},w\sqrt{d'/c'};p)},
$$
for $|q|<1$  we obtain \cite{spi:aa2}
\begin{equation}
f(z;w;a,b;\rho)=\Gamma\left(\frac{pq}{a} z^{\pm 1},
\frac{pq}{b} z^{\pm 1}, \sqrt{\frac{ab}{\rho}}w^{\pm 1} z^{\pm 1};p,q\right),
\label{basis}\end{equation}
up to the multiplication by an arbitrary function
$\varphi(z)$, $\varphi(qz)=\varphi(z)$. Evidently, function \eqref{basis}
coincides with the basis vector \eqref{phi} after a change of parameters.

For $z\to z^{-1}$ invariant functions, $\psi(z)=\psi(z^{-1})$,
we define the scalar product
\begin{equation}
\langle \chi(z),\psi(z)\rangle =\kappa\int_{\T}\frac{\chi(z)\psi(z)}
{\Gamma(z^{\pm 2};p,q)} \frac{dz}{z}.
\label{sp}\end{equation}
Then the formally conjugated to $D$ operator has the form
$$
D^*(a,b,c,d;p;q)=\frac{cd}{q^{1/2}}D\left(\frac{pq^{1/2}}{a},\frac{pq^{1/2}}{b},
\frac{q^{1/2}}{c}, \frac{q^{1/2}}{d};p;q\right).
$$
The dual problem
\be
D^*(a,b,c,d;p;q)g(z;v;a,b;\rho)=\lambda(v) D^*(a,b,c',d';p;q)g(z;v;a,b;\rho)
\lab{dual-gevp}\ee
has a solution
\begin{equation}
g(z;v;a,b;\rho)= \Gamma\left(a z^{\pm 1},
b z^{\pm 1}, \sqrt{\frac{\rho}{ab}}v^{\pm 1} z^{\pm 1}; p,q\right),
\label{basis-dual}\end{equation}
which is also defined up to the multiplication by an arbitrary function
$\varphi(z)$, $\varphi(qz)=\varphi(z)$. Now it is not difficult to see that
the scalar product of functions
\eqref{basis} and \eqref{basis-dual} leads to the $V$-functions:
\begin{eqnarray} \nonumber \makebox[-1em]{}
&& \langle g(z;v;a,b;e),f(z;w;c,d;e)\rangle
\\ && \makebox[2em]{}
 = V\left(a,b,\sqrt{\frac{e}{ab}}v,
\sqrt{\frac{e}{ab}}v^{- 1}\frac{pq}{c},\frac{pq}{d},
\sqrt{\frac{cd}{e}}w,\sqrt{\frac{cd}{e}}w^{- 1};p,q
\right).
\label{V-fact}\end{eqnarray}

Thus, the elliptic analogue of the Gauss hypergeometric function
appears to be directly related to the generalized eigenvalue problem
for a linear combination of the Sklyanin algebra generators.
It is necessary to note that our scalar product \eqref{sp}
is different from the Sklyanin invariant measure \cite{skl:alg2}.
In \cite{ros:sklyanin}, using the latter measure Rosengren has
built an integral representation
for the elliptic $6j$-symbols and proved the Sklyanin conjecture
on the reproducing kernel for representations in the space of
theta functions.

The non-uniqueness in the choices of functions
$f$ and $g$ can be fixed by the requirement that these functions
satisfy simultaneously the equations obtained from \eqref{gevp}
and \eqref{dual-gevp} by the permutation of $p$ and $q$
(because the equations $\varphi(qz)=\varphi(pz)=\varphi(z)$
lead to $\varphi(z)=const$). This means that we introduce
into consideration a second copy of the Sklyanin algebra,
obtained from the first one by the permutation of $\tau$ and $2\eta$:
$$
\tilde S_a=i^{\delta_{a,2}}\frac{\theta_{a+1}(\frac{\tau}{2}|2\eta)}
{\theta_1(2u|2\eta)}\left(\theta_{a+1}(2u-2g|2\eta)e^{\frac{\tau}{2}\partial_u}
-\theta_{a+1}(-2u-2g|2\eta)e^{-\frac{\tau}{2}\partial_u}\right).
$$
For these two algebras the following cross-commutation relations are valid:
\begin{eqnarray*}
&& S_a \tilde S_b=\tilde S_b S_a, \quad a,b\in\{0,3\} \quad
\text{or} \quad a,b\in\{1,2\},
\\
&& S_a \tilde S_b=-\tilde S_b S_a, \quad a\in\{0,3\},\; b\in\{1,2\} \quad
\text{or} \quad a\in\{1,2\},\;b\in\{0,3\}.
\end{eqnarray*}

One can substitute in  equations \eqref{gevp} and
\eqref{dual-gevp} the parameterization $z=e^{2\pi i u/\omega_2}$,
$2\eta=\omega_1/\omega_2$,
$\tau=\omega_3/\omega_2$, and to build their solutions
well defined for $|q|=1$
with the help of the modified elliptic
gamma function $G(u;\mathbf{\omega})$. In this case the uniqueness of
solutions can be reached by the requirement that they satisfy
simultaneously to equations obtained from the original ones by the
permutation of $\omega_1$ and $\omega_2$.
This leads to another copy of the Sklyanin algebra, which is
obtained from the first one by the transformations
 $\eta\to 1/(4\eta), u\to u/(2\eta), \tau\to \tau/(2\eta)$:
$$
\tilde S_a=i^{\delta_{a,2}}
\frac{\theta_{a+1}\left(\frac{1}{4\eta}\Big|\frac{\tau}{2\eta}\right)}
{\theta_1\left(\frac{u}{\eta}\Big|\frac{\tau}{2\eta}\right)}
\left(\theta_{a+1}\left(\frac{u-g}{\eta}\Big|\frac{\tau}{2\eta}\right)
e^{\frac{1}{2}\partial_u} -\theta_{a+1}\left(\frac{-u-g}{\eta}\Big|
\frac{\tau}{2\eta}\right) e^{-\frac{1}{2}\partial_u}
\right).
$$
In this case some of the generators $S_a$ and $\tilde S_a$
anticommute with each other as well.

The described direct products of the Sklyanin algebra pairs
can be considered as elliptic analogues of the Faddeev modular double
$U_q(sl_2)\otimes U_{{\tilde q}^{-1}}(sl_2)$
\cite{fad:discrete1,fad:mod}. The latter double
can be obtained from the second  case in the limit Im$(\tau)\to+\infty$
(in the first case this limit is not defined) \cite{spi:aa2}.

\section{Partial fraction decompositions and determinants}

Expansions of different rational functions defined as ratios of
two polynomials into partial fractions are
used in the proofs of many identities for plain and $q$-hypergeometric
series and integrals. At the elliptic level these rational functions are
replaced by ratios of products of theta functions, and one searches for
their expansions into sums of ratios of theta functions with
the minimal number of poles. If the partial fraction expansion
for an arbitrary rational function is a standard procedure, it is not so
in the theta functions case. The first known relation of such a type
follows from an identity given in \cite{whi-wat:course} as an exercise.

\begin{theorem}
Let  $2n$ variables $a_1,\ldots,a_n,b_1,\ldots,b_n\in\C^*$
satisfy the constraint $\prod_{k=1}^na_k=\prod_{k=1}^nb_k$
and $a_j/a_k\neq p^k,\, k\in\Z,$ for $j\neq k$. Then the following
relation for theta functions is true:
\begin{equation}
\sum_{k=1}^n\frac{\prod_{j=1}^n \theta(a_k/b_j;p)}{\prod_{j=1,\, \neq k}^n
\theta(a_k/a_j;p)}=0.
\label{theta-id1}\end{equation}
\end{theorem}
\begin{proof}
We replace in \eqref{theta-id1} $n$ by $n+1$, denote
$a_{n+1}=t$ and substitute $b_{n+1}=a_1\ldots a_{n}t/b_1\ldots b_n$.
After taking out of the sum the $(n+1)$-st term, this relation
can be rewritten in the form \cite{ros:elliptic}
\begin{equation}
\prod_{k=1}^n\frac{\theta(t/b_k;p)}{\theta(t/a_k;p)}
=\sum_{r=1}^n\frac{\theta(ta_1\cdots a_n/a_rb_1\cdots b_n;p)}
{\theta(t/a_r,a_1\cdots a_n/b_1\cdots b_n;p)}
\frac{\prod_{j=1}^n\theta(a_r/b_j;p)}{\prod_{j=1,\, \neq r}^n
\theta(a_r/a_j;p)}
\label{pf1}\end{equation}
and interpreted as a partial fraction expansion over theta functions.
Then the proof of this identity is rather elementary.
For $n=2$ it is reduced to the addition formula \eqref{ident}.
By induction it follows that the left-hand side can be decomposed
into the sum
$$
\sum_{r=1}^n c_r\frac{\theta(ta_1\cdots a_n/a_rb_1\cdots b_n;p)}{\theta(t/a_r;p)},
$$
where the coefficients $c_r$ are easily found after
the multiplication by $\theta(t/a_r;p)$ and the choice $t=a_r$.
\end{proof}

\begin{theorem} \cite{gus87}
Let $2n$ variables $a_1,\ldots,a_n,b_1,\ldots,b_n\in\C^*$
satisfy the relations $a_ja_k,a_j/a_k\neq p^k,\, k\in\Z,$ for $j\neq k$.
Then the following identity for theta functions is true
\begin{equation}
\sum_{k=1}^n\frac{a_k\prod_{j=1}^{n-2} \theta(a_kb_j^{\pm 1};p)}
{\prod_{j=1,\neq k}^n \theta(a_ka_j^{\pm1};p)}=0.
\label{pf2}\end{equation}
\end{theorem}
\begin{proof}
After the replacements $n\to n+1$ and $a_{n+1}\to t$ in
\eqref{pf2} and singling the $(n+1)$-st term out of
the sum, we obtain the partial fraction expansion of the form
$$
\frac{\prod_{j=1}^{n-1}\theta(tb_j^{\pm1};p)}
{\prod_{j=1}^{n}\theta(ta_j^{\pm1};p)}=
\sum_{k=1}^n\frac{\prod_{j=1}^{n-1}\theta(a_kb_j^{\pm1};p)}
{\theta(ta_k^{\pm1};p)\prod_{j=1,\neq k}^{n}\theta(a_ka_j^{\pm1};p)},
$$
which is easily proved by induction.
\end{proof}

These expansions into ``simple" fractions for theta functions
were used in the papers
\cite{die-spi:selberg,spi-war:inversions,ros:elliptic} for the
proof of some exact summation and integration formulae
for elliptic hypergeometric functions. Let us describe also another
expansion which was used recently in \cite{rs:det}:
$$
\frac{\prod_{j=1}^{n+2}\theta(v_jz;p)}{z\theta(z^2;p)\prod_{i=1}^n\theta(u_iz^{-1};p)}
+(z\to z^{-1})=\sum_{i=1}^n\frac{\prod_{j=1}^{n+2}\theta(u_iv_j;p)}
{u_i\theta(u_iz^{\pm1};p)\prod_{k=1,\neq i}^n\theta(u_ku_i^{-1};p)},
$$
where $\prod_{i=1}^n u_i\prod_{j=1}^{n+2}v_j=p^{n-1}.$

An elliptic analogue of the Cauchy determinant has the form
\begin{eqnarray}\lab{d-cau}
&& \det_{1\le i,j\le n} \left(\frac{1}{a_i^{-1}\theta(a_iz_j^{\pm 1};p)}\right)
\\ && \makebox[2em]{}
= (-1)^{n(n-1)/2}\frac{
\prod_{1\le i<j\le n} a_i^{-1}\theta(a_i a_j^{\pm 1};p)
\prod_{1\le i<j\le n} z_i^{-1}\theta(z_i z_j^{\pm 1};p)}
{\prod_{1\le i,j\le n} a_i^{-1}\theta(a_i z_j^{\pm 1};p)}.
\nonumber\end{eqnarray}

The Frobenius determinant has the form
\begin{equation}
\det_{1\le i,j\le n} \left(\frac{\theta(ta_ib_j;p)}
{\theta(t,a_ib_j;p)}\right)
= \frac{\theta(t\prod_{i=1}^na_ib_i;p)}{\theta(t;p)}
\frac{\prod_{1\le i<j\le n} a_jb_j\theta(a_i/a_j,b_i/b_j;p)}
{\prod_{1\le i,j\le n} \theta(a_ib_j;p)}.
\label{frob}\end{equation}

In \cite{war:summation}, Warnaar suggested a new determinant for theta functions
\begin{eqnarray}\label{e-kratt}
&& \det_{1 \leq i,j \leq n} \left(
\frac{ \theta(ax_i,ac/x_i)_{n-j} }
{ \theta(bx_i,bc/x_i)_{n-j} } \right)
\\ && \makebox[1em]{}
= a^{\binom{n}{2}}q^{\binom{n}{3}}
\prod_{1\leq i<j\leq n}x_j\theta(x_ix_j^{-1},cx_i^{-1}x_j^{-1};p)
\prod_{i=1}^n \frac{\theta(b/a,abcq^{2n-2i})_{i-1}}
{\theta(bx_i,bc/x_i)_{n-1}},
\nonumber \end{eqnarray}
where $\theta(a)_n=\prod_{j=0}^{n-1}\theta(aq^j;p)$.
For $p\to 0$, it reduces to the Krattenthaler determinant \cite{kra:major}.

Formulae \eqref{d-cau}--\eqref{e-kratt} are used in
\cite{spi:theta2,rai:trans,war:summation,kaj-nou,ros-sch1,rai:rec}
and some other papers as auxiliary tools for proving necessary
elliptic hypergeometric identities.
Partial fraction decompositions and determinants are somewhat equivalent to
each other. For instance, if one expands determinant \eqref{frob}
along the last row and evaluates each term by the same
formula \eqref{frob} in the smaller dimension $n-1$, then there appears an identity
equivalent to  \eqref{pf1}. Therefore formula \eqref{frob}
follows from  \eqref{pf1} by induction on $n$, and vice versa.
In a similar way, formula \eqref{d-cau} is equivalent to \eqref{pf2}  \cite{rai:rec}.
Applications of determinants at the level of
$q$-hypergeometric functions are described, for example, in \cite{tv}.
In the paper \cite{ros-sch3}, Rosengren and Schlosser have
systematically considered determinants of theta functions on root systems
(in particular, this paper contains a detailed list of references on
this subject) and constructed a number of new exactly computable cases,
which we skip for brevity.
For applications of elliptic determinants to some problems of the number
theory, combinatorics, and statistical mechanics, see \cite{ros:num,ros:IK}.

\section{The elliptic beta integrals on root systems}

\subsection{Integrals for the root system $C_n$}\

There are two different generalizations of the elliptic beta integral
\eqref{ell-int} to multiple integrals for the root system
$C_n$ (or $BC_n$), suggested by van Diejen and the author
\cite{die-spi:elliptic,die-spi:selberg}. We describe first the
multiparameter integral of type I.

\begin{theorem}\label{C-I}
Let $z_1,\ldots,z_n\in\T$ and complex parameters
$t_1,\ldots,t_{2n+4}$ and $p,q$ satisfy the constraints
$|p|, |q|, |t_j|<1$ and $\prod_{j=1}^{2n+4}t_j=pq$. Then
\begin{eqnarray}\nonumber
&& \kappa_n\int_{\T^n}\prod_{1\leq j<k\leq n}\frac{1}{\eg(z_j^{\pm 1} z_k^{\pm 1};p,q)}
\prod_{j=1}^n\frac{\prod_{m=1}^{2n+4}\eg(t_mz_j^{\pm 1};p,q)}
{\eg(z_j^{\pm2};p,q)}\frac{dz_1}{z_1}\cdots\frac{dz_n}{z_n}
\\ && \makebox[4em]{}
=\prod_{1\leq m<s\leq 2n+4}\eg(t_mt_s;p,q), \qquad
\kappa_{n}=\frac{(p;p)_\infty^n(q;q)_\infty^n}{(4\pi i)^n n!}.
\label{C-typeI}\end{eqnarray}
\end{theorem}
\begin{proof}
We consider the function
\begin{eqnarray}\nonumber
&& \rho(z,t;C_n)=\prod_{1\leq i<j\leq n}\frac{1}{\eg(z_i^{\pm 1} z_j^{\pm 1};p,q)}
\prod_{i=1}^n\frac{\prod_{m=1}^{2n+3}\eg(t_mz_i^{\pm 1};p,q)}
{\eg(z_i^{\pm 2},Az_i^{\pm 1};p,q)}
\\ && \makebox[6em]{} \times
\frac{\prod_{m=1}^{2n+3}
\eg(At_m^{-1};p,q)}{\prod_{1\leq m<s\leq 2n+3}\eg(t_mt_s;p,q)},
\label{kernel-C}\end{eqnarray}
where $A=\prod_{m=1}^{2n+3}t_m$. For all $z_i$ there are poles
of \eqref{kernel-C} in the points
$$
\PP=\{t_mq^ap^b,\, A^{-1}q^{a+1}p^{b+1}
\}_{m=1,\ldots, 2n+3,\, a,b=0,1,\ldots}
$$
converging to zero. The coordinates
of the poles going to infinity form the set $\PP^{-1}$.
Then the statement of the theorem may be rewritten in the form
\begin{equation}\label{ell-int-C}
\int_{C^n}\rho(z,t;C_n)\frac{dz}{z}=\frac{(4\pi i)^nn!}
{(q;q)_\infty^n(p;p)_\infty^n},\qquad \frac{dz}{z}:=\prod_{j=1}^n\frac{dz_j}{z_j},
\end{equation}
where the contour $C$ in an arbitrary deformation of $\T$ separating $\PP$ and
$\PP^{-1}$.

The integral kernel $\rho(z,t;C_n)$ satisfies the equation
analogous to (\ref{eqn}):
\begin{eqnarray} \nonumber
&& \rho(z,qt_1,t_2,\ldots,t_{2n+3};C_n)-\rho(z,t;C_n)
\\ && \makebox[4em]{}
=\sum_{i=1}^n\left(g_i(z_1,\ldots,q^{-1}z_i,\ldots,z_n,t)-g_i(z,t)\right),
\label{eqn-C}\end{eqnarray}
where
\begin{equation}
g_i(z,t)=\rho(z,t;C_n)\prod_{j=1,\neq i}^n\frac{\theta(t_1z_j^{\pm 1};p)}
{\theta(z_iz_j^{\pm 1};p)}\frac{\prod_{m=1}^{2n+3}\theta(t_mz_i;p)}
{\prod_{m=2}^{2n+3}\theta(t_1t_m;p)}
\frac{\theta(t_1A;p)}{\theta(z_i^2,Az_i;p)}\frac{t_1}{z_i}.
\label{g-C}\end{equation}
Dividing equation (\ref{eqn-C}) by $\rho(z,t;C_n)$, we obtain
\begin{eqnarray}\nonumber
\lefteqn{\prod_{i=1}^n\frac{\theta(t_1z_i^{\pm 1};p)}{\theta(Az_i^{\pm 1};p)}
\prod_{m=2}^{2n+3}\frac{\theta(At_m^{-1};p)}{\theta(t_1t_m;p)}-1
=\frac{t_1\theta(t_1A;p)}{\prod_{m=2}^{2n+3}\theta(t_1t_m;p)}
\sum_{i=1}^n\frac{1}{z_i\theta(z_i^2;p)}
}&&
\\ &&
\times
\prod_{j=1,\neq i}^n\frac{\theta(t_1z_j^{\pm 1};p)}{\theta(z_iz_j^{\pm 1};p)}
\left(z_i^{2n+2}\frac{\prod_{m=1}^{2n+3}\theta(t_mz_i^{-1};p)}
{\theta(Az_i^{-1};p)} - \frac{\prod_{m=1}^{2n+3}\theta(t_mz_i;p)}
{\theta(Az_i;p)}\right).
\label{eqn-exp-C}\end{eqnarray}
Both sides of this equality are invariant under the transformation
$z_1\to pz_1$ and have equal sets of poles (singularities at the points
$z_1=z_j,z_j^{-1}, j=2,\ldots, n,$ and $z_1=\pm p^{k/2},\ k\in\Z$
on the right-hand side are cancelled) with their residues.
Therefore the functions on both sides of the equality
(\ref{eqn-exp-C}) differ only by an additive constant,
independent on $z_1$. This constant equals to zero which
follows from a trivial check of equality (\ref{eqn-exp-C}) at $z_1=t_1$.

Integrating (\ref{eqn-C}) over the variables $z\in C^n$, we obtain
\begin{equation}
I(qt_1,t_2,\ldots,t_{2n+3})-I(t)=\sum_{i=1}^n
\left(\int_{C^{i-1}\times(q^{-1}C)\times C^{n-i}}-\int_{C^n}\right)g_i(z,t)\frac{dz}{z},
\label{int-eqn-C}\end{equation}
where $I(t)=\int_{C^n}\rho(z,t;C_n)dz/z$ and
$q^{-1}C$ denotes the contour $C$ dilated with respect to the zero point.

Poles of the function (\ref{g-C}) in variable $z_i$ converge to zero
along the point $z_i=t_mq^ap^b,$  $A^{-1}q^{a}p^{b+1}$ and diverge to infinity
at $z_i=t_m^{-1}q^{-1-a}p^{-b},Aq^{-a}p^{-b-1}$, where $m=1,\ldots,2n+3$,
$a,b=0,1,\ldots$. For $|t_m|<1$ and $|p|<|A|$ the region
$1\leq |z_i|\leq |q|^{-1}$ does not contain the poles, so that we
can set $C=\T$, deform back $q^{-1}\mathbb{T}$ to $\mathbb{T}$
in  (\ref{int-eqn-C}) and obtain the equality $I(qt_1,t_2,\ldots,t_{2n+3})=I(t)$.

Repeating almost literally the procedure of analytical continuation
used in the $n=1$ case, we find that $I$ is a constant, which
depends only on $p$ and $q$. Its value is found by considering
the limits $t_jt_{j+n}\to 1,\, j=1,\ldots,n,$
analogous to the $n=1$ case, which yields the right-hand side of
(\ref{ell-int-C}).
\end{proof}

Formula \eqref{C-typeI} was suggested and partially justified
in \cite{die-spi:selberg}, and it was proved completely by different
methods in \cite{rai:trans,spi:short,rs:det}. In a  special
limit $p\to 0$, it is reduced to one of the integration formulae of
Gustafson \cite{gus:some1}.

The root system $C_n$ consists of the set of vectors from $\R^n$ of the form
$X({C_n})=\{\pm 2e_i, \pm e_i \pm e_j,\, i< j\}\big|_{i,j=1,\ldots,n}$,
where $e_i$ is an orthonormal basis of $\R^n$. Denoting
$z_i=\exp(e_i)$, we see that the denominator of integral's kernel
\eqref{C-typeI} contains a product over roots of $C_n$ of the form
\begin{equation*}
\prod_{\alpha\in X({C_n})}\Gamma(\e^{\alpha};p,q)
=\prod_{i=1}^n\Gamma(z_i^{\pm 2};p,q)\prod_{1\leq i<j\leq n}
\Gamma(z_i^{\pm1}z_j^{\pm1};p,q).
\end{equation*}
The root system $BC_n$ contains additionally the vectors
$\{\pm e_1,\ldots,\pm e_n\}$, but the general rules of the appearance
of these vectors  in integrals' kernels are not established yet.
The Weyl group of these systems $S_{n}\times \Z_2^n$ is
a symmetry of the integral kernel.

The  $C_n$-elliptic beta integral of type II is built with the
help of formula \eqref{C-typeI} by a purely algebraic means
\cite{die-spi:selberg}.

\begin{theorem}
Let complex parameters $t, t_m (m=1,\ldots , 6), p$ and $q$ satisfy
conditions $|p|, |q|,$ $|t|,$ $|t_m| <1,$ and
$t^{2n-2}\prod_{m=1}^6t_m=pq$. Then
\begin{eqnarray}\nonumber
\kappa_n\int_{\T^n} \prod_{1\leq j<k\leq n}
\frac{\eg(tz_j^{\pm 1} z_k^{\pm 1};p,q)}{\eg(z_j^{\pm 1} z_k^{\pm 1};p,q)}
\prod_{j=1}^n\frac{\prod_{m=1}^6\eg(t_mz_j^{\pm 1};p,q)}{\eg(z_j^{\pm2};p,q)}
\frac{dz_1}{z_1}\cdots\frac{dz_n}{z_n}
\\
= \prod_{j=1}^n\left(\frac{\eg(t^j;p,q)}{\eg(t;p,q)}
\prod_{1\leq m<s\leq 6}\eg(t^{j-1}t_mt_s;p,q )\right).
\label{SintB}\end{eqnarray}
\end{theorem}
\begin{proof}
We denote the integral on the left-hand side of \eqref{SintB} as
$I_n(t,t_1,\ldots,t_5)$ and consider the $(2n-1)$-tuple integral
\begin{eqnarray}
&& \makebox[-1em]{}\kappa_{n}\kappa_{n-1}\int_{\T^{2n-1}}
\prod_{1\leq j<k\leq n}\frac{1}{\eg(z_j^{\pm 1} z_k^{\pm 1};p,q)}
 \prod_{j=1}^n\frac{\prod_{r=0}^5\eg(t_rz_j^{\pm 1};p,q)}
{\eg(z_j^{\pm2};p,q)}
\nonumber \\ &&
\times
\prod_{\stackrel{1\leq j\leq n}{1\leq k\leq n-1}}
\eg(t^{1/2}z_j^{\pm 1} w_k^{\pm 1};p,q)
\prod_{1\leq j<k\leq n-1}\frac{1}{\eg(w_j^{\pm 1} w_k^{\pm 1};p,q)}\nonumber \\
&&
\times \prod_{j=1}^{n-1}
\frac{\eg(w_j^{\pm 1} t^{n-3/2}\prod_{s=1}^5t_s;p,q)}
{\eg(w_j^{\pm2},w_j^{\pm 1} t^{2n-3/2}\prod_{s=1}^5t_s;p,q)}
\frac{dw_1}{w_1}\cdots\frac{dw_{n-1}}{w_{n-1}}
\frac{dz_1}{z_1}\cdots\frac{dz_n}{z_n},
\label{compint}
\end{eqnarray}
with $p, q, t$ and $t_r$, $r=0,\ldots,5,$ lying inside the unit circle
such that $t^{n-1}\prod_{r=0}^5t_r=pq$. Integration over the variables
$w_j$ with the help of formula \eqref{C-typeI} brings the expression
\eqref{compint} to the form
$\eg^n(t;p,q)\eg^{-1}(t^n;p,q)$ $I_n(t,t_1,\ldots,t_5)$
(where it is assumed that $t_6=pqt^{2-2n}/\prod_{j=1}^5t_j$).
Because the integrand is bounded on the contour of integration,
we can change the order of integrations. Then the integration over
the $z_k$-variables with the help of formula \eqref{C-typeI}
converts expression \eqref{compint} to
$$
\eg^{n-1}(t;p,q) \prod_{0\leq r< s\leq 5} \eg(t_rt_s;p,q)\,
I_{n-1}(t,t^{1/2}t_1,\ldots,t^{1/2}t_5).
$$
As a result, we obtain a recurrence relation connecting integrals of
different dimension $n$:
$$
I_n(t,t_1,\ldots,t_5)= \frac{\eg(t^n;p,q)}{\eg(t;p,q)}\makebox[-0.5em]{}
\prod_{0\leq r<s\leq 5}\makebox[-0.5em]{}\eg(t_rt_s;p,q)\;
I_{n-1}(t,t^{1/2}t_1,\ldots,t^{1/2}t_5).
$$
Using known initial condition at $n=1$ \eqref{ell-int},
 we find \eqref{SintB} by recursion.
\end{proof}

If one expresses $t_6$ via other parameters and removes the multipliers
 $pq$ in the arguments of the elliptic gamma functions with the help
of the reflection formula, then it is easy to pass to the limit $p\to0$
for fixed parameters. This leads to a multiple $q$-beta integral
of Gustafson \cite{gus:some2}. A number of other limiting transitions
in parameters leads to the Selberg integral -- a fundamentally important
integral with a large number of applications in mathematical physics
\cite{for-war}:
\begin{eqnarray}
\lefteqn{\int_0^1\cdots\int_0^1
\prod_{1\leq j\leq n} x_j^{\alpha-1} (1-x_j)^{\beta-1}
\prod_{1\leq j<k\leq n}|x_j-x_k|^{2\gamma}\; dx_1\cdots dx_n } && \nonumber \\
&& =\prod_{1\leq j\leq n}
\frac{\Gamma(\alpha+(j-1)\gamma) \Gamma(\beta+(j-1)\gamma) \Gamma (1+j\gamma)}
     {\Gamma (\alpha+\beta+(n+j-2)\gamma) \Gamma (1+\gamma)}, \label{Sint}
\end{eqnarray}
where $\text{Re}(\alpha ), \text{Re} (\beta) >0$ and
$\text{Re}(\gamma ) > - \min (1/n,
\text{Re} (\alpha )/(n-1) ,\text{Re} (\beta)/(n-1) )$.

Therefore formula \eqref{SintB} represents an elliptic
analogue of the Selberg integral. It can also be interpreted
as an elliptic generalization of the
Macdonald--Morris constant term identities for the $BC_n$-root system \cite{mac}.
The given proof is taken from the work \cite{die-spi:selberg}.
It models the proof of the Selberg integral suggested by Anderson
\cite{and:short} and represents a generalization of the method used
by Gustafson in \cite{gus:some2} for proving the corresponding $q$-beta integral.

\subsection{Integrals for the root system $A_n$}\

Different elliptic beta integrals on the root system $A_n$
have been proposed in papers \cite{spi:theta2,spi-war:inversions}.
By analogy with the $C_n$-cases integration formulae depending
on $2n+3$ parameters will be considered as the type I integrals.
Exactly computable integrals with a smaller number of parameters,
which are derived with the help of type I integrals,
will be classified as the type II integrals.
Let us list these formulae omitting their derivations.

\begin{theorem}
Let $|p|, |q|<1$ and $2n+4$ parameters $t_m, s_m,\, m=1,\ldots, n+2,$
satisfy the constraints $|t_m|, |s_m|<1$ and $ST=pq$, where
$S=\prod_{m=1}^{n+2}s_m$ and $T=\prod_{m=1}^{n+2}t_m$. Then
\begin{eqnarray}\nonumber
&& \kappa_{n}^A\int_{{\mathbb T}^n}
\prod_{1\le j<k\le n+1}\frac{1}{\eg(z_iz_j^{-1},z_i^{-1}z_j;p,q)}
\,\prod_{j=1}^{n+1}\prod_{m=1}^{n+2}\eg(s_mz_j,t_mz_j^{-1};p,q)
\,\frac{dz}{z}
\\ && \makebox[4em]{}
=\prod_{m=1}^{n+2} \eg(Ss_m^{-1},Tt_m^{-1};p,q)
\prod_{k,m=1}^{n+2} \eg(s_kt_m;p,q),
\label{AI}\end{eqnarray}
where $z_1z_2\cdots z_{n+1}=1$ and
$$
\kappa_{n}^A=\frac{(p;p)_\infty^n(q;q)_\infty^n}{(2\pi i)^n(n+1)!}.
$$
\end{theorem}

In this type I integral we have a split of  $2n+4$ parameters
(homogeneous in the $C_n$-case) with the fixed product into two
groups with $n+2$ elements. This integration formula was proposed
and partially justified in \cite{spi:theta2}, various
different complete proofs are given in \cite{rai:trans,spi:short}.
In the simplest possible $p\to 0$ limit there appears one
of the Gustafson integrals \cite{gus:some1}.

The root system $A_n$ consists of the vectors
$X(A_n)=\{e_i-e_j,\, i\neq j\}\big|_{i,j=1,\ldots,n+1}$,
where $e_i$ is an orthonormal basis of  $\R^{n+1}$.
These vectors lie in the hyperplane orthogonal to the vector
$E=\sum_{i=1}^{n+1}e_i$. Setting $z_i=\exp(e_i-E/(n+1))$,
we obtain $z_1\cdots z_{n+1}=1$. The denominator of the kernel
of integral \eqref{AI} contains a product of the form
\begin{equation*}
\prod_{\alpha\in X(A_n)}\Gamma(\e^{\alpha};p,q)
=\prod_{1\leq i<j\leq n+1}\Gamma(z_i/z_j,z_j/z_i;p,q).
\end{equation*}
The full integral kernel is invariant with respect to the
$A_n$-Weyl group $S_{n+1}$.

There exists an additional independent $A_n$-integral of type I
\cite{spi-war:inversions}.

\begin{theorem}
Let $|p|, |q|<1$ and $2n+3$ parameters $t_k,\, k=1,2,\ldots ,n,$
and $s_m,\, m=1,2,\ldots,n+3,$ satisfy the constraints
$|t_k|<1,\, |s_m|<1$ and $|pq|<|St_k|$, where $S=\prod_{m=1}^{n+3}s_m$.
Then the following explicit integration formula is true:
\begin{eqnarray} \nonumber
&&  \kappa_{n}^A\int_{{\mathbb T}^n}
\prod_{1\le i<j\le n+1}\frac{\eg(Sz_i^{-1}z_j^{-1};p,q)}
{\eg(z_iz_j^{-1},z_i^{-1}z_j;p,q)}
\\ \nonumber
&& \makebox[2em]{}  \times
\prod_{j=1}^{n+1}
\frac{\prod_{k=1}^{n}\eg(t_kz_j;p,q)\prod_{m=1}^{n+3}\eg(s_mz_j^{-1};p,q)}
{\prod_{k=1}^{n}\eg(St_kz_j^{-1};p,q)} \frac{dz}{z}
\\ && \makebox[4em]{}
= \prod_{k=1}^{n} \prod_{m=1}^{n+3}
\frac{\eg(t_ks_m;p,q)}{\eg(St_ks_m^{-1};p,q)}
\prod_{1\le l<m\le n+3} \eg(Ss_l^{-1}s_m^{-1};p,q),
\lab{AI'}\end{eqnarray}
where $z_1\cdots z_{n+1}=1$.
\end{theorem}

Here we have a split of $2n+3$ independent variables into two
homogeneous groups with $n$ and $n+3$ elements.
In the limit  $p\to0$ there appears a $q$-beta integral which was
not considered in the literature until the derivation of the
elliptic case, as well as its plain hypergeometric degeneration
appearing in the $q\to 1$ limit.

There are several $A_n$-elliptic beta integrals of type II \cite{spi:theta2}.
We denote
\begin{eqnarray} \nonumber
&& I^{II}({\bf t},{\bf s};p,q;A_n)=\kappa_{n}^A
\int_{\T^n}\prod_{1\leq i<j\leq n+1}
\frac{\eg(tz_iz_j;p,q)}{\eg(z_iz_j^{-1},z_i^{-1}z_j;p,q) }
 \\ && \makebox[4em]{} \times
 \prod_{j=1}^{n+1} \prod_{k=1}^{n+1}\eg(t_kz_j^{-1};p,q)
\prod_{i=1}^4\eg(s_iz_j;p,q)\frac{dz_1}{z_1}\cdots\frac{dz_n}{z_n},
\lab{AIIb}\end{eqnarray}
where $|t_k|,\, |s_i|<1,\, k=1,\ldots,n+1,\, i=1,\ldots,4,$ and
$$
t^{n-1}\prod_{k=1}^{n+1}t_k\prod_{i=1}^4s_i=pq, \qquad
\prod_{j=1}^{n+1}z_j=1.
$$

\begin{theorem}
Under the described restrictions on the parameters
the following $A_n$-integration formulae are true. For odd $n$, we have
\begin{eqnarray} \nonumber
&& I^{II}({\bf t},{\bf s};p,q;A_n)
=\prod_{1\leq j<k\leq n+1} \eg(tt_jt_k;p,q)
\prod_{k=1}^{n+1}\prod_{i=1}^4\eg(t_ks_i;p,q)
 \\ && \makebox[4em]{} \times
\frac{\eg(t^{\frac{n+1}{2}},A;p,q)}{\eg(t^{\frac{n+1}{2}}A;p,q)}
\prod_{1\leq i<m\leq 4}\eg(t^{\frac{n-1}{2}}s_is_m;p,q),
\lab{AIIb-even}\end{eqnarray}
where $A=\prod_{k=1}^{n+1}t_k$; for even $n$
\begin{eqnarray} \nonumber
&& I^{II}({\bf t},{\bf s};p,q;A_n)
= \prod_{1\leq j<k\leq n+1} \eg(tt_jt_k;p,q)
\prod_{k=1}^{n+1}\prod_{i=1}^4\eg(t_ks_i;p,q)
\\ && \makebox[6em]{} \times
\eg(A ;p,q) \prod_{i=1}^4\frac{\eg(t^{\frac{n}{2}}s_i;p,q)}
{\eg(t^{\frac{n}{2}}As_i;p,q) }.
\lab{AIIb-odd}\end{eqnarray}
\end{theorem}

These formulae contain only $n+5$ free parameters.
They can be derived as direct consequences of the elliptic beta
integrals of type I \eqref{C-typeI} and \eqref{AI}  \cite{spi:theta2}.
In the simplest limit $p\to 0$ they are reduced to the
Gustafson-Rakha $q$-beta integrals \cite{gus-rak:beta}.

In order to describe a different $A_n$-integral of type II,
we need ten complex parameters $p,q,$ $t, s,$ $t_1,t_2,t_3,$ $s_1,s_2, s_3$
with the constraints
$(ts)^{n-1}\prod_{k=1}^3t_ks_k=pq$ and $|p|, |q|, |t|, |s|, |t_k|, |s_k|<1$.
We define now the integral
\begin{eqnarray} \lab{AIIa}
&& I^{II}({\bf t},{\bf s};p,q;A_n)=
\\ && \makebox[1em]{}
\kappa_{n}^A\int_{\T^n} \prod_{1\leq i<j\leq n+1}
\frac{\eg(tz_iz_j,sz_i^{-1}z_j^{-1};p,q)}
{\eg(z_iz_j^{-1},z_i^{-1}z_j;p,q)}
\prod_{j=1}^{n+1}\prod_{k=1}^3\eg(t_kz_j,s_kz_j^{-1};p,q)
\frac{dz}{z}.
\nonumber\end{eqnarray}

\begin{theorem}
Under the indicated restrictions on the parameters
the following integration formulae for the root system $A_n$ are true.
For odd $n$, we have
\begin{eqnarray} \nonumber
&&
I^{II}({\bf t},{\bf s};p,q;A_n)=
\eg(t^{\frac{n+1}{2}},s^{\frac{n+1}{2}};p,q)
\prod_{1\leq i<k\leq 3} \eg(t^{\frac{n-1}{2}}t_it_k,s^{\frac{n-1}{2}}s_is_k;p,q)
\\  && \makebox[4em]{} \times
\prod_{j=1}^{(n+1)/2}\prod_{i,k=1}^3\eg((ts)^{j-1}t_is_k;p,q)
\lab{AIIa-odd} \\ && \makebox[4em]{}  \times
\prod_{j=1}^{(n-1)/2}\left(\eg((ts)^j;p,q)\prod_{1\leq i<k\leq 3}
\eg(t^{j-1}s^jt_it_k,t^js^{j-1}s_is_k;p,q)\right);
\nonumber\end{eqnarray}
for even $n$
\begin{eqnarray} \nonumber
&& I^{II}({\bf t},{\bf s};p,q;A_n)
= \prod_{i=1}^3\eg(t^{\frac{n}{2}}t_i,s^{\frac{n}{2}}s_i;p,q) \qquad\qquad\qquad
\\  &&\makebox[4em]{} \times
\eg(t^{\frac{n}{2}-1}t_1t_2t_3,s^{\frac{n}{2}-1}s_1s_2s_3;p,q)
\prod_{j=1}^{n/2}\Biggl(\eg((ts)^j;p,q)
\lab{AIIa-even}\\ && \makebox[4em]{}
\times \prod_{i,k=1}^3\eg((ts)^{j-1}t_is_k;p,q)
\prod_{1\leq i<k\leq 3}\eg(t^{j-1}s^jt_it_k,t^js^{j-1}s_is_k;p,q)\Biggr).
\nonumber \end{eqnarray}
\end{theorem}

These formulae can be derived as direct consequences of the $C_n$-elliptic beta
integrals of type I and II and the $A_n$-integral of type I \eqref{AI}
\cite{spi:theta2}.
In the simplest possible limit $p\to 0$ these integrals are reduced
to the Gustafson $q$-hypergeometric integrals \cite{gus:some2}.

\section{Some multiple series summation formulae}

There are several known summation formulae for multiple elliptic hypergeometric
series representing multivariate extensions of the Frenkel--Turaev sum
\cite{spi:theta2,rai:trans,war:summation,ros:elliptic,ros-sch1,
die-spi:modular,cos-gus}.
A direct connection with the sums of residues of particular
pole sequences of the kernels of elliptic beta integrals on root systems
was established for some of them.
Let us describe some of these formulae without proofs.

\begin{theorem}\label{mgsum:thm}
Let $|p|<1$ and the parameter $q\in\C$ is not equal to an integer power
of $p$. Then the following summation formula is true
\begin{eqnarray}\label{mg:sum}
&& \sum_{\stackrel{0\leq \lambda_j\leq N_j}{j=1,\ldots ,n}}
q^{\sum_{j=1}^n j\lambda_j}
\prod_{1\leq j<k\leq n}
\frac{\theta (t_jt_kq^{\lambda_j+\lambda_k},
              t_jt_k^{-1}q^{\lambda_j-\lambda_k} ;p)}
     {\theta (t_jt_k,t_jt_k^{-1} ;p)}  \\
&& \quad\qquad\qquad\qquad\times \prod_{1\leq j\leq n} \Biggl(
\frac{\theta (t_j^2q^{2\lambda_j};p)}{\theta (t_j^2;p)}
\prod_{0\leq r\leq 2n+3} \frac{\theta (t_jt_r)_{\lambda_j}}
     {\theta (qt_jt_r^{-1})_{\lambda_j}}\Biggr)   \nonumber\\
&&= \theta (qa^{-1}b^{-1},qa^{-1}c^{-1},qb^{-1}c^{-1})_{N_1+\cdots +N_n}
\nonumber\\
&&\qquad\times\prod_{1\leq j<k\leq n}
\frac{\theta (qt_jt_k)_{N_j}\theta (qt_jt_k )_{N_k}}
     {\theta (qt_jt_k)_{N_j+N_k}} \nonumber\\
&&\qquad \times \prod_{1\leq j\leq n}
\frac{\theta (qt_j^2)_{N_j}}
     {\theta (qt_ja^{-1},qt_jb^{-1},
    qt_jc^{-1},q^{1+N_1+\cdots+N_n-N_j}t_j^{-1}a^{-1}b^{-1}c^{-1})_{N_j}},
    \nonumber
\end{eqnarray}
where $\theta(a)_\lambda=\Gamma(aq^\lambda;p,q)/\Gamma(a;p,q)$ and
\begin{equation*}
\begin{array}{lll}
q^{-1}\prod_{r=0}^{2n+3}t_r=1 & & (\text{the\ balancing\ condition}) ,\\
q^{N_j}t_jt_{n+j}=1, & j=1,\ldots ,n & (\text{the\ termination\ condition}),
\end{array}
\end{equation*}
with $N_j=0,1,\ldots$, $j=1,\ldots ,n$, and $a=t_{2n+1}, b=t_{2n+2}, c=t_{2n+3}$.
\end{theorem}

This formula was deduced in \cite{die-spi:modular} from the $C_n$-elliptic
beta integral of type I (which was not proven at that moment yet).
Its first recursive proof was obtained in \cite{ros:elliptic}.
In the limit $p\to 0$ it degenerates to a multivariate
$_8\varphi_7$-sum, which was found in \cite{den-gus:q-beta}.

\begin{theorem}
Let $N=0,1,\ldots$ and the parameters $p,q,t, t_0,\ldots,t_5\in\C$
satisfy the restrictions $|p|<1$ and
\begin{equation}\label{pc}
\begin{array}{rl}
q^{-1}t^{2n-2}\prod_{r=0}^5t_r=1 & \text{(the balancing condition)}, \\
q^Nt^{n-1}t_0t_4=1 & \text{(the termination condition)}.
\end{array}
\end{equation}
Then the following summation formula for a multiple elliptic hypergeometric
series is true
\begin{eqnarray}\nonumber
&& \makebox[-3em]{}
\sum_{0\leq \lambda_1\leq \lambda_2\leq \cdots \leq \lambda_n\leq N}
q^{\sum_{j=1}^m\lambda_j} t^{2 \sum_{j=1}^m (n-j)\lambda_j}
\prod_{1\leq j<k\leq m}
\Biggl(
\frac{\theta(\tau_k\tau_jq^{\lambda_k+\lambda_j},
\tau_k\tau_j^{-1}q^{\lambda_k-\lambda_j};p)}
{\theta(\tau_k\tau_j,\tau_k\tau_j^{-1};p)} \\ \nonumber
&& \makebox[8em]{}\times
\frac{\theta(t\tau_k\tau_j)_{\lambda_k+\lambda_j}}
     {\theta(qt^{-1}\tau_k\tau_j)_{\lambda_k+\lambda_j}}
\frac{\theta(t\tau_k\tau_j^{-1})_{\lambda_k-\lambda_j}}
     {\theta(qt^{-1}\tau_k\tau_j^{-1})_{\lambda_k-\lambda_j}}
\Biggr) \\ \nonumber
&&\makebox[3em]{} \times \prod_{j=1}^m
\Biggl( \frac{\theta(\tau_j^2q^{2\lambda_j};p)}{\theta(\tau_j^2;p)}
\prod_{r=0}^5 \frac{\theta(t_r\tau_j)_{\lambda_j}}
     {\theta(qt_r^{-1}\tau_j)_{\lambda_j}}\Biggr)
\\ \makebox[2em]{}
&= &\prod_{j=1}^n
\frac{\theta(qt^{n+j-2}t_0^2)_N
      \prod_{1\leq r <s \leq 3} \theta(qt^{1-j} t_r^{-1}t_s^{-1})_N}
      {\theta(qt^{2-n-j}\prod_{r=0}^3t_r^{-1})_N
      \prod_{r=1}^3 \theta(qt^{j-1}t_0t_r^{-1})_N}.
\label{bcn-sum}\end{eqnarray}
Here we use the notation $\tau_j=t_0t^{j-1}$, $j=1,\ldots ,n$.
\end{theorem}

This summation formula was suggested by Warnaar \cite{war:summation}.
It follows from the $C_n$-elliptic beta integral of the type II
\cite{die-spi:elliptic} and was proven for the first time recursively
in the paper \cite{ros:proof}.

For the root system $A_n$ the type I elliptic hypergeometric
sum has the following form:
\begin{eqnarray} \nonumber
\lefteqn{
\sum_{0\leq \lambda_j \leq N_j \atop j=1,\ldots, n}
q^{\sum_{j=1}^nj\lambda_j}
\prod_{j=1}^n\frac{\theta(t_jq^{\lambda_j+|\lambda|};p)}{\theta(t_j;p)}
\prod_{1\leq i<j \leq n} \frac{\theta(t_it_j^{-1}q^{\lambda_i-\lambda_j};p)}
{\theta(t_it_j^{-1};p)}  } &&
\\ \nonumber
&& \times \prod_{i,j=1}^n\frac{\theta(t_it_j^{-1}q^{-N_j})_{\lambda_i}}
{\theta(qt_it_j^{-1})_{\lambda_i}}
\prod_{j=1}^n\frac{\theta(t_j)_{|\lambda|}}
{\theta(t_jq^{1+N_j})_{|\lambda|}}
\\ \nonumber
&& \times \frac{\theta(b,c)_{|\lambda|}}
{\theta(q/d, q/e)_{|\lambda|}}
\prod_{j=1}^n \frac{\theta(dt_j, et_j)_{\lambda_j}}
{\theta(t_jq/b, t_jq/c)_{\lambda_j}}
\\
&&
=\frac{\theta(q/bd,q/cd)_{|N|}}{\theta(q/d,q/bcd)_{|N|}}
\prod_{j=1}^n\frac{\theta(t_jq,t_jq/bc)_{N_j}}
{\theta(t_jq/b,t_jq/c)_{N_j}},
\label{A_n}\end{eqnarray}
where $|\lambda|=\lambda_1+\cdots+\lambda_n$,
$|N|=N_1+\cdots+N_n$ and $bcde=q^{1+|N|}$.
This formula was proven recursively in \cite{ros:elliptic}. It follows also
from the analysis of residues for the type I integral
\eqref{AI} \cite{spi:theta2} and for $p\to0$ reduces to
the multiple $_8\varphi_7$-sum of Milne \cite{mil:multiple}.

Omitting a number of other established summation formulae for multiple
elliptic hypergeometric series, we describe a hypothesis from
 \cite{spi:theta2}.

\begin{conjecture}
Let $|p|<1$, $N=0,1,\ldots,$ and $n$ parameters $t_k\in\C$, $k=1,\ldots,n,$
satisfy the balancing condition $\prod_{k=1}^n t_k=q^{-N}$.
Then the following summation formula is true:
\begin{eqnarray}\nonumber
&&
\sum_{\stackrel{\lambda_k=0,\ldots,N}{\lambda_1+\ldots+\lambda_n=N}}
\frac{
\prod_{1\leq i<j\leq n}\theta(tt_it_j)_{\lambda_i+\lambda_j}
\prod_{i=1}^n\prod_{j=n+1}^{n+3}\theta(tt_it_j)_{\lambda_i}
\prod_{i,j=1}^n\theta(t_it_j^{-1})_{-\lambda_j}
}{
\prod_{i,j=1;i\neq j}^n\theta(t_it_j^{-1})_{\lambda_i-\lambda_j}
\prod_{j=1}^n\theta(t^{n+1}t_j^{-1}\prod_{k=1}^{n+3}t_k)_{-\lambda_j}
} \\ && \makebox[2em]{}
= \left\{ \begin{aligned} 
\frac{\theta(1)_{-N}}{\theta(t^{n/2})_{-N}
\prod_{n+1\leq i<j\leq n+3}\theta(t^{(n+2)/2}t_it_j)_{-N}},
& \quad n \quad \text{even,}  \\
\frac{\theta(1)_{-N}}{\prod_{i=n+1}^{n+3}\theta(t^{(n+1)/2}t_i)_{-N}
\theta(t^{(n+3)/2}\prod_{i=n+1}^{n+3}t_i)_{-N} },
& \quad n \quad \text{odd.}  \\
\end{aligned} \right.
\label{new-sums} \end{eqnarray}
\end{conjecture}

It is conjectured also that this formula is related to sums of residues
for the $A_n$-integral \eqref{AIIb}. For $p\to0$ there appears a
Gustafson--Rakha summation formula for a multiple $_8\varphi_7$-series
\cite{gus-rak:beta}.

\section{Symmetry transformations for multiple integrals}

General elliptic hypergeometric integrals of type I for the root system
 $C_n$ have the form
\begin{eqnarray*} &&
I_n^{(m)}(t_1,\ldots,t_{2n+2m+4})
\\ && \makebox[1em]{}
=\kappa_n\int_{\T^n}
\prod_{1\leq i<j\leq n} \frac{1}{\eg(z_i^{\pm 1}z_j^{\pm 1};p,q)}
\prod_{j=1}^n\frac{\prod_{i=1}^{2n+2m+4}\eg(t_iz_j^{\pm 1};p,q)}
{\eg(z_j^{\pm 2};p,q)}\frac{dz_j}{z_j},
\end{eqnarray*}
where $|t_j|<1$,
$$
\prod_{j=1}^{2n+2m+4}t_j=(pq)^{m+1},\qquad
\kappa_n=\frac{(p;p)_\infty^n(q;q)_\infty^n}{(4\pi i)^n n!}.
$$
For $n=0$ we set $I^{(m)}_0=1$. The integral $I_1^{(1)}$
coincides with the $V$-function, and $I_n^{(0)}$ -- with the elliptic
beta integral \eqref{C-typeI}. After degeneration of $I^{(m)}_n$-functions
to the level of ordinary beta integral, they reduce to the
Dixon integrals \cite{rai:limits,Dixon}.
In \cite{rai:trans}, Rains proved the following transformation formula:
\begin{equation}
I_n^{(m)}(t_1,\ldots,t_{2n+2m+4})=\prod_{1\leq r<s\leq 2n+2m+4}\eg(t_rt_s;p,q)\;
I_m^{(n)}\left(\frac{\sqrt{pq}}{t_1},\ldots,\frac{\sqrt{pq}}{t_{2n+2m+4}}\right).
\label{trafo}\end{equation}
It represents a direct generalization of the third symmetry transformation for
the $V$-function described in \eqref{E7-3}.

In \cite{rs:det}, the following determinant representation
has been found:
\begin{eqnarray}\label{detrep} &&\makebox[-1em]{}
I_n^{(m)}(t_1,\ldots,t_{2n+2m+4})= \prod_{1\leq i<j\leq n}
\frac{1}{a_j\theta(a_ia_j^{\pm 1};p)b_j\theta(b_ib_j^{\pm 1};q)}
\\ && \makebox[-1em]{} \times
\det_{1\le i,j\le n}\left(\kappa\int_\T \frac{\prod_{r=1}^{2n+2m+4}
\eg(t_r z^{\pm 1};p,q)}{\eg(z^{\pm2};p,q)}
\prod_{k\ne i} \theta(a_kz^{\pm 1};p)
\prod_{k\ne j} \theta(b_kz^{\pm 1};q)
\; \frac{dz}{z} \right),
\nonumber \end{eqnarray}
where $a_j, b_j$ are some arbitrary parameters.
For the choice $a_i=t_j$, $b_j=t_{n+j}$, $j=1,\ldots, n$,
there appears the determinant of the matrix
\[
T_q(t_i)^{-1} T_p(t_{n+j})^{-1}
I^{(m+n-1)}_1(qt_1,\dots,qt_n,pt_{n+1},\dots,pt_{2n},t_{2n+1},\dots,t_{2n+2m+4}),
\]
where $T_q(t_k)$ is the $q$-shift operator, $T_q(t_k)f(t_k)=f(qt_k)$.
For $m=0$, we obtain thus an exactly computable determinant
of univariate elliptic hypergeometric integrals.

In order to prove \eqref{detrep}, it is necessary to write
$$
\prod_{1\leq i<j\leq n} \frac{1}{\eg(z_i^{\pm 1}z_j^{\pm 1};p,q)}
=\prod_{1\le i<j\le n} z_i^{-1}\theta(z_i z_j^{\pm 1};p)
\prod_{1\le i<j\le n}z_i^{-1}\theta(z_i z_j^{\pm 1};q)
$$
and connect both multipliers on the right-hand side
with the elliptic analogue of the Cauchy determinant \eqref{d-cau}.
Using the Heine formula
$$
\frac{1}{n!}\int\det_{1\le i,j\le n} \phi_i(z_j)
\det_{1\le i,j\le n} \psi_i(z_j)
\prod_{1\le i\le n} d\mu(z_i)
= \det_{1\le i,j\le n}\int \phi_i(z)\psi_j(z) d\mu(z)
$$
with $\phi_i(z_j)=a_i/\theta(a_iz_j^{\pm 1};p)$ and
$\psi_i(z_j)=b_i/\theta(b_iz_j^{\pm 1};q)$, one obtains the necessary result.

The following $(n+2)$-term recurrence relation takes place
\cite{die-spi:selberg,rs:det}
\begin{equation}
\sum_{i=1}^{n+2}
\frac{t_i}{\prod_{j=1,\, \ne i}^{n+2} \theta(t_i t_j^{\pm 1};p)}
T_q(t_i) I^{(m)}_n(t_1,\dots,t_{2n+2m+4})
=0,
\label{rec-rel1}\end{equation}
where $\prod_{j=1}^{2n+2m+4}t_j=(pq)^mp$.
Indeed, if we write the same recurrence relation for the $I^{(m)}_n$-integral kernel,
then it reduces to identity \eqref{pf2} after taking away a common
multiplier. Integrating the latter relation over appropriate contours,
we obtain \eqref{rec-rel1}.

The determinant representation permits one to derive another
$(m+2)$-term recurrence relation
\begin{equation}
\sum_{k=1}^{m+2}
\frac{\prod_{m+3\le i\le 2n+2m+4} \theta(t_it_k/q;p)}
     {t_k\prod_{1\le i\le m+2;i\ne k} \theta(t_i/t_k;p)}
T_q(t_k)^{-1}
I^{(m)}_n(t_1,\dots,t_{2n+2m+4}) =0,
\label{rec-rel2}\end{equation}
where $t_1\cdots t_{2m+2n+4} = (pq)^{m+1}q$.
Note that transformation \eqref{trafo} maps equation
\eqref{rec-rel1} to \eqref{rec-rel2}, and vice versa.
Therefore an analysis of common solutions of these equations
provides an alternative way of proving transformation \eqref{trafo} \cite{rs:det}.

Let us consider symmetry transformations for other multidimensional
integrals. For the parameters $t,t_1,\ldots,t_8\in\C$ satisfying
the restrictions $|t|, |t_j|<1$ and $ t^{2n-2}\prod_{j=1}^8t_j=p^2q^2$,
we define the integral
\begin{eqnarray}\label{type2} &&
I(t_1,\ldots,t_8;t;p,q)=\prod_{1\leq j<k\leq 8}\Gamma(t_jt_k;p,q,t)
\\ && \makebox[1em]{} \times
\kappa_n \int_{\T^n}\prod_{1\le j<k\le n}\!
\frac{\eg(tz_j^{\pm 1}z_k^{\pm1};p,q)}
{\eg(z_j^{\pm1}z_k^{\pm1};p,q)}
\prod_{j=1}^n \frac{\prod_{k=1}^8\eg(t_kz_j^{\pm1};p,q)}
{\eg(z_j^{\pm 2};p,q)}\frac{dz_j}{z_j}.
\nonumber\end{eqnarray}
Here
$$
\Gamma(z;p,q,t)=\prod_{j,k,l=0}^\infty(1-zt^jp^kq^l)(1-z^{-1}t^{j+1}p^{k+1}
q^{l+1})
$$
is the elliptic gamma function of a higher order, which is connected to
the Barnes multiple gamma function $\Gamma_4(u;\mathbf{\omega})$
and satisfies the equation
$$
\Gamma(tz;p,q,t)=\Gamma(z;p,q)\Gamma(z;p,q,t).
$$
Function \eqref{type2} generalizes the $V$-function to the type II
integrals for the root system $BC_n$. In \cite{rai:trans}, Rains
has proved the following symmetry transformation formula:
$$
I(t_1,\ldots,t_8;t;p,q)=I(s_1,\ldots,s_8;t;p,q),
$$
where
\begin{eqnarray*}
\left\{
\begin{array}{cl}
s_j =\rho^{-1} t_j,&   j=1,2,3,4  \\
s_j = \rho t_j, &    j=5,6,7,8
\end{array}
\right.;
\quad \rho=\sqrt{\frac{t_1t_2t_3t_4}{pqt^{1-n}}}
=\sqrt{\frac{pqt^{1-n}}{t_5t_6t_7t_8}},
\quad |t|,|t_j|,|s_j|<1.
\end{eqnarray*}
It describes a generalization of the key relation for the $E_7$-group \eqref{E7-1},
i.e. integral  \eqref{type2} is invariant with respect to all
transformations of this group.
As shown in \cite{spi:cs}, the function $I(t_1,\ldots,t_8;t;p,q)$
emerges in the quantum multiparticle model of \cite{die:integrability}
under certain restrictions on the parameters (the balancing condition)
as a normalization of a special eigenfunction of the Hamiltonian.

General elliptic hypergeometric integrals of type I on the root
system $A_n$ are defined as
\begin{eqnarray*}
&& I_{n}^{(m)}(s_1,\ldots,s_{n+m+2};t_1,\ldots,t_{n+m+2};A)
\\ && \makebox[2em]{}
= \kappa_{n}^A\int_{\T^n} \prod_{1\leq j<k\leq n+1}
\frac{1}{\eg(z_jz_k^{-1},z_j^{-1}z_k;p,q)}\prod_{j=1}^{n+1}
\prod_{l=1}^{n+m+2}\eg(s_lz_j^{-1},t_lz_j;p,q)\; \frac{dz}{z},
\end{eqnarray*}
where $|t_j|, |s_j|<1$,
$$
\prod_{j=1}^{n+1}z_j=1,\qquad \prod_{l=1}^{n+m+2}s_lt_l=(pq)^{m+1},
$$
and we set $I_{0}^{(m)}=\prod_{l=1}^{m+2}\Gamma(s_l,t_l;p,q)$.
Using the result of exact computation of the $I_{n}^{(0)}$-integral \eqref{AI},
one can derive the following recursive relation for the
 $I_{n}^{(m)}$-integrals in the variable $m$:
\begin{eqnarray*}
&&
I_{n}^{(m+1)}(s_1,\ldots,s_{n+m+3};t_1,\ldots,t_{n+m+3};A)
\\  && \makebox[1em]{}
=\frac{\kappa_{n}^A}{\eg(v^{n+1};p,q)}
\prod_{l=1}^{n+2}\frac{\eg(t_{n+m+3}s_l;p,q)}{\eg(v^{-n-1}t_{n+m+3}s_l;p,q)}
\int_{\T^n} \prod_{1\leq j<k\leq n+1}
\frac{1}{\eg(w_jw_k^{-1},w_j^{-1}w_k;p,q)}
\\  && \makebox[3em]{} \times
\prod_{j=1}^{n+1}
\eg(v^{-n}t_{n+m+3}w_j;p,q)\prod_{l=1}^{n+2} \eg(v^{-1}s_lw_j^{-1};p,q)\;
\\  && \makebox[3em]{} \times
I_{n}^{(m)}(vw_1,\ldots,vw_{n+1},s_{n+3},\ldots,s_{n+m+3};
t_1,\ldots,t_{n+m+2};A)\;\frac{dw_1}{w_1}\cdots\frac{dw_n}{w_n},
\end{eqnarray*}
where
$$
v^{n+1}=\frac{t_{n+m+3}}{pq}\prod_{k=1}^{n+2}s_k
=\frac{(pq)^{m+1}}{\prod_{k=1}^{n+m+2}t_k\prod_{l=n+3}^{n+m+3}s_l}.
$$
For $m=0$ the $I_{n}^{(0)}$-integral on the right-hand side
is computable and  yields the symmetry transformation \cite{spi:theta2}
\begin{eqnarray*}
&&
I_{n}^{(1)}(s_1,\ldots,s_{n+3};t_1,\ldots,t_{n+3};A)
=\prod_{k=1}^{n+2}\eg\left(t_{n+3}s_k,\frac{\prod_{i=1}^{n+2}s_i}{s_k},
s_{n+3}t_k,\frac{\prod_{i=1}^{n+2}t_i}{t_k};p,q\right)
\\  && \makebox[3em]{} \times
I_{n}^{(1)}(v^{-1}s_1,\ldots,v^{-1}s_{n+2},v^{n}s_{n+3};
vt_1,\ldots,vt_{n+2}, v^{-n}t_{n+3};A),
\end{eqnarray*}
where  $\prod_{k=1}^{n+3}t_ks_k=(pq)^2$ and
$$
v^{n+1}=\frac{t_{n+3}}{pq}\prod_{k=1}^{n+2}s_k
=\frac{pq}{ s_{n+3}\prod_{k=1}^{n+2}t_k}.
$$
Since the left-hand side of this
relation is symmetric in the parameters $t_k$ or $s_k$,
the same $S_{n+3}\times S_{n+3}$ symmetry is valid for the right-hand side which
leads to additional nontrivial transformations \cite{rai:trans}.
It is convenient to denote
\begin{eqnarray*}
&& I(v)=\prod_{k=1}^{n+2}\eg(v^{-n-1}s_{n+3}s_k,v^{n+1}t_{n+3}t_k;p,q)
\\ && \makebox[2em]{} \times
I_{n}^{(1)}(vs_1,\ldots,vs_{n+2}, v^{-n}t_{n+3};
v^{-1}t_1,\ldots,v^{-1}t_{n+2},v^ns_{n+3}  ;A),
\end{eqnarray*}
where the arguments of $I_{n}^{(1)}$ lie inside the unit circle,
$\prod_{k=1}^{n+3}t_k=\prod_{k=1}^{n+3}s_k=pq$, and $v$ is an
arbitrary free parameter (the total number of free parameters is
equal to $2n+5$). Then the derived relation can be rewritten
in the form $I(v)=I(v^{-1})$.

Another type of transformations for the
$I_{n}^{(m)}$-integrals was found by Rains \cite{rai:trans}.
We denote $T=\prod_{j=1}^{n+m+2}t_j,\, S=\prod_{j=1}^{n+m+2}s_j$,
so that $ST=(pq)^{m+1}$, and let all $|t_k|,\, |s_k|,\, |T^{\frac{1}{m+1}}/t_k|,\,
|S^{\frac{1}{m+1}}/s_k|<1$. Then the following symmetry transformation is true:
\begin{eqnarray}\nonumber
&&
I_{n}^{(m)}(t_1,\ldots,t_{n+m+2};s_1,\ldots,s_{n+m+2} ;A)
=\prod_{j,k=1}^{n+m+2}\eg(t_js_k;p,q) \\ && \makebox[1em]{} \times
I_{m}^{(n)}\left(\frac{T^{\frac{1}{m+1}}}{t_1},\ldots,
\frac{T^{\frac{1}{m+1}}}{t_{n+m+2}};
\frac{S^{\frac{1}{m+1}}}{s_1},\ldots,\frac{S^{\frac{1}{m+1}}}{s_{n+m+2}};
A\right).
\label{rains2}\end{eqnarray}
This relation generalizes transformation \eqref{E7-2}, and there are no
natural analogues of other $E_7$-reflections.
We see thus that the $V$-function symmetries, generated by different
elements of the Weyl group for the exceptional root system $E_7$,
have some multidimensional analogues. However, different reflections
are generalized to  different integrals whose kernels obey
symmetries in the integration variables related to different
root systems.

\section{Conclusion}

Despite of a rather large volume, many problems were
considered in this review only fragmentarily, a significant number of
statements was given without
proofs, and a number of interesting questions was not touched at all.
Let us list some of the skipped achievements of the theory of elliptic
hypergeometric functions and indicate several important open problems.

Suppose that there exists a finite difference operator of the first order
which maps given rational functions to different  rational functions
with a smaller number of poles (the ``lowering" operator). In the paper
\cite{spi-zhe:grids} it was shown that this is possible only under the
condition that the poles of these rational functions are parameterized
by a general elliptic function of the second order, and the problem itself
is related to the Poncelet mapping.

We skipped description of the connection between biorthogonal rational
functions and the Pad\'e approximation with prescribed zeros and poles
\cite{zhe:gevp,spi-zhe:grids,zhe:pad}.
So, in the paper \cite{zhe:pad}  it is shown that the
$_{12}V_{11}$-elliptic hypergeometric series appears in the Pad\'e
interpolation tables of some functions.

It is natural to expect that the multiple elliptic beta integrals
define a measure in biorthogonality relations for some functions
of many variables, which would generalize the univariate relations
 \eqref{2ib}. The first system of such multivariate functions,
which is based on the elliptic analogue of the Selberg integral
 \eqref{SintB}, was built by Rains \cite{rai:trans,rai:abelian}
with the help of raising and lowering operators.
In certain limits, these functions degenerate to orthogonal polynomials
of Macdonald \cite{mac}, Koornwinder  \cite{koo}, or interpolating
polynomials of Okounkov \cite{oku} (on the connection with the latter polynomials
see also \cite{cos-gus}). The results of the papers \cite{rai:trans,rai:abelian}
represent the most advanced achievements of the theory of elliptic
hypergeometric functions of many variables. Another type of generalization
of the Macdonald polynomials to the level of theta functions was
suggested in \cite{fel-var:mac}.

There is a beautiful geometric interpretation of some of the elliptic
hypergeometric functions in terms of the dynamics on
algebraic surfaces. In \cite{sak}, Sakai gave a classification of
discrete Painlev\'e equations connected with the affine Weyl groups.
On the top of this scheme one has the elliptic Painlev\'e equation
related to the root system $\hat E_8$. In the paper \cite{kmnoy}
this equation was considered in detail and a reduction to the elliptic
hypergeometric equation was found. Respectively, it was indicated
that the elliptic hypergeometric series $_{12}V_{11}$ provides
a particular solution of the elliptic Painlev\'e equation.
A similar role is played by the general solution of the elliptic
hypergeometric equation \cite{spi:thesis,spi:cs} and by the $BC_n$-elliptic
hypergeometric integral of type II for some special values of the
parameters \cite{rai:rec}.

The Pochhammer and Horn approach to the functions of hypergeometric
type \cite{aar,ggr}, which we used in the case of univariate functions,
is not generalized yet to the level of elliptic hypergeometric
functions of many variables. For that it is necessary to learn how
to solve systems of difference equations of the first order
for the kernels of multiple series or integrals with the coefficients
which are elliptic functions of all summation or integration variables
\cite{spi:theta1,spi:theta2}. This task is quite complicated,
and all the examples, which were considered above,
are built on the basis of different constructive ideas.
On this route there appears a general problem of classification of
all types of the elliptic beta integrals and of their multiparameter extensions
to the higher order functions. It is expected, in particular, that
there exist elliptic generalizations of multiple $q$-beta integrals
for the exceptional root systems \cite{ito}.
Let us remark also that under appropriate restrictions
on the parameters all the integrals considered above
are the integrals over polycycles. It would be interesting to
consider integrals over more complicated regions of the complex integration
variables.

In the paper \cite{spi-zhe:spectral}, a nonlinear discrete integrable system
was considered and a self-similar reduction of the corresponding equations
was suggested, which leads to elliptic solutions with many parameters.
Using this result,
a system of discrete biorthogonal functions was built, which are
expressed in terms of the $_{12}V_{11}$-series. Biorthogonal functions
described in section 7.1 represent only a particular subcase of these
more general functions, which
are expressed as linear combinations of several $_{12}V_{11}$-series
and contain three additional parameters. A detailed analysis of these
functions is not performed yet. Let us mention among open problems
related to the general solution of the elliptic hypergeometric
equation a search for an explicit form of the
non-terminating elliptic hypergeometric continued fraction
 and a buildup of elliptic analogues of the associated
Askey--Wilson polynomials \cite{ism-rah:associated}.

Different generalizations of the integral transformation \eqref{f}
to  multidimensional integrals on root systems were suggested in
\cite{spi-war:inversions}, but some of the corresponding inversion
formulae are not proven yet.
The problem of convergency of infinite elliptic hypergeometric series
requires a deep analysis. It is necessary to understand in what sense
such functions can exist. At this moment it is completely unclear
what are the elliptic analogues of the number theoretical properties
of the plain hypergeometric functions. It is necessary also to clarify
whether it is possible to build nontrivial functions of the
hypergeometric type for Riemann surfaces of a higher genus (the
simplest example of such functions is given in  \cite{spi:macd}).

In conclusion, we can state that the main structural elements of the
theory of plain and $q$-hypergeometric functions have their natural
elliptic analogues. Moreover, various ``old" hypergeometric notions
acquire a new meaning connected with the properties of
the elliptic functions. Sufficiently many
applications of elliptic hypergeometric functions in mathematical physics
are known at present:
in the exactly solvable models of statistical mechanics related to
the elliptic solutions of the Yang--Baxter equation  \cite{ft}
and the Sklyanin algebra \cite{rai:abelian,spi:aa2,ros:sklyanin,konno},
in nonlinear integrable discrete time chains \cite{spi-zhe:spectral},
in relativistic quantum multiparticle models of the Calogero--Sutherland
type \cite{spi:cs}, and in nonlinear discrete equations of the Painlev\'e
type \cite{rai:rec,kmnoy}.
It is natural to expect that with time the number of such applications
will grow and, besides that, there will appear new conceptual intersections
with other parts of mathematics.

\medskip
I am deeply indebted to my coauthors J. F. van Diejen,
E. M. Rains, S. O. Warnaar and A. S. Zhedanov for a fruitful
collaboration and many stimulating discussions. During the work on
the theory of elliptic hypergeometric functions discussions of
various problems with H. Rosengren, S. Ruijsenaars and M. Schlosser,
who made their own essential contribution to the development of
this theory, as well as with G. E. Andrews, R. Askey, C. Krattenthaler,
A. Levin, Yu.~I.~Manin, V.~B.~Priezzhev, D. Zagier, and W. Zudilin
 were quite useful. This review is partially based on author's habilitation
thesis \cite{spi:thesis} and lecture notes to an introductory course read
at the Independent University (Moscow) in the fall of 2005.
Some of the presented results were obtained during the visits
to the Max Planck Institute for Mathematics (Bonn), to the directorate
of which I am grateful for the hospitality. This work is supported
in part by the Russian foundation for basic research, grant no. 08-01-00392.

\appendix
\section{Elliptic functions and the Jacobi theta functions}\

The periods of a periodic function are called primitive, if their
linear combinations with integer coefficients yield all periods of
this function. It is well known that nontrivial meromorphic functions
cannot have more than two primitive periods \cite{akh}.
Functions of a real variable can have only one primitive period.
These statements were used in the derivation of the expression
\eqref{mod-e-gamma} and in the construction of elliptic analogues
of the Meijer function. The meromorphic functions $f(u)$ with two
primitive periods are called elliptic functions, that is
there exist $\omega_1,\omega_2\in \C$, Im$(\omega_1/\omega_2)\neq 0$, and
$ f(u+\omega_1)=f(u+\omega_2)=f(u). $

Primitive periods of an elliptic function $f(u)$ form
parallelograms of periods and $f(u)$ is determined by
its values inside of them and on a pair of adjacent edges.
We call as a fundamental domain $D$ the interior of one
of such parallelograms chosen in such a way that on its
boundary $\partial D$ there are no divisor points of $f(u)$.
Clearly,  $f(u)$ (as a meromorphic function) has a finite
number of zeros and poles inside $D$.

The integral $(2\pi i)^{-1}\int_{\partial D}f'(u)/f(u)du$
defines the number of zeros in a domain $D$  for entire functions
and the difference of the numbers of
zeros and poles for meromorphic functions.
Because of the periodicity, this integral is equal to zero
for elliptic functions, that is the number of zeros
equals to the number of poles in $D$.
This number of zeros $s$ (or poles) is called the order of elliptic function.
The equality $f(u)=C$, where $C$ is an arbitrarily chosen constant,
is satisfied in $D$ precisely $s$ times (i.e., $f(u)$ is a
$s$-sheeted function). This statement follows from the fact that
the elliptic function $f(u)-C$ has $s$ poles in $D$ and,
consequently, precisely $s$ zeros.

The sum of residues of the poles of $f(u)$ in $D$ is equal to zero,
which follows from the equality $\int_{\partial D} f(u) du=0$.
Therefore there are no elliptic functions of the order $s=1$,
and an elliptic function of the zeroth order is   constant
(Liouville's theorem). This gives a method of proving elliptic
functions' identities: if the difference of functions $f_1(u)-f_2(u)$
(or the ratio $f_1(u)/f_2(u)$) is elliptic and contains no more than
than one pole in the fundamental domain, then
$f_1(u)-f_2(u)=const$ (or $f_1(u)/f_2(u)=const$).

The well known elliptic function of Weierstrass
 $\wp (u|\omega_1,\omega_2)$ \cite{akh}
has in the fundamental domain one pole of the second order, i.e. $s=2$.
The pair $(x,y)=(\wp (u),\wp '(u))$ defines a uniformization of
an elliptic curve $y^2=4x^3-g_2x-g_3$.
Elliptic functions form a differential field, and any two elliptic functions
 $f(u)$ and $g(u)$ are related by an algebraic relation
$P(f,g)=0$, where $P(f,g)$ is a polynomial of its arguments.
For the choice $g(u)=f'(u)$, one sees that any elliptic function satisfies
some nonlinear differential equation of the first order $P(f,f')=0$.
Taking $g(u)=f(u+y)$, we obtain
$P(f,g)=\sum_{k=0}^Np_{k}(f(u),y)$ $f(u+y)^k=0$, where $p_k(f(u),y)$
are some polynomials in $f(u)$ with the coefficients depending on $y$.
Permuting $u$ and $y$, we see that $p_k(f(y),u)=p_k(f(u),y)$ are
symmetric polynomials of $f(u)$ and  $f(y)$ with constant coefficient.
The condition $P(f,g)=0$ can be rewritten therefore as $Q(f(u),f(y),f(u+y))=0$ for
some polynomial of its arguments $Q$. When such a condition is satisfied,
one says that a function  $f(u)$ obeys an algebraic addition theorem.
As shown by Weierstrass, a meromorphic function $f(u)$ obeying
such an  addition theorem must be either an elliptic function
or its degeneration to a trigonometric or a rational function.

For a more explicit representation of elliptic functions one needs
theta functions. Arbitrary entire functions $f(u)$ are called
(elliptic) theta functions, if
\begin{equation}
f(u+\omega_1)=e^{au+b}f(u),\qquad  f(u+\omega_2)=e^{cu+d}f(u),
\label{def-theta}\end{equation}
for some $a,b,c,d\in\C$ and Im$(\omega_1/\omega_2)\neq 0.$
Replacing $u$ by $\omega_2u$ and multiplying  $f(u)$ by $e^{\alpha u^2 +\beta u}$
with some specially chosen constants $\alpha$ and $\beta$, one can reach
the equalities
$$
f(u+1)=f(u),\qquad f(u+\tau)=e^{au+b}f(u). 
$$
In this appendix we use the parameterization $\tau=\omega_1/\omega_2$
and denote $q=e^{2\pi i\tau}$ (in difference from the main body
of the review, where $\tau=\omega_3/\omega_2$ and $p=e^{2\pi i\tau}$).
For $a\neq 0$, the parameter $b$ can be removed by the shift
$u\to u-b/a$. For a parallelogram $D$ with the vertices $(0,1,1+\tau,\tau)$,
we find $\int_{\partial D}f'(u)/f(u)du= -a$. Therefore
$a=-2\pi is$, where the quantity $s=0,1,\ldots$ determines the
number of zeros of $f(u)$ in $D$. The key characteristics of a
theta functions $s$ is called its order.

The periodicity $f(u+1)=f(u)$ permits to expand $f(u)$ into the Fourier
series $f(u)=\sum_{j=-\infty}^\infty c_je^{2\pi i ju}.$
Substituting it into the second equation with $a=-2\pi is$ and solving
the emerging recurrence relation for the coefficients $c_j$, we find
$$
f(u)=\sum_{l=0}^{s-1}c_lz^l\sum_{k\in\Z}q^{sk(k-1)/2}(q^lz^s)^k,\quad z=e^{2\pi iu}.
$$
The coefficients $c_0,\ldots,c_{s-1}$ are arbitrary, i.e. theta functions
of the order $s$ form an $s$-dimensional vector space.

If we restore arbitrary quasiperiodicity multipliers, then a theta
function without zeros is equal to $e^{P_2(u)}$, where $P_2(u)$
is a polynomial of the second order.  Theta functions of the first order
with one zero in the fundamental parallelogram of quasiperiods are
called the Jacobi theta functions, and the functions with $s>1$
are called theta functions of the higher level. For $s=1$, it is convenient
to work with four theta functions with characteristics
$$
\theta_{ab}(u)=\sum_{k\in \Z}e^{\pi i\tau(k+a/2)^2}
e^{2\pi i(k+a/2)(u+b/2)},
$$
where the variables $a$ and $b$ take the values $0$ or $1$.
The standard Jacobi theta functions are defined as \cite{erd:higher}:
\begin{eqnarray*}
&& \theta_1(u|\tau)=\theta_1(u)=-\theta_{11}(u),
\\ &&
\theta_2(u|\tau)=\theta_2(u)=\theta_{10}(u)=\theta_1(u+1/2),
\\ &&
\theta_3(u|\tau)=\theta_3(u)=\theta_{00}(u)
=e^{\pi i \tau/4+\pi i u}\theta_1(u+1/2+\tau/2),
\\ &&
\theta_4(u|\tau)=\theta_4(u)=\theta_{01}(u)
=-ie^{\pi i \tau/4+\pi i u}\theta_1(u+\tau/2).
\end{eqnarray*}
Note that all of them have the form $\sum_{k\in\Z}c_k$ with
 $h(k)=c_{k+1}/c_k=q^ky$ for some constant $y$.
Since $h(k)$ is rational in $q^k$, Jacobi theta functions
represent a special class of $q$-hypergeometric functions.

The $\theta_1(u)$-function is odd, $\theta_1(-u)=-\theta_1(u),$
and satisfies the quasiperiodicity conditions
$$
\theta_1(u+1)=-\theta_1(u),\qquad
\theta_1(u+\tau)=-e^{-\pi i\tau -2\pi iu}\theta_1(u).
$$
It is related to the shortened theta function
$\theta(z;q)=(z;q)_\infty(qz^{-1};q)_\infty$
by the Jacobi triple product identity
\begin{equation}
\theta_1(u)=iq^{1/8}e^{-\pi iu}(q;q)_\infty\theta(e^{2\pi iu};q).
\label{tji}\end{equation}
Transformation properties of the $\theta_1$-function with respect to
the $PSL(2,\Z)$-group of modular transformations
$\tau\to (a\tau+b)/(c\tau+d)$, $a,b,c,d\in\Z,\, ad-bc=1$,
are determined by the relations \cite{eic-zag:theory}
\begin{eqnarray}
&& \theta_1(u|\tau+1)=e^{\pi i/4}\theta_1(u|\tau), \quad
\theta_1\left(\frac{u}{\tau}\Big|\frac{-1}{\tau}\right)
=-i\sqrt{-i\tau}\; e^{\pi iu^2/\tau}\theta_1(u|\tau).
\label{mod}\end{eqnarray}

It is convenient to use notation $\theta_a(u_1,\ldots,u_k):=\theta_a(u_1)\cdots
\theta_a(u_k)$ and $\theta_a(x\pm y):=\theta_a(x+y,x-y)$.
Then the argument duplication formula has the form
$$
\theta_1(2u)=\frac{iq^{1/8}}{(q;q)_\infty^3}
\theta_1\left(u,u+\frac{1}{2},u+\frac{\tau}{2},u-\frac{1+\tau}{2}\right).
$$
The addition theorem for theta functions, which is called sometimes a
Riemann relation, uses products of four theta functions
\begin{equation}
\theta_1(u\pm a,v\pm b)-\theta_1(u\pm b,v\pm a)=\theta_1(a\pm b,u\pm v)
\label{add-add}\end{equation}
or
\begin{eqnarray}
\theta(xw^{\pm 1},yz^{\pm 1};p) -\theta(xz^{\pm 1},yw^{\pm 1};p)
=yw^{-1}\theta(xy^{\pm 1},wz^{\pm 1};p).
\label{ident}\end{eqnarray}
The proof of this equality if elementary. The ratio of expressions
standing on its left- and right-hand sides is a bounded function of
$x\in\C^*$ (it is invariant with respect to the transformation
$x\to px$, and it does not contain poles in the annulus
$|p|\leq |x|\leq 1$). By the Liouville theorem this ratio
does not depend on $x$, but for $x=w$ it equals to 1.

Any theta function $f(u)$ of the order
$s$ with the quasiperiods $\omega_1,\omega_2$
and coordinates of the zeros in the fundamental domain
$a_1,\ldots, a_s$ can be represented in the form
\begin{equation}
f(u)=e^{P_2(u)} \prod_{k=1}^s\theta_1\left(\frac{u-a_k}{\omega_2}\Big|
\frac{\omega_1}{\omega_2}\right),
\label{theta-fact}\end{equation}
where $P_2(u)$ is some polynomial of $u$ of the second order.
Indeed, the function
$$
g(u)=\prod_{k=1}^s\theta_1\left(\frac{u-a_k}{\omega_2}\Big|
\frac{\omega_1}{\omega_2}\right)
$$
is a theta function of the order $s$ with the same zeros
in the parallelogram of quasiperiods $\omega_{1}$ and $\omega_{2}$
as the function $f(u)$. Therefore the ratio
$f(u)/g(u)$ is an entire function without zeros and poles,
namely, a theta function of the zeros order, i.e. $e^{P_2(u)}$.
In \cite{eic-zag:theory}, a self-contained theory of Jacobi forms ---
the functions obeying transformation properties similar to those of
functions  \eqref{theta-fact}, was formulated.
Note that all vectors  $f_j(u),\, j=1,\ldots,s,$ of any
basis of the space of theta functions of the order $s$ can be represented in the
indicated form with the matrix of zeros $a_{jk}$ satisfying
a number of constraints (e.g., $\sum_{k=1}^sa_{jk}=const$).

We call as meromorphic theta functions ratios of theta functions
of an arbitrary finite order. It is easy to see that they define
meromorphic solutions of  equations  \eqref{def-theta}. For this we
denote as $a_1,\ldots,a_n$ coordinates of the zeros and as $b_1,\ldots,b_m$
coordinates of the poles of the corresponding function $f(u)$ in the
fundamental domain. Then the ratio
$$
\frac{f(u)\prod_{k=1}^m\theta_1\left(\frac{u-b_k}{\omega_2}\big|
\frac{\omega_1}{\omega_2}\right)}
{\prod_{k=1}^n\theta_1\left(\frac{u-a_k}{\omega_2}\big|
\frac{\omega_1}{\omega_2}\right)}
$$
is an entire function without zeros satisfying the equations
$f(u+\omega_1)=e^{a'u+b'}f(u)$ and  $f(u+\omega_2)=e^{c'u +d'}$
with some $a',b',c',d'$, that is a theta function of the
zeroth order $e^{P_2(u)}$.

Any elliptic function $f(u)$ of the finite order  $s$
with the periods $\omega_1,\omega_2$ can be represented in the form
\begin{equation}
f(u)=C\; \prod_{k=1}^s\frac{\theta_1\left(\frac{u-a_k}{\omega_2}\big|
\frac{\omega_1}{\omega_2}\right)}{\theta_1\left(\frac{u-b_k}{\omega_2}\big|
\frac{\omega_1}{\omega_2}\right)},
\label{ell-fact}\end{equation}
where $C$ is some constant, and $a_k$ and $b_k$ denote coordinates of
some zeros and poles  of $f(u)$ congruent to the zeros and poles
in the fundamental parallelogram of periods. The following constraint
should be satisfied by $a_k$ and $b_k$:
\begin{equation}
a_1+\ldots+a_s=b_1+\ldots+b_s \mod \omega_2.
\label{balance-add}\end{equation}
This follows from the fact that both parts of equality \eqref{ell-fact} are
meromorphic and doubly periodic. Therefore their ratio
defines a bounded entire function, i.e. a constant.
The linear constraint on the values of
$a_k$ and $b_k$ \eqref{balance-add}, which we call
the balancing condition, follows from the requirement of
cancellation of the quasiperiodicity multipliers
of the $\theta_1$-functions appearing from the
$u\to u+\omega_1$ shift.

Substituting expression \eqref{tji} in \eqref{ell-fact} and denoting
$z=e^{2\pi iu/\omega_2}$, $t_k=e^{-2\pi ia_k/\omega_2}$,
$w_k=e^{-2\pi ib_k/\omega_2}$, we obtain
\begin{equation}
f(u)=\pm C\; \prod_{k=1}^s\frac{\theta(t_kz;q)}{\theta(w_kz;q)},
\qquad \prod_{k=1}^st_k=\prod_{k=1}^sw_k,
\label{ell-fact-m}\end{equation}
where the sign ambiguity appears from
the factor $e^{-\pi iu}$ in \eqref{tji}.

\end{document}